\documentclass{siamltex}

\usepackage{amssymb,amsmath,amscd,verbatim}
\usepackage{algorithmic}
\usepackage{bm}
\usepackage{bbm}
\usepackage{tikz}
\usetikzlibrary{shapes}

\usepackage{algorithmic}
\usepackage{algorithm}
\usepackage{subfigure}
\newcommand{\bi}{\begin{itemize}}
\newcommand{\ei}{\end{itemize}}

\newcommand{\caliber}{{c}}

\newcommand{\finevar}{F}
\newcommand{\coarsevar}{C}

\usepackage{color}

\numberwithin{equation}{section}
\numberwithin{table}{section}
\numberwithin{figure}{section}

\title{Algebraic distance for anisotropic diffusion problems: multilevel results}
\author{A.~Brandt\footnotemark[2], J.~Brannick\footnotemark[3], K.~Kahl\footnotemark[4], and I.~Livshits\footnotemark[5]}

\begin{document}
\maketitle

\begin{abstract}
{
In this paper we motivate, discuss the implementation and present the resulting numerics for a new definition of strength of connection which is based on the notion of {\em algebraic 
distance}. This algebraic distance measure, combined with compatible relaxation, is used to choose suitable coarse grids and accurate interpolation operators for algebraic multigrid algorithms. The main tool of the proposed measure is the least squares functional defined using a set of relaxed test vectors. The motivating application is the anisotropic diffusion problem, in particular problems with non-grid aligned anisotropy.  
We demonstrate numerically that the measure yields a robust technique for determining strength of connectivity among variables, for both two-grid and multigrid solvers. 
The proposed algebraic distance measure can also be used in any other coarsening procedure, assuming a rich enough set of the near-kernel components of the matrix for the targeted system is known or computed.   
} 
\end{abstract}


\footnotetext[3]{Department of Mathematics,
Pennsylvania State University, University Park, PA 16802, USA (brannick@psu.edu). \thanks{Brannick's work was supported by the National Science Foundation under grants OCI-0749202 and DMS-810982.}}
\footnotetext[4]{Fachbereich Mathematik und Naturwissenschaften,
Bergische Universit\"at Wuppertal, D-42097 Wuppertal, Germany, 
(kkahl@math.uni-wuppertal.de). \thanks{Kahl's work was supported by the
  Deutsche Forschungsgemeinschaft through the Collaborative Research
  Centre SFB-TR 55 ``Hadron physics from Lattice QCD''}}
\footnotetext[5]{Department of Mathematical Sciences, Ball State University, Muncie IN, 47306, USA (ilivshits@bsu.edu), \thanks{Livshits' work was supported by DOE subcontract B591416.}}
\renewcommand{\thefootnote}{\arabic{footnote}}

\maketitle

\section{Introduction}

Multigrid (MG) is a methodology for designing numerical iterative methods for solving sparse matrix equations.  One of the key issues in developing an efficient multigrid method is 
the selection of the coarse grids, which must consist of a smaller number of degrees of freedom than the fine grid, yet still be rich enough to allow 
for the accurate representation of smooth error.  In the Algebraic Multigrid (AMG) framework considered here, the coarse grids, the fine-to-coarse residual transfer (restriction) operators, the coarse-to-fine correction transfer (prolongation) operators, and the coarse grid operators are all computed in a two-level setup
algorithm that proceeds recursively.  A main tool used in the AMG setup algorithm is strength of connectivity between variables.  
For many problems, e.g., those arising from discretization of a partial differential equation (PDE), the connectivity among variables and thereby a proper choice of these multigrid components are known and 
well understood theoretically.  This is, however, not the case for the anisotropic diffusion problems that are the focus of this paper.

Generally speaking, multigrid methods for solving systems of sparse linear equations  
$Ax=b$ are all based on the smoothing property of relaxation.
An error vector $e$ is called $\tau$-smooth if
its 
residual
is smaller than $\tau
\|e\|$. The basic observation \cite{B83} is that the convergence
of a proper relaxation process
 slows down only when the current error is
$\tau$-smooth with $\tau \ll 1$, the smaller the $\tau$ the slower
the convergence.  In other words, the convergence rate of relaxation deteriorates 
as the error becomes dominated by eigenvectors with small eigenvalues in magnitude. 
This observation and the assumption that when
relaxation slows down, the error vector $e$ can be approximated in
a much lower-dimensional subspace, are the main ideas behind the multigrid 
methodology. 
{
Very efficient {geometric
multigrid} solvers have been developed for the case that this
lower-dimensional subspace { of smooth errors} corresponds to functions defined 
on a well-structured grid (the {\it coarse} level) and can be  approximated by easy-to-derive 
equations, based for example on discretizing the same continuum 
operator that has given rise to the fine-level equations $Ax=b$.} The coarse-level
equations are solved  recursively using a similar combination of
relaxation and still-coarser-level approximations to the
resulting smooth error.

To deal with more general problems, e.g., the ones where the
fine-level system may not be defined on a well-structured grid
nor perhaps arise from any continuum problem, {algebraic
multigrid}  methods have been developed.  These methods
require techniques for deriving the
set of coarse-level {\it variables} and the coarse-level {\it
equations}, based solely on the (fine-level) matrix $A$. The basic
approach, developed in \cite{BMR83, oAMG, JWRuge_KStuben_1987a} and called today
classical AMG or RS-AMG,  involves the following two steps:
\begin{itemize}
\item[(1)] Choosing the coarse-level variables as a subset $\coarsevar$ of
the set $\Omega$ of fine variables, such that each variable in $\Omega$ is
{\it strongly connected} to variables in $\coarsevar$.
 \item[(2)]
Approximating the fine-level residual problem $Ae=r$ by the
coarse-level equations $A_c e_c=r_c$ using the Galerkin
prescription $A_c = RAP$ and $r_c=Rr$,  yielding an approximation $ P e_c$ to $e$.
\end{itemize}
The {\it interpolation matrix} $P$ and the {\it restriction matrix} $R$
are both defined directly in terms of the entries of the matrix
$A$.  Their construction relies on the notion 
of  {\it strong connections}, developed originally 
for $M$-matrices, to provide a measure
of the coupling of the variables used explicitly
to coarsen the problem.  

In the past two decades, numerous extensions of the classical AMG algorithm have
been introduced, including modifications to the coarse-grid selection
algorithms
\cite{AJCleary_RDFalgout_VEHenson_JEJones_1998a,HDeSterck_UMYang_JJHeys_2005a,VEHenson_UMYang_2002a}, 
and the definition of interpolation \cite{Brannick_Trace_06,MBrezina_etal_2000a,MBrezina_etal_2005a,TChartier_etal_2003a,PVanek_JMandel_MBrezina_1995a}.
These works are motivated by the fact that
 the applicability of the original AMG algorithm
is limited by the $M$-matrix heuristics upon which its 
strength of connectivity measure among variables is based.
In particular, the measure and, hence, the classical AMG approach yields an efficient solver for problems where the matrix $A$ has a dominant diagonal and, with small possible exceptions, all its off-diagonal elements have the same sign.  Even then, the produced interpolation can have limited accuracy, insufficient for full multigrid efficiency.  An example where the
performance of AMG methods can deteriorate, and the one we consider here, is given by non-grid aligned scalar anisotropic diffusion problems.   {\color{black}We refer to the paper~\cite{Olson_10} in which a smoothed aggregation variant of AMG is applied to this same test problem we consider and its performance is shown to deteriorate for non-grid aligned constant coefficient cases.  Table VI gives results for two of the tests we consider in the numerical results section.}

Advances in the design of AMG methods for anisotropic problems include the development of new notions of strength of connection that are used to  choose coarse variables and the sparsity pattern and coefficients of interpolation~\cite{Olson_10,Bran_EBSOC,brannick_amli_2011,brannick_local_stab_2011}.
In~\cite{Bran_EBSOC},
an energy-based measure of strength of connection was developed and tested for anisotropic diffusion problems.  The approach
uses a measure of strength of connectivity among variables 
based on approximations of the columns of the inverse of the system 
matrix that are computed using local relaxation.
\textcolor{black}{In~\cite{Olson_10}, a related approach based on an evolution measure 
is developed and then applied to similar anisotropic model problems as we consider here.}   The basic approach considers the connection between weighted-Jacobi-relaxation and the time integration of ordinary differential equations for the specific case of evolving $\delta$-functions to form
the proposed strength measure.  The works in~\cite{brannick_amli_2011,brannick_local_stab_2011} develop a strength measure based 
on local energy estimates for interpolation to form aggregates for non-grid aligned anistropic problems.

Closely related to these works is the  approach of compatible relaxation
\cite{ABrandt_2000a}, which 
uses a modified relaxation scheme to expose the character of the
slow-to-converge, i.e., algebraically smooth, error.  Coarse-grid points are then selected where
this error is the largest, thus avoiding explicit use of a measure of strength of connection
in choosing the set of coarse variables.  A theoretical framework for CR-based coarsening is 
introduced in~\cite{PanayotRob_2003}, and extensive tests for anisotropic diffusion problems on structured and unstructured meshes are found in~\cite{Brannick_Trace_06,OLivne_2004a,JBrannick_RFalgout,JBrannick_2005a}. 
A related smoothed aggregation-based aggressive coarsening algorithm that uses 
the evolution measure from~\cite{Olson_10}
is integrated with an energy minimizing form of interpolation in~\cite{Jacob},
yielding an effective solver for two-dimensional non-grid aligned anisotropic problems.
In this work, the CR serves as an analysis tool to develop the aggregation
scheme. 

While these developments have resulted in marked improvements
in certain cases, generally speaking all existing methods require a substantial overlap of the coarse grid basis functions (columns of $P$) to obtain fast multigrid convergence for general 
non-grid aligned anisotropic problems.
In the BAMG context, this amounts to using a high-caliber interpolation (i.e, with large interpolatory sets), which leads to a rapid fill-in of the coarse-grid operators.  

In the present article, we focus on developing an alternative compatible relaxation coarsening algorithm that uses an {\em algebraic distance} measure to select coarse variables and interpolatory sets.     
Our aim is to investigate the suitability of such an approach together with a low-caliber (i.e, one with small interpolatory sets) interpolation constructed using the bootstrap algebraic multigrid (BAMG) framework for anisotropic problems.  

The algebraic distance measure, we propose, is based on a notion of strength of connectivity among variables  that is derived from  the local least squares (LS) formulation for computing  caliber-one interpolations~\cite{B00,BAMG2010,DB_11}. The approach first constructs a caliber-one LS interpolation for a given set of test vectors (representatives of algebraically smooth error) and then defines the algebraic distance between a fine point and its neighboring points  
in terms of the values of the local least squares functionals resulting from the so defined interpolation.  The algebraic distance measure thus aims to address the issue of strength of connections in a general context -- the goal being to determine explicitly those degrees of freedom from which a high quality least squares interpolation for some given set of test vectors can be constructed.  

The resulting measure of distance (connectivity) among variables is used to derive a strength graph which is then passed to a coloring algorithm~\cite[see Chap. 8]{WLBriggs_VEHenson_SFMcCormick_2000a} 
to coarsen the unknowns at each stage of a compatible relaxation 
coarsening algorithm.  
Given the set of coarse variables, a similar measure is used in defining interpolation.  The idea is to approximately minimize the values of the LS functional locally~\cite{B00,BAMG2010}.   
This is accomplished by forming the LS-based interpolation for a given fine point and various candidate  interpolatory sets (neighboring pairs (or sets) of coarse points) and then by choosing the set with smallest value of the associated LS functional to define the row of interpolation.  

The remainder of the paper is organized as follows. 
 Section \ref{sec:bamg} contains an introduction to the bootstrap algebraic multigrid components, with an emphasis on the least squares interpolation approach and compatible relaxation coarsening algorithm.   Then, in Section \ref{sec:distances}, a general definition of strength of connection and the notion of algebraic distance,  as well as its connection to compatible relaxation, are discussed.   
The diffusion equation with anisotropic coefficients and its discretizations are introduced in Section \ref{sec:numerics} as are  the results of numerical experiments of our approach applied to these systems. Conclusions are presented in Section \ref{sec:conclusions}.

\section{Bootstrap AMG}
\label{sec:bamg}
The bootstrap AMG framework, introduced in \cite[\S
17.2]{B00}, was proposed to extend AMG to general (non $M$-matrix)
problems.
The framework combines the following two general devices to
inexpensively construct high-quality interpolation. \bi
\item[(A)] {The interpolation is derived to provide the best  least squares
fit  to a set of $\tau$-smooth 
test vectors (TVs) } obtained by a process described
below. 
\ei Denote by $C_i$ the set of coarse variables used in
interpolating to fine grid variable  $i$. 
It follows from the
satisfaction of the compatible relaxation criterion (see the next section) that, with proper choices of
$C_i$ for all $i$, there exists an interpolation operator that approximates {\it any} vector $x$ which
is $\tau$-smooth with error proportional to $\tau$, i.e, there exists an interpolation operator that satisfies the 
so-called weak approximation property.  The proof of this result for the special case of the so-called ideal interpolation
operator is given in \cite{PanayotRob_2003}. 
The size $|C_i|$ of this set should in principle
increase as $\tau$ decreases (and smaller $\tau$ means overall better multigrid convergence), but in practice a pre-chosen and sufficiently small interpolation {\it caliber} \[ \displaystyle{\caliber := \max_{i \in \Omega}|\coarsevar_i|}\]
often yields small enough
errors, as demonstrated by the extensive numerical tests presented in~\cite{JBrannick_RFalgout}.

The set $C_i$ can often be adequately chosen by natural
considerations, such as the set of geometrical neighbors with
$i$ in their convex hull. If the chosen set is inadequate, the
least squares procedure will show a {\it poor fitness} (interpolation
errors large compared with $\tau$), and the set can then be extended. The least squares procedure can also detect variables in
$C_i$ that can be discarded without a significant accuracy loss.
Thus, this approach allows creating interpolation with whatever
needed accuracy which is {\it as sparse as possible}. 
\bi \item[(B)] {Generally, the test vectors 
are constructed in a bootstrap manner}, in which several tentative
AMG levels are generated by interpolation fitted to only moderately smooth TVs; this tentative (multilevel)  structure  is then  used to produce {\it improved} (much
smoother) TVs, yielding a more accurate interpolation operator.  The process continues  if needed until fully efficient AMG
levels have been generated. \ei

\noindent
 {\color{black}The first test vectors are each produced by relaxing the homogeneous system
\begin{equation} Av=0,
\label{eq:hom}
\end{equation}  using  different initial vectors,  leading to an initial 
$\tau$-smooth test set.  A mixture of random
vectors and/or diverse geometrically smooth vectors can generally be
used as initial approximations. The latter test vectors may not require relaxation at all or, perhaps, only near boundaries or other regions where singularities or discontinuities are present.  In the case of discretized isotropic
PDEs, if geometrically smooth vectors that satisfy the homogeneous
boundary conditions are used, relaxation may not be needed at all.
For the anisotropic problem considered, we use a constant vector and bootstrapped test vectors generated
from initially random initial guesses 
to define an algebraic distance notion of strength of connection and the least squares interpolation operator.}

\subsection{Least squares interpolation}\label{sec:lsp}

The basic idea of the least squares interpolation approach is to approximate
a set of test vectors, $\{v^{(1)}, \ldots, v^{(k)}\}
\subset \mathbb{R}^n$, 
minimizing the interpolation error for these
vectors in a least squares sense.  In the context considered here,
namely, applying the least squares process to construct a classical AMG form of interpolation,  \textcolor{black}{each row of $P$, denoted by $p_i$,}
is defined as the minimizer of a local least squares functional. Given
a splitting of variables $ \coarsevar$ and $\finevar = \Omega \setminus \coarsevar$, 
for each $i \in \finevar$ and for various $\coarsevar_i \subset C$, find $p_{i} \in \mathbb{R}^{n_{c}}, n_c = 
|\coarsevar|$, 
that minimizes
\begin{equation}\label{eq:LSfuncrowi}
\mathcal{LS}(p_i) = 
\sum_{\kappa=1}^k\omega_\kappa\left(v_{i}^{(\kappa)} - \sum_{j\in \coarsevar_i} \left(p_{i}\right)_{j} v_j^{(\kappa)}\right)^{2} \mapsto \min.
\end{equation} 
To guarantee the uniqueness to the of solution (\ref{eq:LSfuncrowi}), the number of test vectors, $k$ should be greater or equal to the caliber $c$. Further, in practice we have observed that the accuracy of the least
squares interpolation operator and, hence, the performance of the resulting solver generally improves with increasing $k$~\cite{BAMG2010}, up to some value proportional to caliber $\caliber$. 
The weights $\omega_{\kappa}>0$ are chosen to reflect the global algebraic smoothness
of the test vectors.  We give their specific choice in the numerical experiments section.  

The quality of the LS  interpolation can be further improved by applying a sweep of local Jacobi relaxation to the test vectors  prior to computing the LS-functional. Equivalently, the action of such pre-smoothing can be directly built-in by considering a residual-based LS process, see~\cite{iBAMG}:  
\begin{equation}\label{eq:modls1}
  \mathcal{LS}(p_i) = 
\sum_{\kappa = 1}^k\omega_{\kappa} \left(v_{i}^{(\kappa)} + \frac{1}{a_{ii}}r^{(\kappa)}_{i} - \sum_{j\in \coarsevar_{i}} (p_{i})_j v_{j}^{(\kappa)}\right)^{2} \mapsto \min, 
\end{equation}
where $r^{(\kappa)}$ is the residual of $v^{(\kappa)}$ in (\ref{eq:hom}).

\subsection{Compatible relaxation}\label{sec:CR}
A general criterion for choosing an adequate set of coarse
variables is the fast convergence of  {\sl compatible relaxation}
(CR), as introduced in \cite{ABrandt_2000a} and further developed in  \cite{OLivne_2004a,PanayotRob_2003,Brannick_Trace_06,JBrannick_RFalgout}.	

Compatible relaxation can be generally employed in two capacities. First,  given a relaxation scheme and a coarse grid  set $\coarsevar$, it can serve  to predict multigrid convergence (habituated CR in \cite{OLivne_2004a} gives the best estimates). Second,  given a relaxation scheme and a desired convergence rate, it can be used to construct an adequate coarse grid. A brief introduction to the details of CR relevant to the content of this paper is presented next.

Given a matrix $A\in \mathbb{R}^{n\times n}$ and a suitable relaxation process with error propagation matrix $E = I - M^{-1}A$, assume that a classical AMG coarse-and-fine level splitting $\Omega = \coarsevar \cup \finevar$ has been selected. 
Then, one choice of compatible relaxation is given by $\finevar$-relaxation for the homogeneous system, i.e, relaxation applied only to the set 
of $\finevar$ variables.
In other words, ordering the unknowns according to the partitioning of $\Omega$ into $\finevar$ and $\coarsevar$:
$$
u = \begin{pmatrix} u_f \\ u_c \end{pmatrix}, \enspace 
A = 
\begin{pmatrix}
    A_{ff}  & A_{fc}   \\
    A_{cf} & A_{cc} 
\end{pmatrix},
\enspace \mbox{and} \enspace
M = 
\begin{pmatrix}
    M_{ff}  & M_{fc}   \\
    M_{cf} & M_{cc} 
\end{pmatrix},
$$
the $\finevar$-relaxation form of compatible relaxation is then defined by 
\begin{equation} \label{frelax:eq}
   u_f^{\nu+1} = E_{f}u_{f}^{\nu} \quad \mbox{with} \quad E_{f} = I-M_{ff}^{-1}A_{ff}.
\end{equation}
If $A$ is symmetric and positive definite and $M$ is symmetric, then the asymptotic convergence rate of compatible relaxation,   
\[
 \rho_{f}
=\rho(E_{f}), 
\]
where $\rho$ denotes the spectral radius, provides a measure of the quality of the
coarse space, that is, a measure of the ability of the set of coarse
variables to represent error not eliminated by the given fine-level relaxation.
This measure can be approximated using $\finevar$-relaxation for the homogeneous 
system with a random initial guess $u_f^{0}$. Since $\lim_{\nu \to \infty} 
\|E_f^\nu\|^{1/\nu} = \rho(E_f)$ for any norm $\| \cdot \|$, the measure
\begin{equation} \label{meas1:eq}
\varrho_f  = \left(\|u_f^\nu\| / \|u_f^0\| \right)^{1/\nu}
\end{equation}
estimates $\rho_f$ and provides a measure of the quality of the coarse variable set. \\

As a tool for choosing $\coarsevar$, CR can be used in the following way. Starting with an initial set of coarse variables $C_0$,
 a few sweeps of compatible relaxation will detect slowness if the current set is inadequate and expose viable candidates to be added to $\coarsevar$, those that are slow to converge to zero in the CR process. 
A subset of these variables, is then added to $C_0$ to form a new coarse set, $C_1$.  The process  repeats until fast convergence of CR is obtained.  In practice,  $C_0$ may be empty or for structured problems standard coarsening can be used.  
The process continues until the desired $\rho_f$ is achieved. \\

Current CR algorithms do not use  
a strength of connection measure in constructing the coarse variable set. Instead, they rely on the error
produced by CR to form candidate sets of potential $\coarsevar$-points.  Our aim  is to develop
a more general notion of strength of connections based on algebraic distances 
and to explore its use in a compatible relaxation coarsening scheme
and in defining the interpolatory set for each row of interpolation.
  
\section{Selecting coarse variables and interpolation via algebraic distances}
\label{sec:distances}
The classical definition of strength of connection is intended for the case
of diagonally dominant M-matrices; it can break down when applied to problems involving  
more general classes of matrices, such as the anisotropic problems considered in this paper.
The near null space of a diagonally dominant M-matrix is typically well approximated locally
by the constant vector, and the AMG strength of connections measure succeeds when this assumption is 
reflected in the coefficients of the system matrix, $A$. (Similar observations motivated 
the work in \cite{Olson_10}.)  If either the near
null space cannot be accurately characterized as locally constant or
this is not reflected in the matrix coefficients, then performance of classical 
AMG  suffers.

 As a more general measure of strength of connection, we consider a variable's ability to
interpolate $\tau$-smooth error for small $\tau$ to its neighbors.
Specifically, for a given fine-grid point $i \in F$, we construct a row of LS interpolation
from each of its neighbors $j \in N_i$, $N_i$ being defined in terms of the graph of powers of $A$, and monitor the values of the LS functional.  If for some $j\in N_i$, 
the LS functional is small relative to its size for other neighbors, then $j$ is determined to be strongly connected to $i$. 
This process, repeated for each $i \in F$,  allows  to identify
suitable coarse-grid points as those from which it is possible to build a high-quality LS interpolation
to its fine-grid neighbors.  Next, we describe in detail the idea of  measuring strength between neighbors using algebraic distance and then discuss how this measure can be incorporated in  CR coarsening algorithm and in computing the nonzero sparsity pattern of interpolation.

\subsection{Strength of connection by algebraic distances}\label{sec:alg_dist}
In the simplest form, the definition of algebraic distances is
straightforward. For any given pair of fine variables $i,j \in \Omega$ compute
\begin{equation}
 \mu_{ij}  = \frac{1}{
\sum_{\kappa = 1}^k\omega_{\kappa} \left(v_{i}^{(\kappa)} + \frac{1}{a_{ii}}r^{(\kappa)}_{i} -  p_{ij} v_{j}^{(\kappa)}\right)^{2}} ,
\label{eq:measure}
\end{equation}
where $p_{ij}$ is the minimizer of the least squares problem for  a
single variable $j$. 
We note that there are many possible choices for defining a measure of strength of connection based on $\mathcal{LS}$. We choose to  take the reciprocal value of $\mathcal{LS}$ to define strength since this coincides with the usual notion of strength between unknowns in that two unknowns are classified as being strongly coupled by a large value of the given measure. Further, our  {experience} suggests that this measure is robust for the targeted anisotropic problems.   We note that the definition of $\mu_{ij}$ is not symmetric since $\mu_{ji}$ involves different entries of the test vectors used in computing LS interpolation than does $\mu_{ij}$. 

For a given SPD matrix $A$, we define its connectivity graph as $\mathcal{G}=(\mathcal{V},\mathcal{E})$,
where $\mathcal{V}$ and $\mathcal{E}$ are the sets of vertices and
edges.
Here, an edge $(i,j) \in E \iff (A)_{ij} \neq0$.   
Let $\mathcal{G}_d(\mathcal{V}_d,\mathcal{E}_d)$ denote the graph of the matrix $A^d$ and
define $\mathcal{G}_{d, i}(\mathcal{V}_{d,i},\mathcal{E}_{d,i})$ as  the subgraph associated with the $i$th vertex and its algebraic neighbors:
\begin{equation}\label{graph}
\mathcal{V}_{d,i} := \{\:  j \:  | \: (A^d)_{ij} \neq 0\} \quad \mbox{and} \quad \mathcal{E}_{d,i} := \{ (i,j) |  \: (A^d)_{ij} \neq 0\}.
\end{equation}  
Then, given a search depth $d$ and a fine variable $i \in \Omega$, we compute $\mu_{ij}$ for all $j \in \mathcal{V}_{d,i}$.  This simplification, combined  with the idea of deriving strength according to an algebraic distance measure based on simple caliber-one interpolation, allows us to control the complexity of the algorithm.
More generally, the algebraic distance measure can be computed for sets of neighboring coarse points and, again, the LS functional can be used as an a posteriori measure of the quality of the interpolatory set $\coarsevar_i$.  We use these observations in the design of our algorithm
for computing the interpolation coefficients discussed in Section \ref{sec:greedy}.  

\subsection{Compatible relaxation coarsening and algebraic distances}\label{sec:cr}
In selecting $\coarsevar$, we integrate the simplified variant of the algebraic distance notion of strength of connection based on caliber-one interpolation~\eqref{eq:measure} into the CR-based coarsening 
algorithm developed in~\cite{JBrannick_RFalgout}.   The notion of algebraic distances is used to form a subgraph of the graph of  
the matrix $A^d$, $d=1,2,\dots$
Specifically, the algebraic distance between any two adjacent vertices 
in the graph $\mathcal{G}_d$
of $A^d$ is computed using~\eqref{eq:measure}.  Then, for each vertex, $i \in \finevar$, we remove edges
adjacent to $i$ with small weights relative to the largest weight of all edges connected to $i$:
\begin{equation}\label{eq:strength_graph}
\mathcal{V}_{\mathcal{M}} = \finevar \; ; \quad 
\mathcal{E}_{\mathcal{M}} = \{ (i,j) \:  | \: i,j \in \finevar  \quad \mbox{and} \quad \mu_{ij} > \theta_{ad}  \max_k \mu_{ik} \}, 
\end{equation}
with $ \theta_{ad} \in (0,1)$.
This, in turn,
produces the graph of strongly connected vertices, $\mathcal{M}^d(\mathcal{V}_{\mathcal{M}}, \mathcal{E}_{\mathcal{M}}) $.  
Note that by definition, the strength graph is restricted to vertices $i\in F$ so that it can be used in the subsequent CR coarsening stages. (Although the measure (\ref{eq:measure}) is able to determine the coupling between any given two points, in order to make the idea practical we restrict its 
use to local neighborhoods; in cases where $\mu_{ij}=\infty$, variables $i$ and $j$ are defined as strong neighbors.  In such cases, 
we do not consider these links in the later stages of the algorithm, that is, they are not used in computing the 
maximum in  \eqref{eq:strength_graph}.)
The strength graph, $\mathcal{M}^d$, is then passed to a coloring
algorithm~\cite[Chap. 8]{WLBriggs_VEHenson_SFMcCormick_2000a}, as in the classical AMG approach, in which coarse points are selected based on their number of strongly connected
neighbors. 

A general description of the overall CR coarsening approach is given by Algorithm 1. For additional specific details of the CR algorithm we refer the reader to~\cite{JBrannick_RFalgout}, in which this scheme was developed for diffusion problems similar to those we consider here.

\begin{algorithm}
\caption{compatible\_relaxation \hfill \{Computes $\coarsevar$ using Compatible Relaxation\}\label{alg:CR}}
\begin{algorithmic}
\STATE \textit{Input:}\,  $u^0$, \,  $\coarsevar_{0}$ \,  \COMMENT{$\coarsevar_{0} = \emptyset$ is allowed}.
\STATE \textit{Output:} $\coarsevar$
\STATE Initialize $\coarsevar = \coarsevar_{0}$
\STATE Initialize ${\finevar} = \Omega \setminus \coarsevar$
\STATE Perform $\nu$ CR iterations starting with an initial guess $u^0$
\WHILE{$\rho_{f} > \delta$}
\STATE $\coarsevar = \coarsevar \cup \{\text{\: maximum independent set of} \: \mathcal{\finevar} \: \}$ based on $\mathcal{M}^d$
\STATE ${\finevar} = \{i \in \Omega \setminus \coarsevar: 
\sigma_i  > \mathbf{tol} \}$
\STATE Perform $\nu$ CR iterations of  Eq.~(\ref{frelax:eq}) with initial guess $u^0$
\ENDWHILE
\end{algorithmic}
\end{algorithm}

Here, $\delta \in (0,1)$ is the tolerance for the approximate CR convergence rate $\rho_{f} $ computed using \eqref{meas1:eq}, and $u^0$ is an initial guess for the compatible relaxation 
iterations.  Further, $d$ is fixed in advance, but more generally it can be adapted at each stage of the CR algorithm.  
We note also that the matrix $A^d$, $d>1$ is not implicitly formed, except for its binary adjacency matrix, as this is all that is needed to construct $\mathcal{M}^d$.
Our choices for these and other parameters of the algorithm used in our numerical experiments are given in Section
\ref{alg:cf}. 
 
Various choices of the candidate set measure, $\sigma_i$, used in determining potential $\coarsevar$-points have been studied in the literature~\cite{ABrandt_2000a,OLivne_2004a,JBrannick_RFalgout}; we follow the  definition in~\cite{JBrannick_RFalgout}:
$$ \sigma_i := \frac{\displaystyle |{u}_{i}^{\nu}|}{\displaystyle \|{u}^{\nu}\|_{\infty}}.$$  
In practice, this measure gives the best overall results for smaller values of $\nu$, say $\nu=5$.  

\subsection{Defining interpolation by algebraic distances}\label{sec:greedy}
Given a set of coarse variables $\coarsevar$ and a set of $\tau$-smoooth test
vectors, $\{v^{(1)},  ..., v^{(k)}\}$, we use algebraic distance to find the coarse interpolatory sets $C_i$, with  
a cardinality  bounded by 
a given caliber $c$. 
More specifically, we consider all possible 
sets of $\coarsevar$-points, $W$ of cardinality up to $c$,  in the $d_{LS}$-ring coarse point neighborhood of a
given $\finevar$-point, $i$,  defined as 
\begin{equation}\label{eq:nbhd}
\mathcal{N}_{d_{LS}, i} :=  \coarsevar \cap \mathcal{V}_{d_{LS}, i}, 
\end{equation}
where $d_{LS}$ is the search depth for constructing least squares interpolation.  That is, $\mathcal{V}_{d_{LS}}$, $\mathcal{E}_{d_{LS}}$, and $\mathcal{G}_{d_{LS}, i}$ are defined as in~\eqref{graph}, where  $d_{LS}$ is a fixed positive integer. 
Thus, the sets of possible interpolatory points can be written as 
$$\mathcal{W} := \{ W \: | \;  \: W \subseteq \mathcal{N}_{d_{LS}, i} \; \mbox{and} \;  |W| \leq c  \}.
$$ Using an exhaustive search of this set, we find the minimizer
 of the least squares functional~\eqref{eq:LSfunc}: for each $i \in \finevar$
$$
C_i = \arg\min_{W \in \mathcal{W}} \mathcal{LS}(p_i(W)) , $$
where
$$
 \mathcal{LS}(p_i(W)) = \sum_{\kappa = 1}^k\omega_{\kappa} \bigg(v_{i}^{(\kappa)} - \frac{1}{a_{ii}}r^{(\kappa)}_{i} - \sum_{j\in W} (p_{i})_j v_{j}^{(\kappa)}\bigg)^{2}, 
$$
 here $p_i$ denotes the minimizer of the least squares problem \eqref{eq:modls1} for the given set $W$.  Thus, we must compute $p_i$ and evaluate $\mathcal{LS}(p_i(W))$ for all possible choices of $W \in \mathcal{W}$.  The row of interpolation, $p_i$, is then chosen as the one that minimizes $\mathcal{LS}(p_i(W))$.

An additional detail of the approach is the penalization of large
interpolatory sets. It is easily shown that for two sets $W^{\prime}
\subset W^{\prime \prime}$, their corresponding minimal least squares
functional values fulfill $\mathcal{LS}^{\prime} \geq
\mathcal{LS}^{\prime \prime}$. In order to keep the interpolation
operator as sparse as possible, we require that the least squares
functional is reduced by a certain factor when increasing the
cardinality of $W$.
That is, a new set of points ${W^{\prime \prime}}$ is preferred
over a set of points ${W^{\prime}}$ with $|W^{\prime \prime}|
> |W^{\prime}|$ if
\begin{equation*}
  \mathcal{LS}^{\prime \prime} < \left(\mathcal{LS}^{\prime}\right)^{\gamma \left( |W^{\prime\prime}|-|W^{\prime}|\right)}.
\end{equation*} Based on numerical experience,  we  choose $\gamma
= 1.5$  which  tends to produce accurate and sparse interpolation
operators for a large class of problems.

The above exhaustive search of all possible combinations of interpolatory sets with cardinality up to some given caliber is one of many possible strategies for selecting $\coarsevar_i$.  In our experience, this is often not necessary, and scanning a small number of  possibilities based on 
the algebraic distance strength measure and caliber one interpolation is sufficient for many problems.
That is, the exact minimization of the LS functional is often not required to obtain sufficiently accurate interpolation.  However, for the anisotropic diffusion problem we consider here, in the case of strong non-grid aligned anisotropies an exhaustive search leads to the best overall results for the low-caliber interpolation, and  we therefore we use this strategy in our numerical experiments.  

\section{Numerical Results}
\label{sec:numerics}
In this section, we present numerical tests that demonstrate the effectiveness of the algebraic distance measure of strength of connection when combined with CR coarsening and LS interpolation.  The tests of the approach consist of a variety of 2D anisotropic diffusion problems discretized using finite differences and finite elements on a $(N+1)\times(N+1)$ uniform grid.

\subsection{Model problem and discretizations}\label{sec:disc}
We consider finite difference and bilinear finite element discretizations of the two-dimensional diffusion operator 
\begin{equation}\label{eq:diff} \mathcal{L}\, u =  \nabla \cdot \mathcal{K} \nabla u, 
\end{equation}
with anisotropic diffusion coefficient 
\begin{equation}  \mathcal{K}  =       \biggl( \begin{matrix}   \cos \alpha & -\sin \alpha \\
      \sin \alpha & \cos \alpha\ \\
   \end{matrix}\biggr)     \biggl( \begin{matrix} 
         1 & 0 \\
         0 & \epsilon \\
       \end{matrix} \biggr)   \biggl( \begin{matrix}   \cos \alpha & \sin \alpha \\
      -\sin \alpha & \cos \alpha\ \\
   \end{matrix}\biggr), 
   \label{eq:a}
   \end{equation} where
$0 <\epsilon<1$ and $0 \leq \alpha  < 2 \pi$. Changing   variables
\begin{equation}   \biggl( \begin{matrix}       \xi \\
      \eta  \\
   \end{matrix} \biggr)
   = \biggl( \begin{matrix}   \cos \alpha & \sin \alpha \\
      -\sin \alpha & \cos \alpha\ \\
   \end{matrix}\biggr) \biggl( \begin{matrix}       x \\
      y  \\
   \end{matrix} \biggr),
   \label{eq:change_var}
   \end{equation}
yields  strong connections aligned with direction $\xi$:
  \begin{equation} \mathcal{L} \, u(\xi, \eta) = u_{\xi \xi} + \epsilon u_{\eta \eta}.
  \end{equation}
An equivalent     
formulation in the $(x,y)$ coordinates  is given by 
\begin{eqnarray}\label{eq:scal} \mathcal{L}u(x, y) = a u_{xx} + 2b u_{xy} + c u_{yy},  \end{eqnarray}
where $\displaystyle{a  = \cos^2 \alpha + \epsilon \sin^2 \alpha, \, b = \frac{(1-\epsilon)}{2} \,\sin 2\alpha}$, and $c = \sin^2 \alpha + \epsilon \cos^2 \alpha$.

In formulating a finite difference  discretization of (\ref{eq:scal}), we consider a standard five-point discretization for the Laplacian term and then define the discretization of the $u_{xy}$ term using lower-left  and upper-right  neighbors and not the lower-right  and upper-left  ones  for each fine-grid point $i \in \Omega$, yielding an overall stencil that includes seven nonzero values.   
This gives a suitable discretization for $\alpha = \pi /4$ 
with stencil $u_{xy}$ given by
\begin{equation}
\label{fd_seven} \tilde{S}_{xy} = \dfrac{1}{2h^2}  
\Biggl(
\begin{matrix}
0 & - 1 & 1 \\
0 & 1 & - 1 \\
0 & 0 & 0 \\
\end{matrix}
\Biggr) + \dfrac{1}{2h^2}  \Biggl(
\begin{matrix}
0 & 0 & 0 \\
-1 & 1 & 0 \\
1 & -1 & 0 \\
\end{matrix}
\Biggr).               
\end{equation}        
In contrast, for the $\alpha = - \pi/4$ case, the direction of the anisotropy can not be captured by our chosen discretization and yields a matrix that is far from an M-matrix.
For example, taking $\epsilon = .1$, the resulting system matrix for this choice of $\alpha$ has stencil
$$
S_A =  \dfrac{1}{2h^2} \Biggl(
\begin{matrix}
 & -1 & .45 \\
-1 & 3.1 & -1 \\
.45 & -1 &  \\
\end{matrix}
\Biggr),
$$ and, thus, is not even approximately an $M$-matrix, which makes MG solution of this problem challenging.  (\textcolor{black}{We note that the $\alpha = \pi/8$ case we consider in our tests also yields non M-matrices for strong anisotropies.})  Here, some of the off-diagonal entries in $A$ are positive and, hence, the heuristics motivating the classical definition of strength of connection are not applicable.  We consider this extreme case, although it is unlikely to arise in practice, to demonstrate the robustness of the BAMG approach as a coarsening strategy for the targeted anisotropic problems.   
We mention that on a structured grid, our chosen seven-point discretization is equivalent to the finite element discretization of the same
elliptic boundary value problem, using triangular finite elements.

In addition, we consider the bilinear finite element discretization of the same model problem, again on a $(N+1)\times(N+1)$ uniform grid.  Letting $\phi_j(x,y)$ denote
a standard bilinear basis function that is one at node $j$ and zero at all other nodes, and writing the solution as $\sum_{j=1}^n u_j \phi_j$, the weak form of \eqref{eq:diff}
is given by
\begin{equation}\label{eq:modelFE}
- \sum_{j=1}^n u_j \int_0^1 \int_0^1 (\nabla \phi_i) \cdot \bigg[ \left(\begin{matrix}
    a  & b   \\
     b &  c
\end{matrix}\right)
\nabla \phi_j\bigg] dx \; dy =  \int_0^1 \int_0^1 f \phi_i dx \; dy, \quad i = 1,\hdots, n.
\end{equation}
And, the corresponding nine point stencil for the global stiffness matrix $A$ is given by
\begin{equation}\label{eq:modelFEStencil}
S_A =  \left(\begin{matrix}
    -a +3b-c  & 2(a-2c) & -a-3b-c   \\
     2(-2a+c) &  8(a+c) & 2(-2a+c) \\
     -a-3b-c & 2(a-2c)  & -a +3b-c  
\end{matrix}\right).
\end{equation}
As an example, when $\theta = \frac{\pi}{4}$, we have 
\begin{equation}\label{eq:modelFEStencilpi4}
S_A =  \left(\begin{matrix}
    \frac12 -\frac52 \epsilon  & -1 - \epsilon &  -\frac52 +\frac12 \epsilon   \\
    -1 - \epsilon &  8 + 8 \epsilon & -1 - \epsilon \\
      -\frac52 +\frac12 \epsilon & -1 - \epsilon  & \frac12 -\frac52 \epsilon 
\end{matrix}\right) \quad \longrightarrow \quad  \left(\begin{matrix}
    \frac12   & -1 &  -\frac52    \\
    -1  &  8  & -1  \\
      -\frac52 & -1  & \frac12 \\
\end{matrix}\right) \quad \text{as} \quad \epsilon \longrightarrow 0,
\end{equation}
which again is far away from an $M$ matrix.

%

\subsection{Formulating a two-level coarsening algorithm}\label{alg:cf}
Our aim is to study the robustness of the algebraic distance notion of
 strength of connection for grid and non-grid aligned anisotropy.  We focus our tests on 
 the seven-point finite difference and nine-point finite element discretizations introduced earlier in this section for various values
 of the anisotropy angle, $\alpha$, and the anisotropy coefficient, $\epsilon$.
 
In all tests, coarse grids are chosen using the compatible relaxation approach discussed in Section \ref{sec:cr},  interpolation  operators are then computed using the LS approach given in Section \ref{sec:greedy} with the caliber set to $c=1$ or $c=2$.  
The stopping tolerance for CR steps is set as $\rho_f < \delta = 0.7$, the number of CR sweeps is chosen as $\nu = 5$, and the tolerance for the candidate set measure is  ${ \bf tol} = 1-\rho_f$.  The larger choice of stopping tolerance generally leads to more aggressive coarsening whereas our choice of $\bf{tol}$ ensures that points are added to the candidate set more sparingly at later stages of the CR coarsening process.
We fix the strength of connection parameter used in forming the strength graph at $\theta_{ad} = .5$ and vary the graph distance, $d$, used to define the graph $\mathcal{G}_d$, from which $\mathcal{M}^d$ defined as in \eqref{eq:strength_graph} is constructed.  
We note that generally the overall quality of the grids the algorithm produces depends only mildly on the choice of $\theta_{ad}$.

The search depth, used in defining the greedy approach for choosing LS interpolation as in \eqref{eq:nbhd}, is fixed as $d_{LS} = d + 2$.  Taking the value of the search 
depth larger than the coarsening depth allows the approach to scan a larger number of possible interpolatory sets in forming long-range interpolation (whenever the problem requires it).  In this way, the LS scheme for constructing interpolation is able to better follow a wider range of values of the anisotropy angle, $\alpha$.   We mention, in addition, that for the anisotropic problems we consider, it is important that the set of test vectors used in guiding the coarsening
consists of the {\it characteristic components}, eigenvectors with small eigenvalues, induced by the anisotropy angle $\alpha$, especially for angles for which longer-range interpolation is needed in order 
to follow the anisotropy. 
In our two-grid tests, where only relaxation is used to compute the test vectors, we thus apply 40 iterations of Gauss Seidel to seven distinct random initial guesses plus the constant vector in order to compute the eight test vectors used to construct LS interpolation in~\eqref{eq:modls1}.  The weights in \eqref{eq:modls1} are defined as 
\newcommand{\innerprod}[3][2]{\langle #2,#3 \rangle_{#1}}
\begin{equation}
\omega_{\kappa} =   \frac{\innerprod[]{v^{(\kappa)}}{v^{(\kappa)}}}{\innerprod[]{Av^{(\kappa)}}{v^{(\kappa)}}}.
\label{eq:ls:weights1}
\end{equation} 

\textcolor{black}{The number of test vectors and especially the amount
of relaxation required to compute them can be reduced by replacing
a single-grid relaxation procedure by a multilevel bootstrap cycling scheme~\cite{BAMG2010}, 
as we show in the next section where we reduce the number of relaxation steps
from 40 to eight in the multilevel setup.}

When presenting results of the solver constructed by the resulting BAMG setup algorithm we use two pre- and post- Gauss Seidel relaxation sweeps on all grids except the coarsest, where a direct solver is used.  Here, the use of two pre- and post-smoothing steps is motivated by the fact that the setup algorithm generally produces aggressive coarsening.  The estimates of the asymptotic convergence rates are computed as 
\begin{equation*}
\rho = \frac{\|e^{\eta}\|_A}{\|e^{\eta-1}\|_A},
\end{equation*}
where $e^\eta$ is the error after $\eta=100$ MG iteration applied to the homogeneous system, 
starting with random initial approximations.  We also report the operator complexity ratio
 $$\gamma_o =  \dfrac{\sum_l nnz(A_l)}{nnz(A)},$$
 and the grid complexity ratio 
$$ \gamma_g =  \dfrac{\sum_l |C_l|}{|\Omega|},$$
where $C_0 = \Omega$. 

\subsection{Bootstrap coarsening}
\begin{figure}[htb!]
  \subfigure[$\alpha = 0, \rho = .23, \rho_f = .27, \gamma_o = 1.482, \gamma_g = .491$ \label{fig:aggs10FD}]{\includegraphics[scale = 0.40]{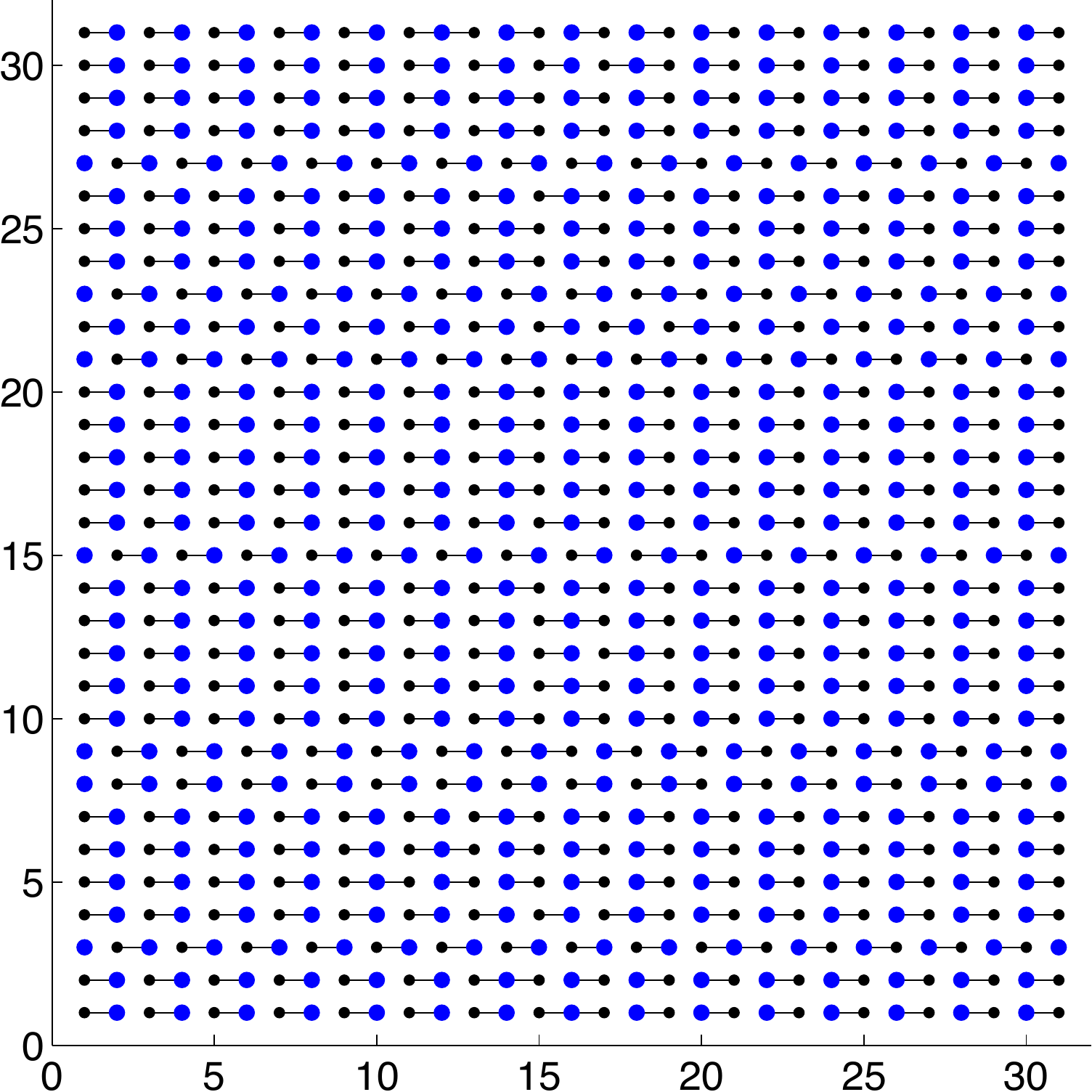}}
  \hfill
  \subfigure[$\alpha = \pi/4, \rho = .17, \rho_f = .65, \gamma_o =1.473, \gamma_g = .488$\label{fig:aggs1p4FD}]{\includegraphics[scale = 0.40]{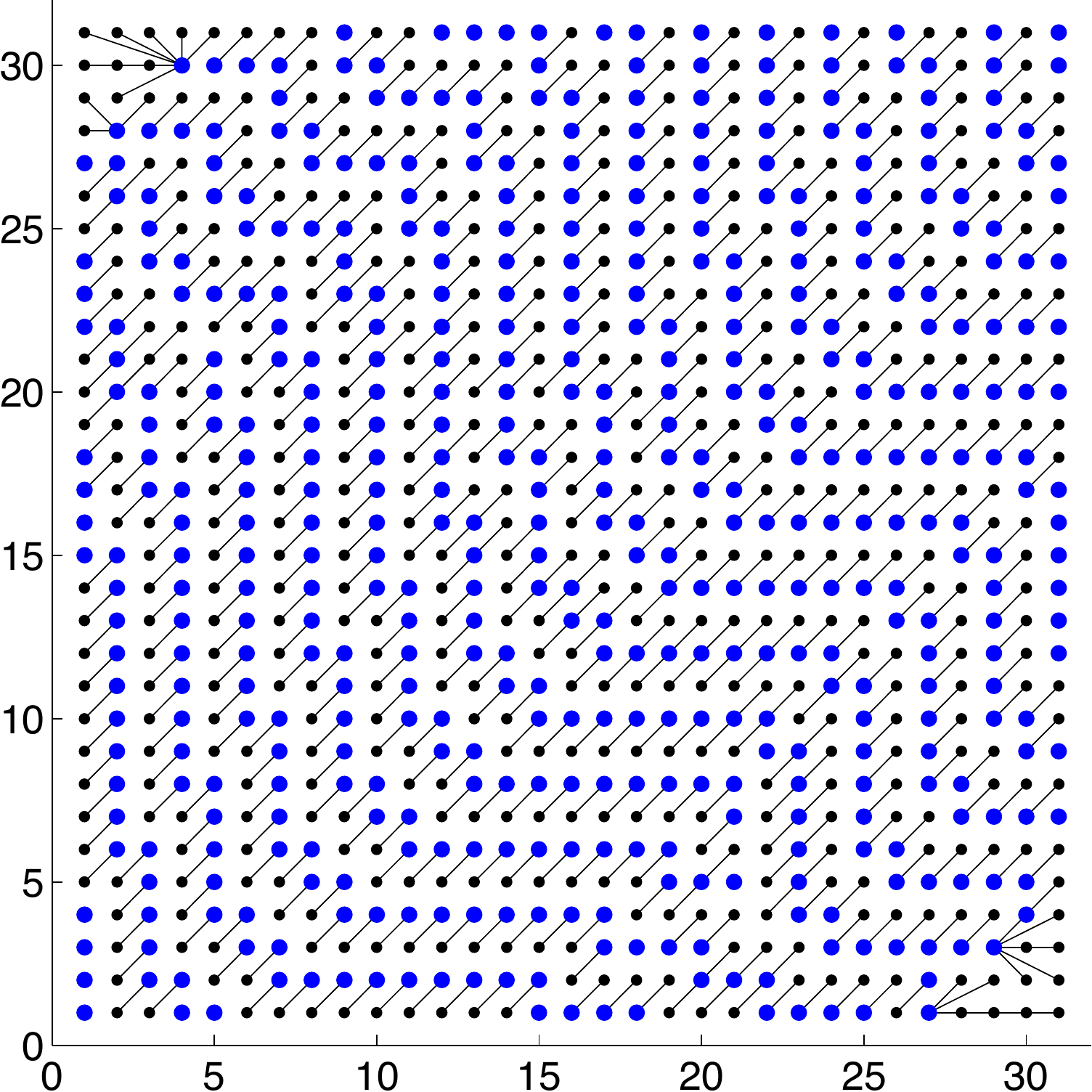}}
\hfill
  \subfigure[$\alpha = -\pi/4, \rho = .72, \rho_f = .72, \gamma_o = 1.361, \gamma_g = .362 $ \label{fig:aggs1mp4FD}]{\includegraphics[scale = 0.40]{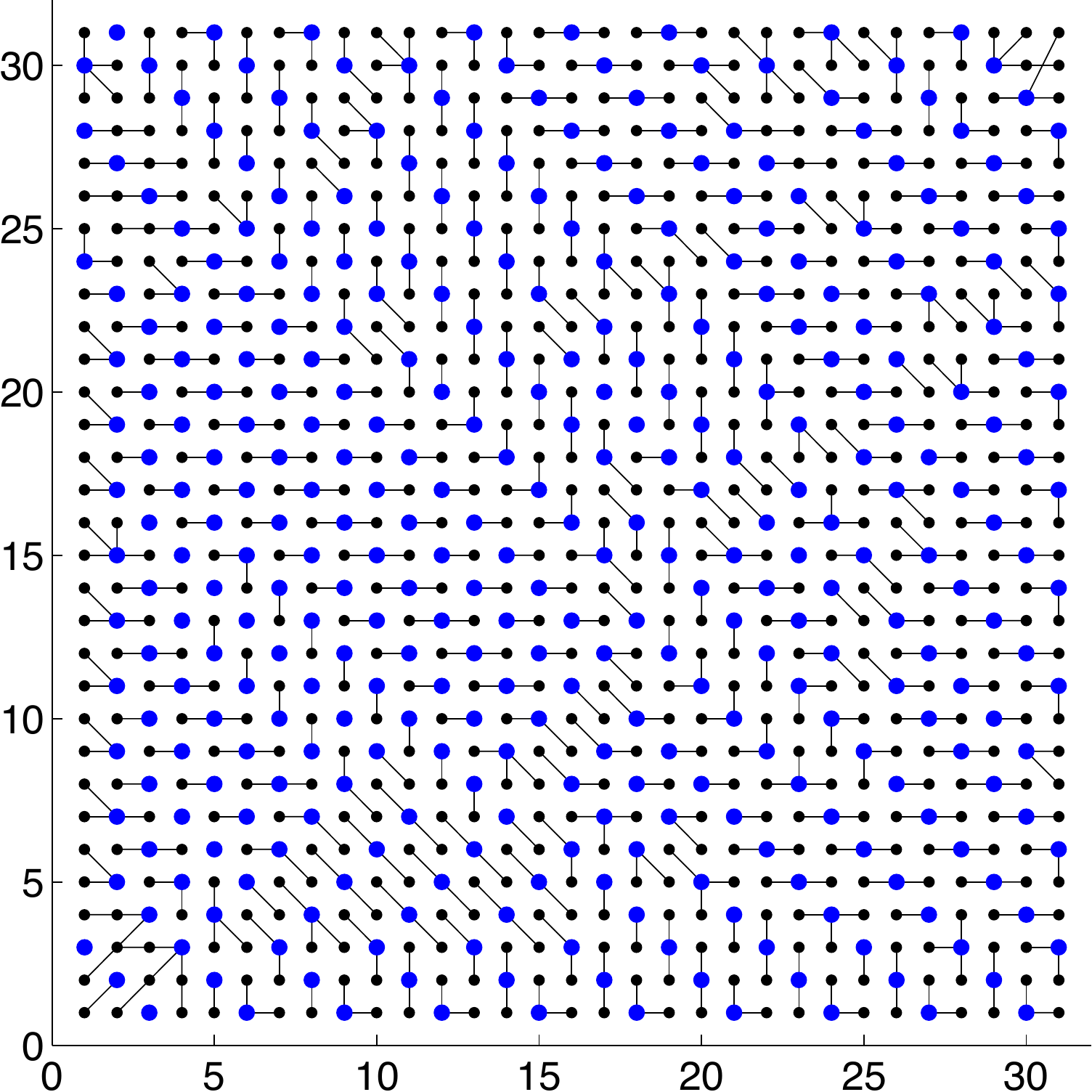}}
 \hfill
    \subfigure[$\alpha = \pi/8,   \rho = .55, \rho_f = 0.59, \gamma_o  = 1.404, 
\gamma_g = 0.405
 $\label{fig:aggs1p16FD}]{\includegraphics[scale = 0.40]{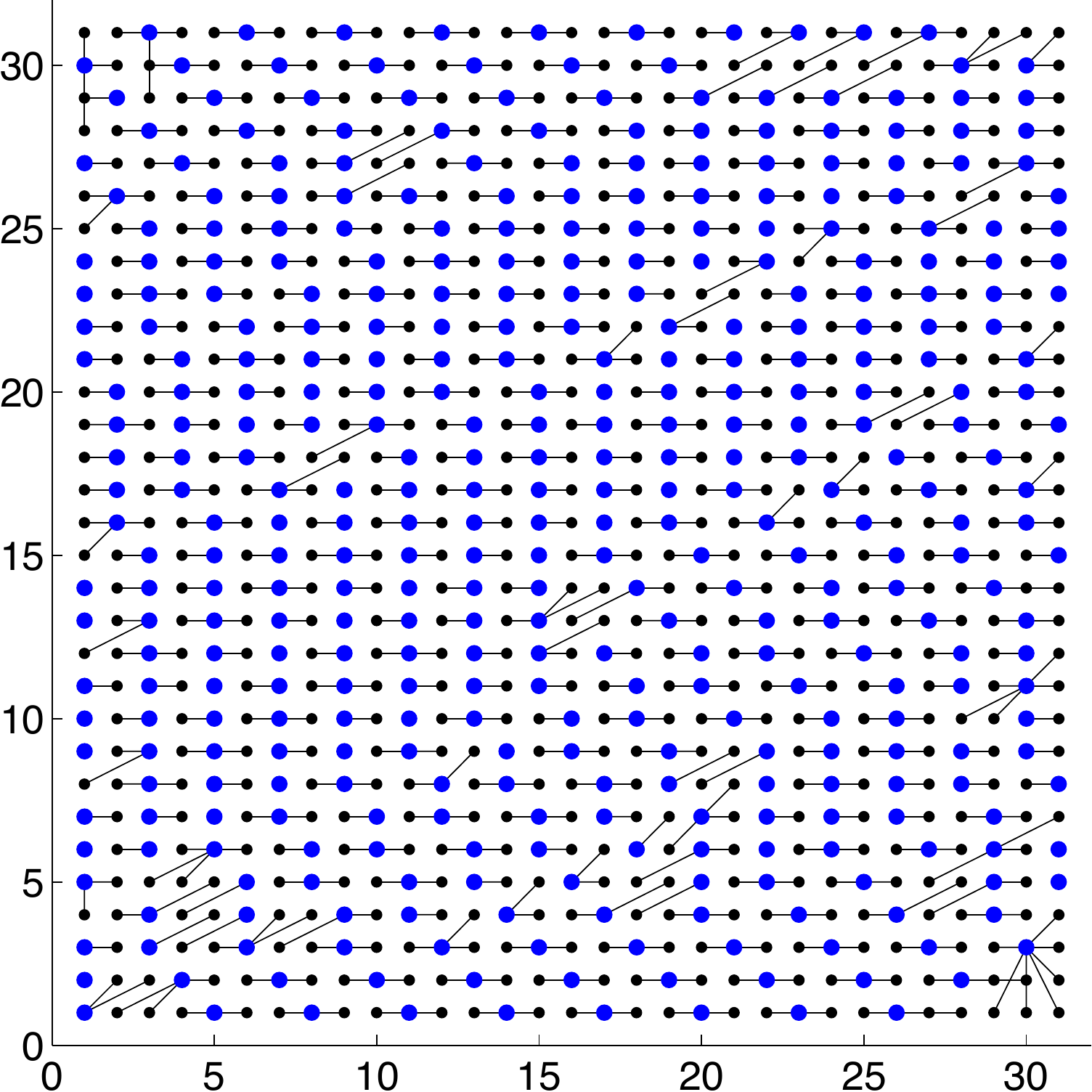}}
    \caption{Coarse grids and caliber $c=1$ interpolation patterns for the finite difference discretization with $h= 1/32$ for various choices of $\alpha$, using the graph of $A$, i.e., $d = 1$ and $d_{LS} = 3$, to define the strength matrix.  Here, the smaller circles are $F$-points and the larger circles are $C$-points \label{fig:aggs1FD}}
\end{figure}

\begin{figure}[htb!]
  \subfigure[$\alpha = 0, \rho = .36, \rho_f = .33, \gamma_o = 1.442, \gamma_g = .499$ \label{fig:aggs10}]{\includegraphics[scale = 0.40]{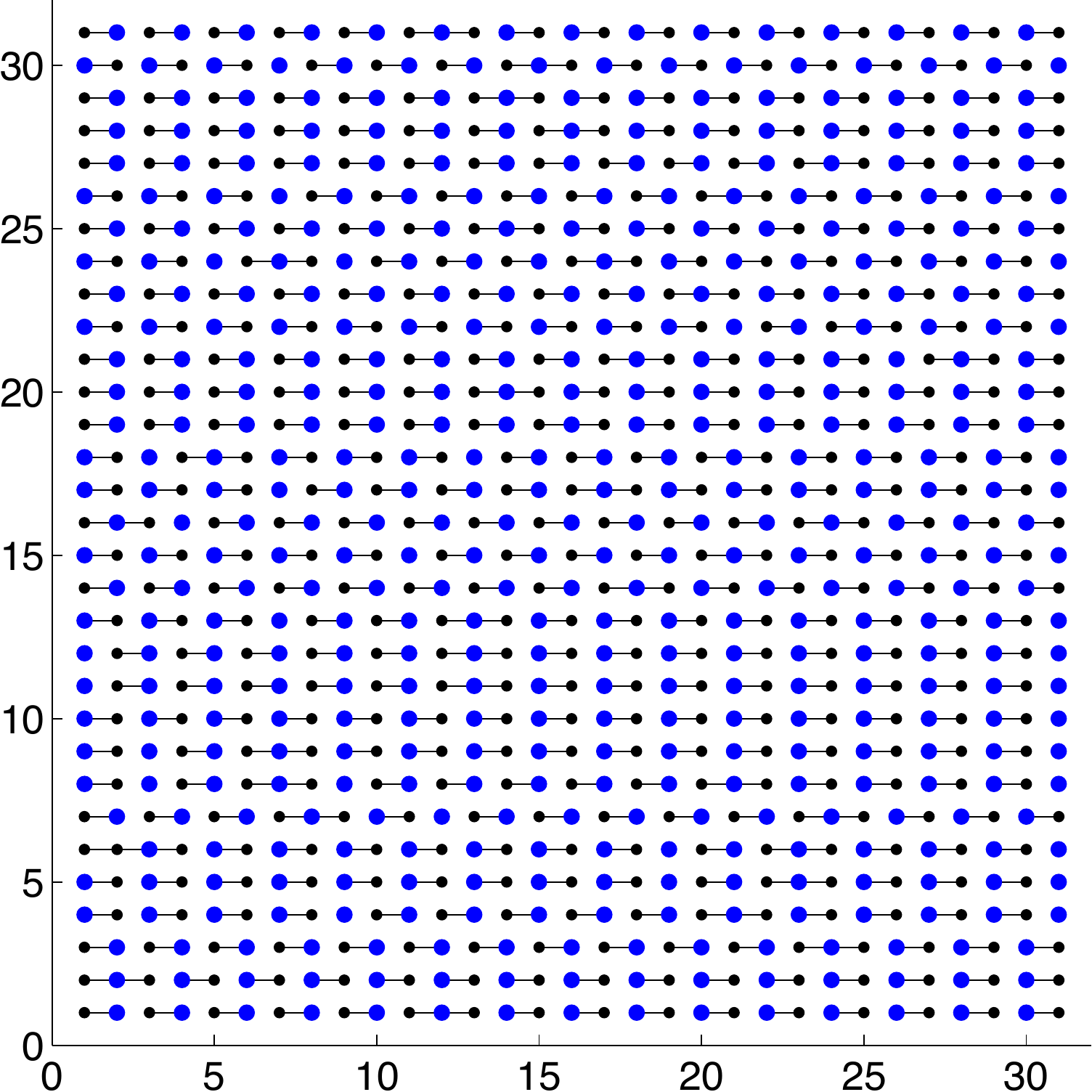}}
  \hfill
  \subfigure[$\alpha = \pi/4, \rho = .28, \rho_f = .41, \gamma_o =1.536, \gamma_g = .495$\label{fig:aggs1p4}]{\includegraphics[scale = 0.40]{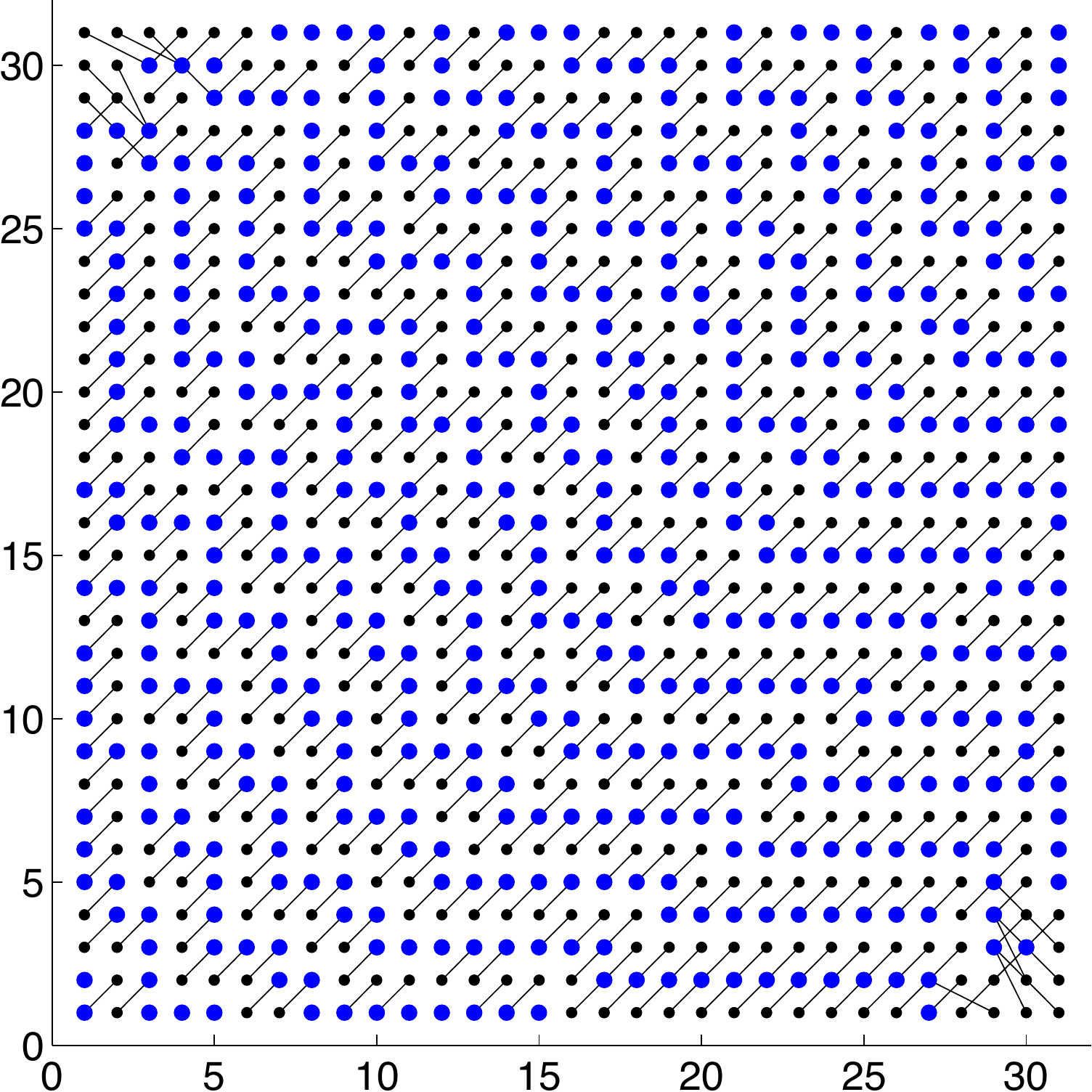}}
\hfill
  \subfigure[$\alpha = -\pi/4, \rho = .24, \rho_f = .33, \gamma_o = 1.533, \gamma_g = .487 $\label{fig:aggs1mp4}]{\includegraphics[scale = 0.40]{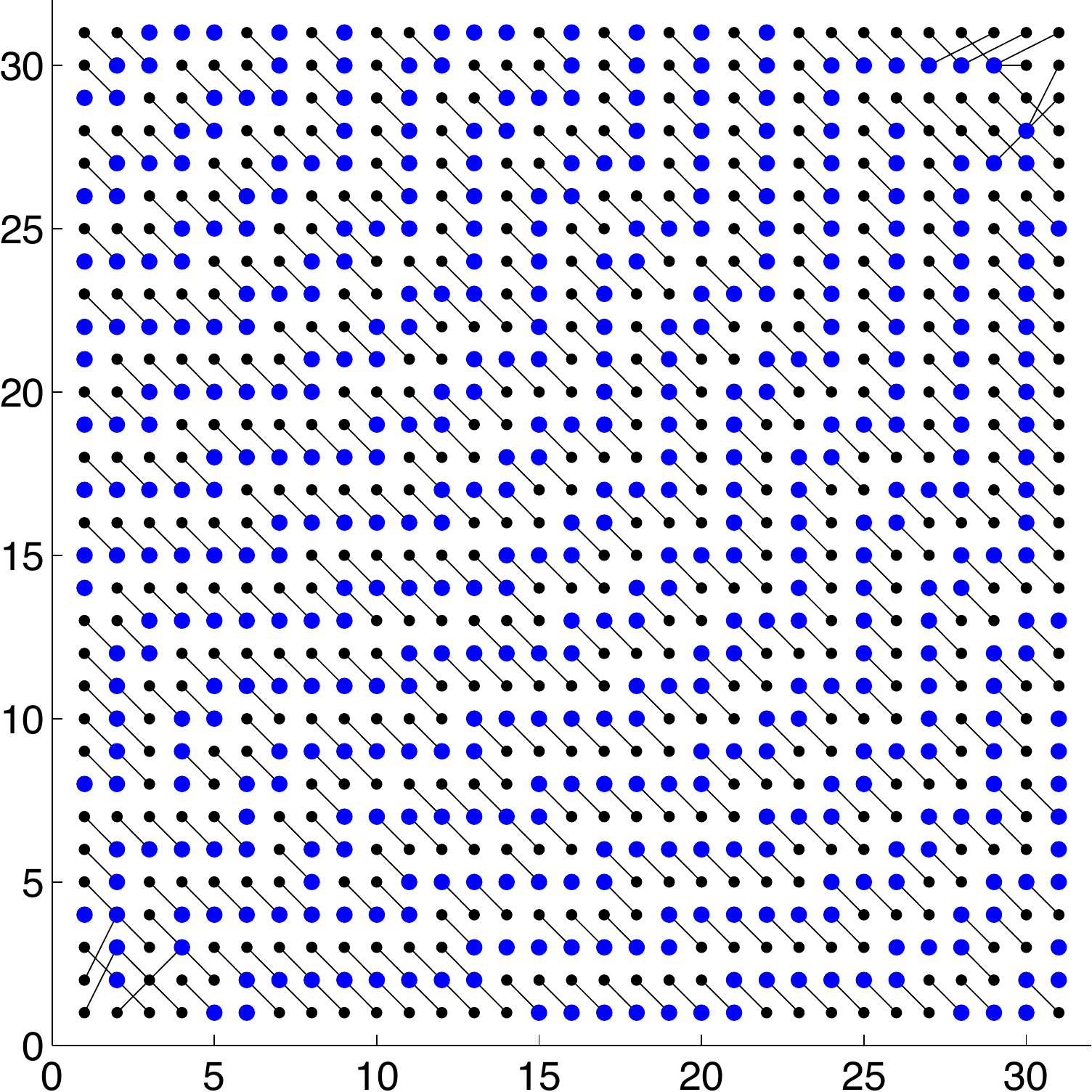}}
 \hfill
    \subfigure[$\alpha = \pi/8,   \rho = 0.54, \rho_f = 0.66, \gamma_o  = 1.360, 
\gamma_g = ..403
 $\label{fig:aggs1p16}]{\includegraphics[scale = 0.40]{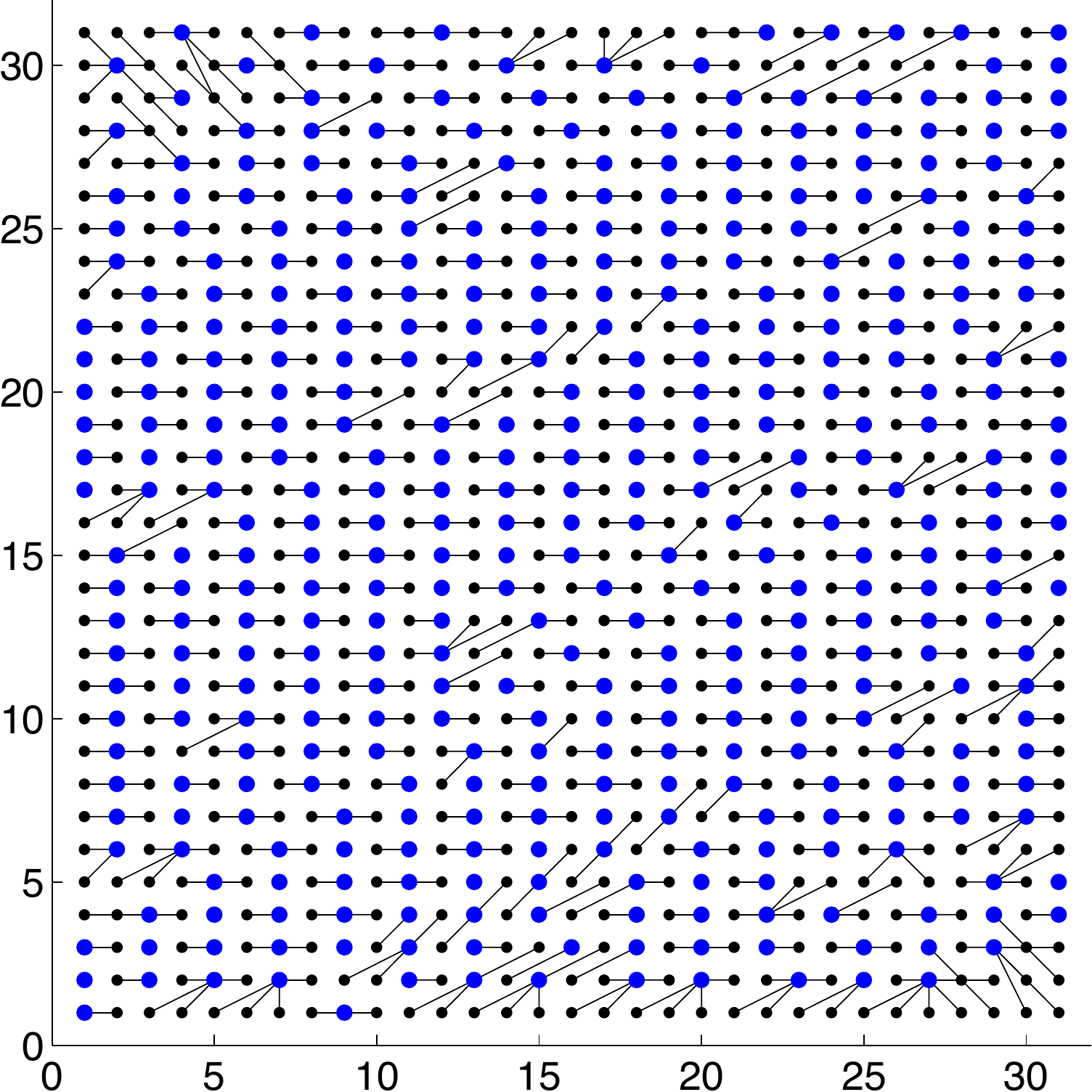}}
    \caption{Coarse grids and caliber $c=1$ interpolation patterns for the bilinear finite element discretization with $h= 1/32$ for various choices of $\alpha$, using the graph of $A$, i.e., $d = 1$ and $d_{LS} = 3$, to define the strength matrix.  Here, the smaller circles are $F$-points and the larger circles are $C$-points \label{fig:aggs1FE}}
\end{figure}

\begin{figure}[htb!]
  \subfigure[$\alpha = 0, \rho = .01, \rho_f = .27, \gamma_o = 1.610, \gamma_g = .491$ \label{fig:s10}]{\includegraphics[scale = 0.40]{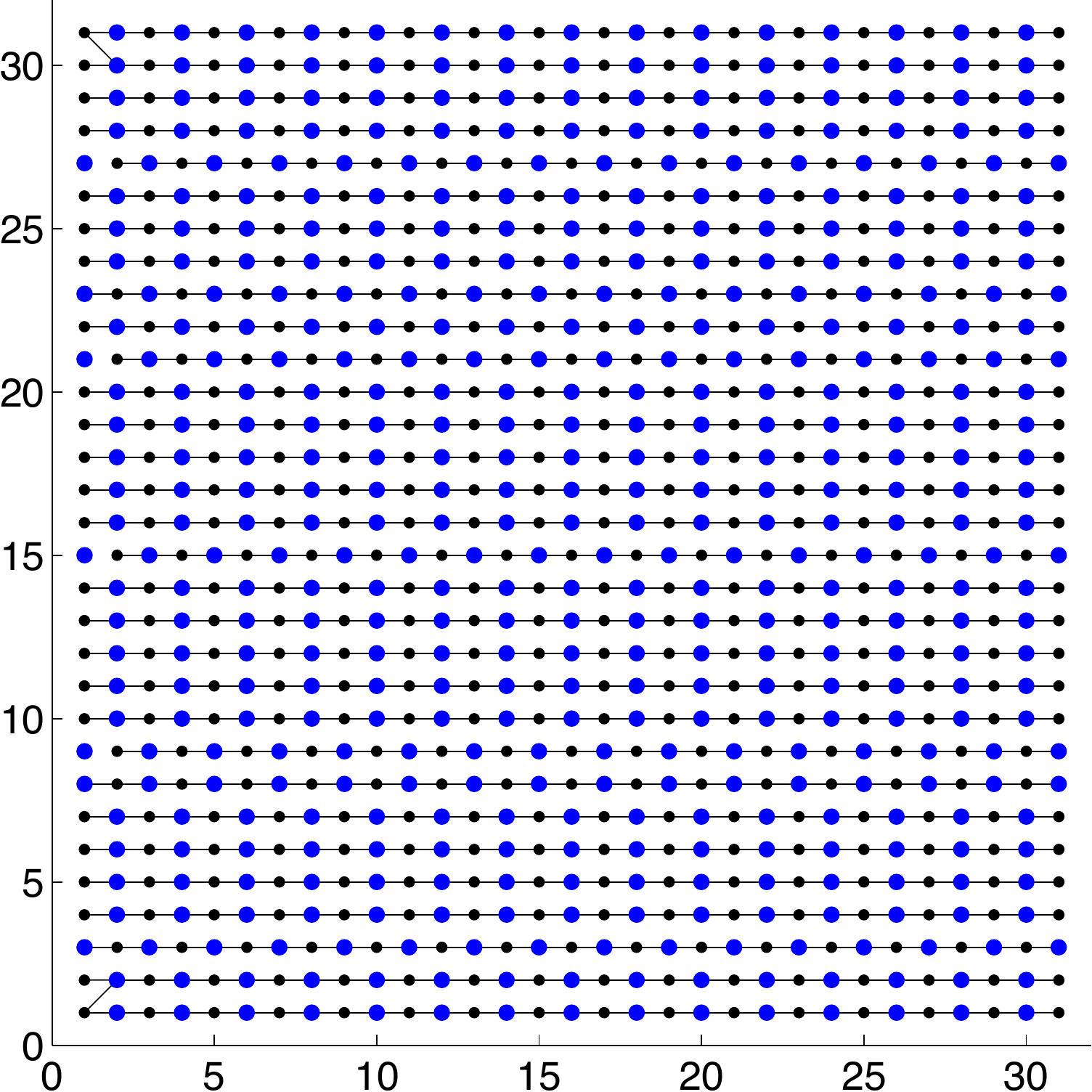}}
  \hfill
  \subfigure[$\alpha = \pi/4, \rho = .03, \rho_f = .65, \gamma_o =1.574, \gamma_g = .488$\label{fig:s1p4}]{\includegraphics[scale = 0.40]{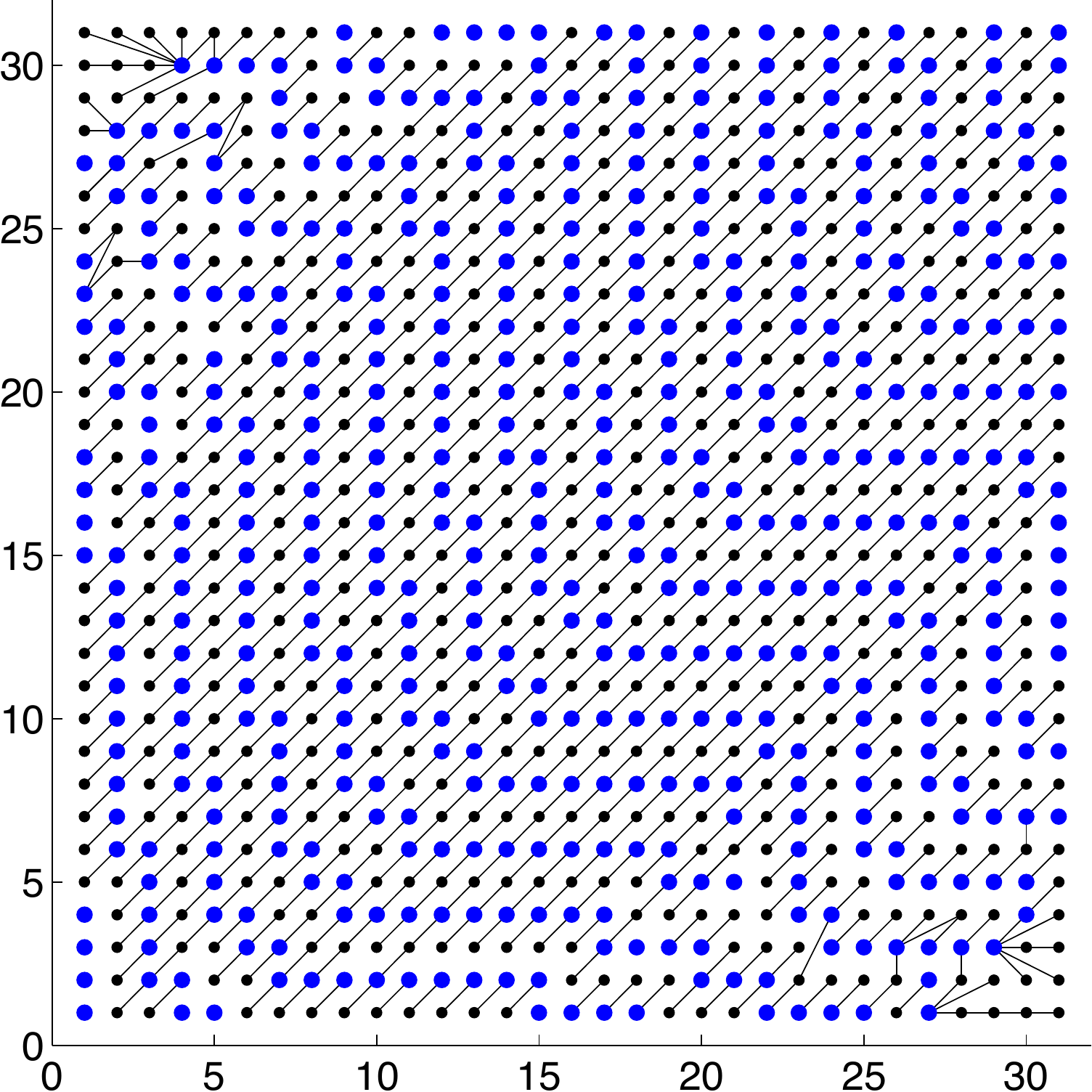}}
\hfill
  \subfigure[$\alpha = -\pi/4, \rho = .32, \rho_f = .72, \gamma_o = 1.509, \gamma_g = .362 $ \label{fig:s1mp4}]{\includegraphics[scale = 0.40]{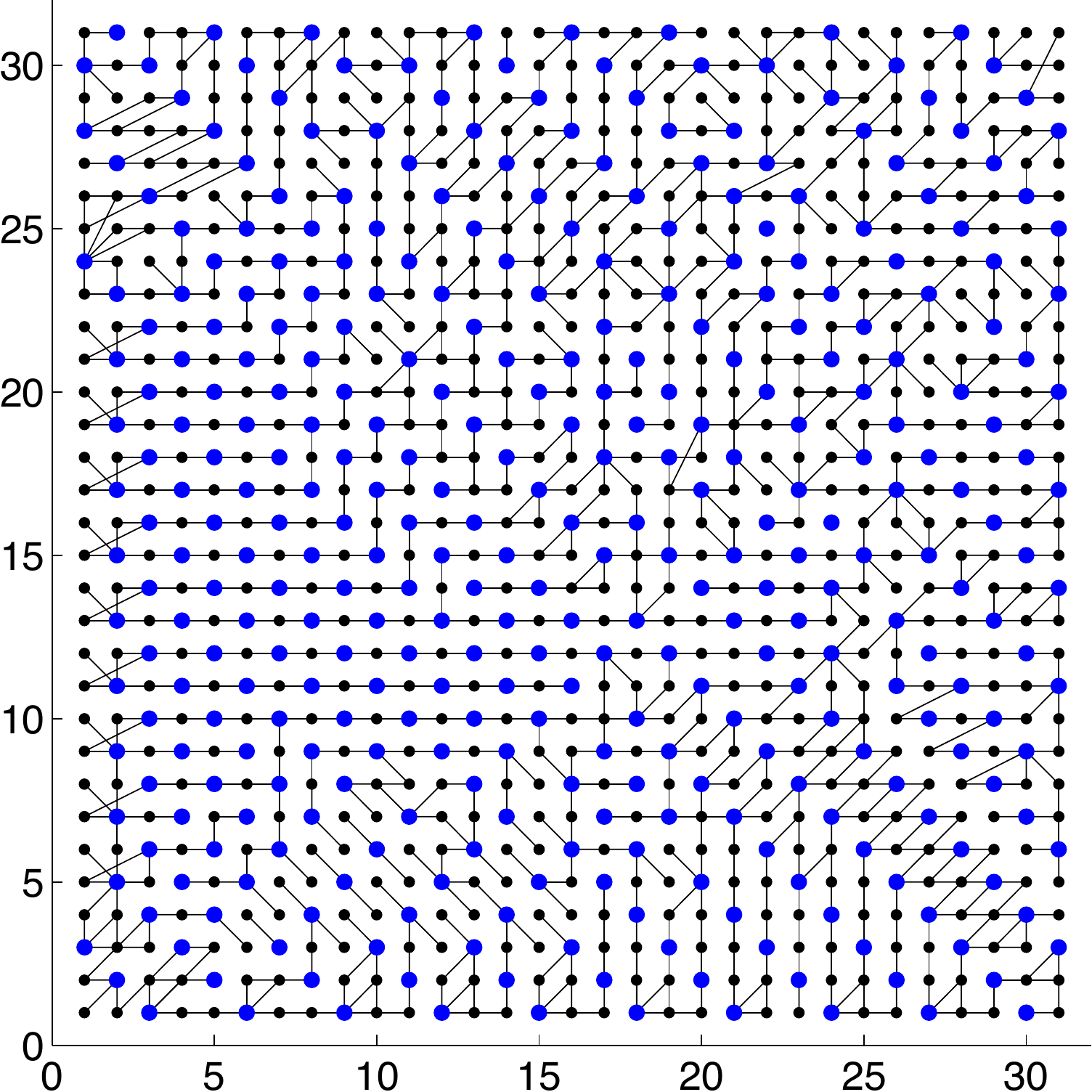}}
 \hfill
    \subfigure[$\alpha = \pi/8,   \rho = 0.18, \rho_f = 0.59, \gamma_o  = 1.573, 
\gamma_g = .405
 $\label{fig:s1p16}]{\includegraphics[scale = 0.40]{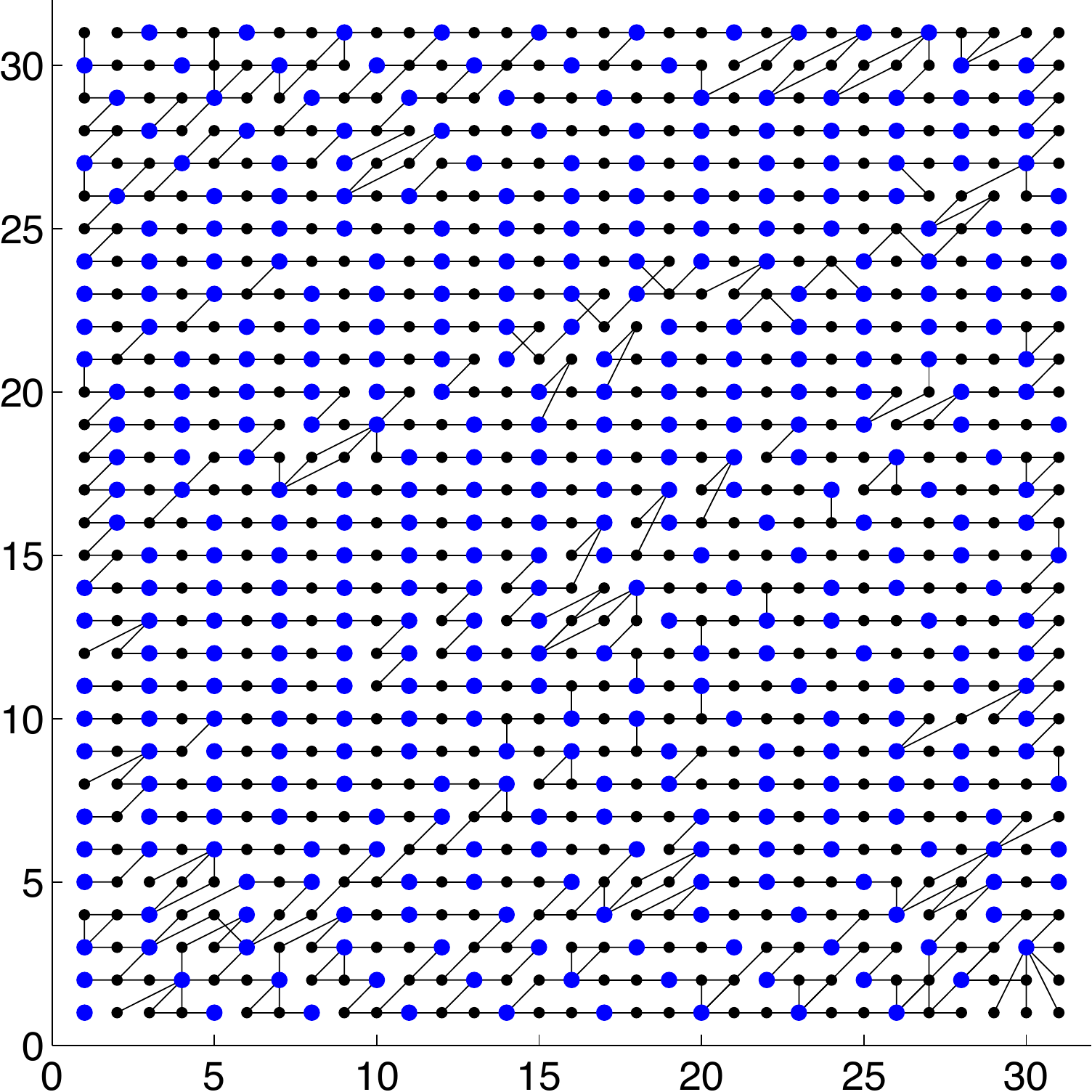}}
    \caption{Coarse grids and caliber $c=2$ interpolation patterns for the finite difference discretization with $h= 1/32$ for various choices of $\alpha$, using the graph of $A$, i.e., $d = 1$ and $d_{LS} = 3$, to define the strength matrix.  Here, the smaller circles are $F$-points and the larger circles are $C$-points \label{fig:s1}}
\end{figure}

\begin{figure}[htb!]
  \subfigure[$\alpha = 0, \rho = .08, \rho_f = .33, \gamma_o = 1.531, \gamma_g = .499$ \label{fig:s10FE}]{\includegraphics[scale = 0.40]{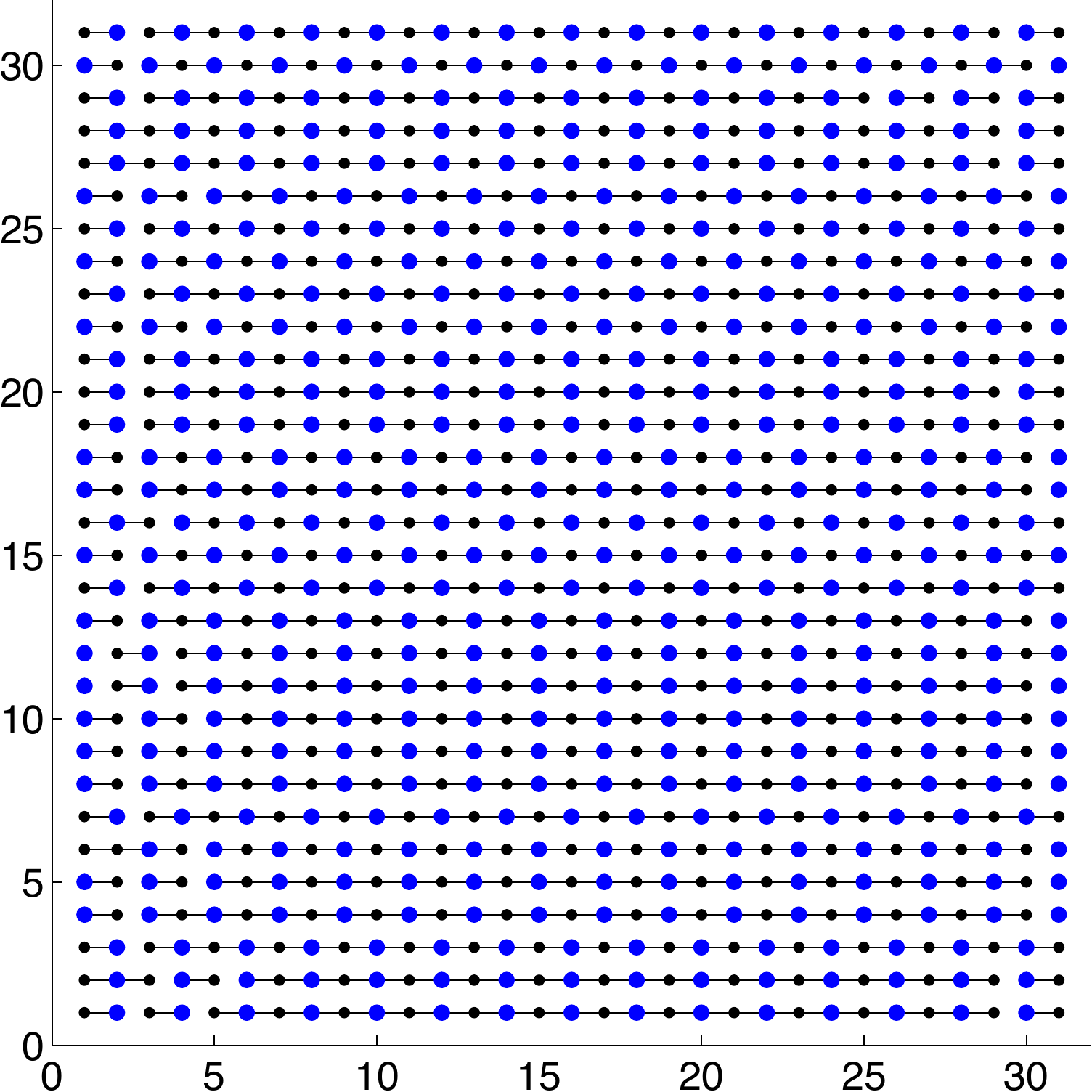}}
  \hfill
  \subfigure[$\alpha = \pi/4, \rho = .03, \rho_f = .41, \gamma_o = 1.660, \gamma_g = .495$\label{fig:s1p4FE}]{\includegraphics[scale = 0.40]{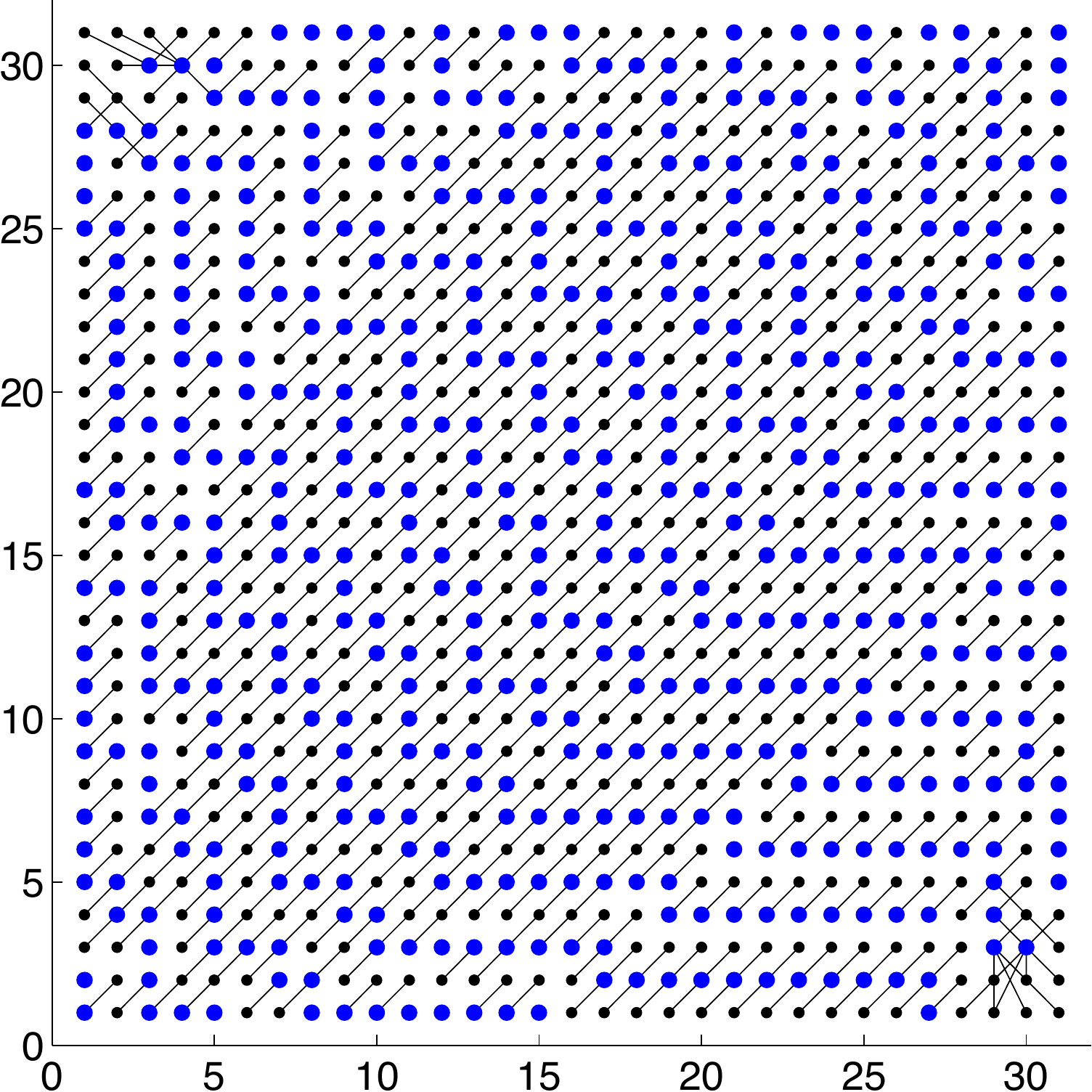}}
\hfill
  \subfigure[$\alpha = -\pi/4, \rho = .06, \rho_f = .33, \gamma_o = 1.648, \gamma_g = .487 $ \label{fig:s1mp4FE}]{\includegraphics[scale = 0.40]{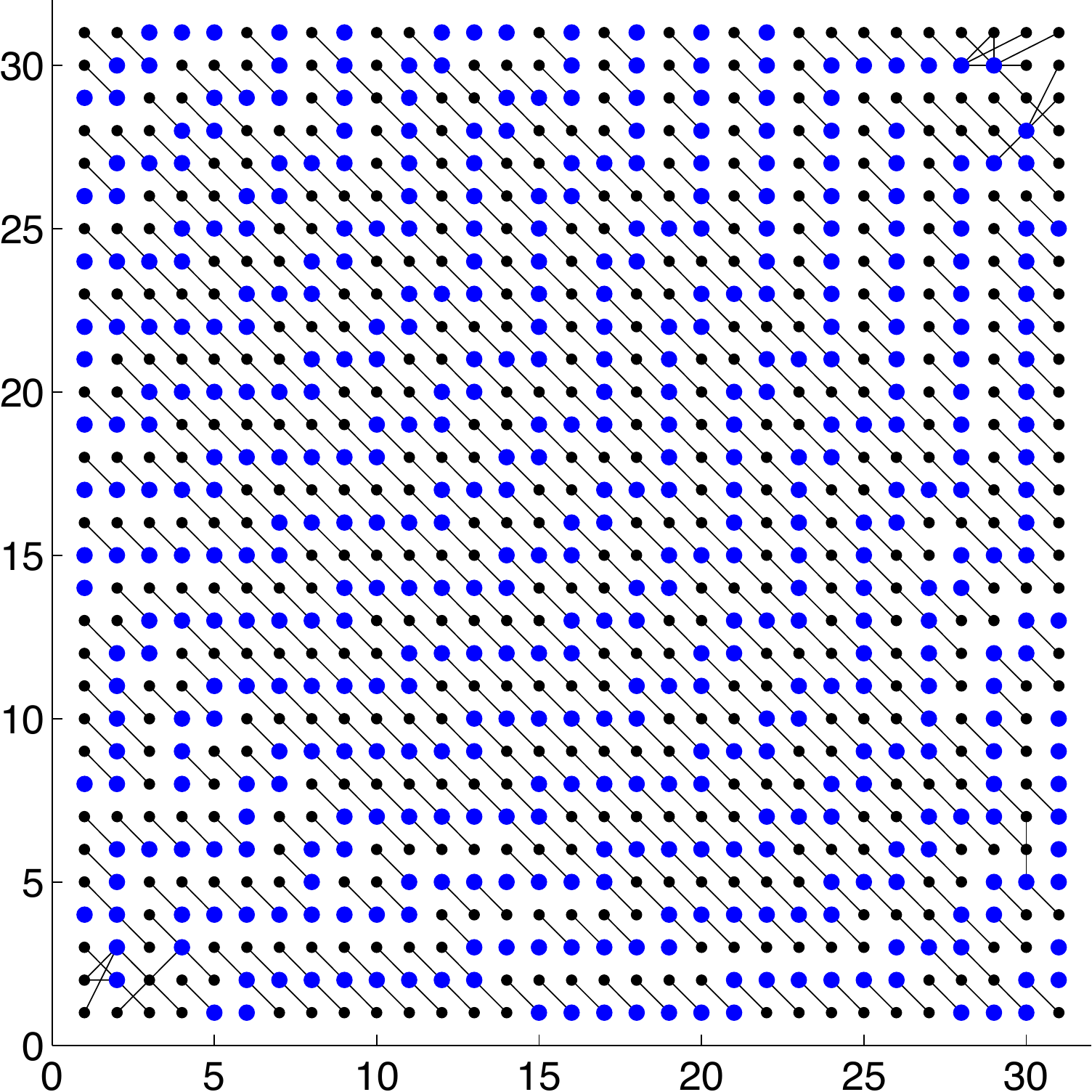}}
 \hfill
    \subfigure[$\alpha = \pi/8,   \rho = 0.24, \rho_f = .59, \gamma_o  = 1.456, 
\gamma_g = .403 $
\label{fig:s1p16FE}] {\includegraphics[scale = 0.40]{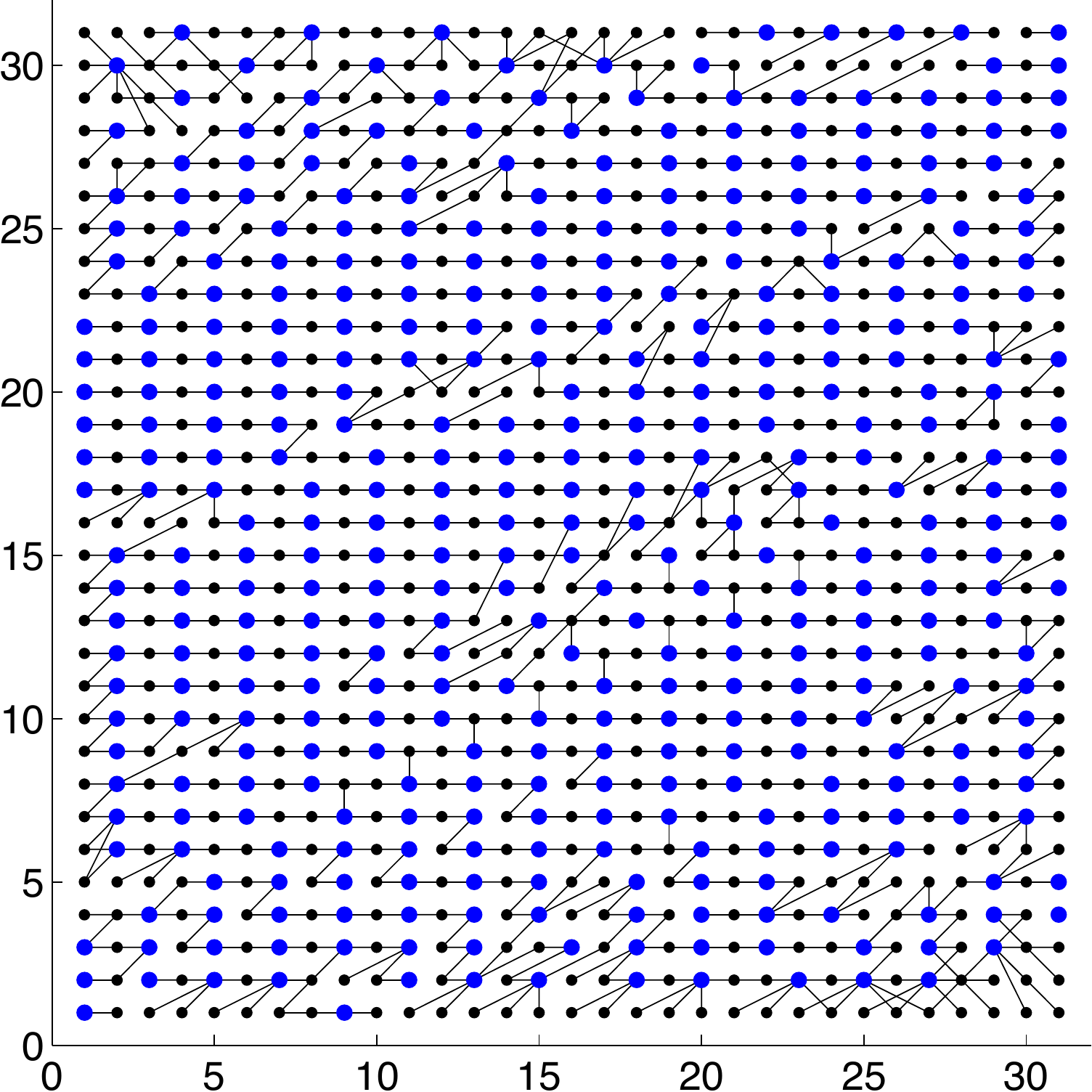}}
    \caption{Coarse grids and caliber $c=2$ interpolation patterns for the bilinear finite element discretization with $h= 1/32$ for various choices of $\alpha$, using the graph of $A$, i.e., $d = 1$ and $d_{LS} = 3$,   to define the strength matrix.  Here, the smaller circles are $F$-points and the larger circles are $C$-points \label{fig:s1FE}}
\end{figure}

\begin{figure}[htb!]
  \subfigure[$\alpha = 0, \rho = .18, \rho_f = .25, \gamma_o = 1.339, \gamma_g = .352$ \label{fig:aggs20FD}]{\includegraphics[scale = 0.40]{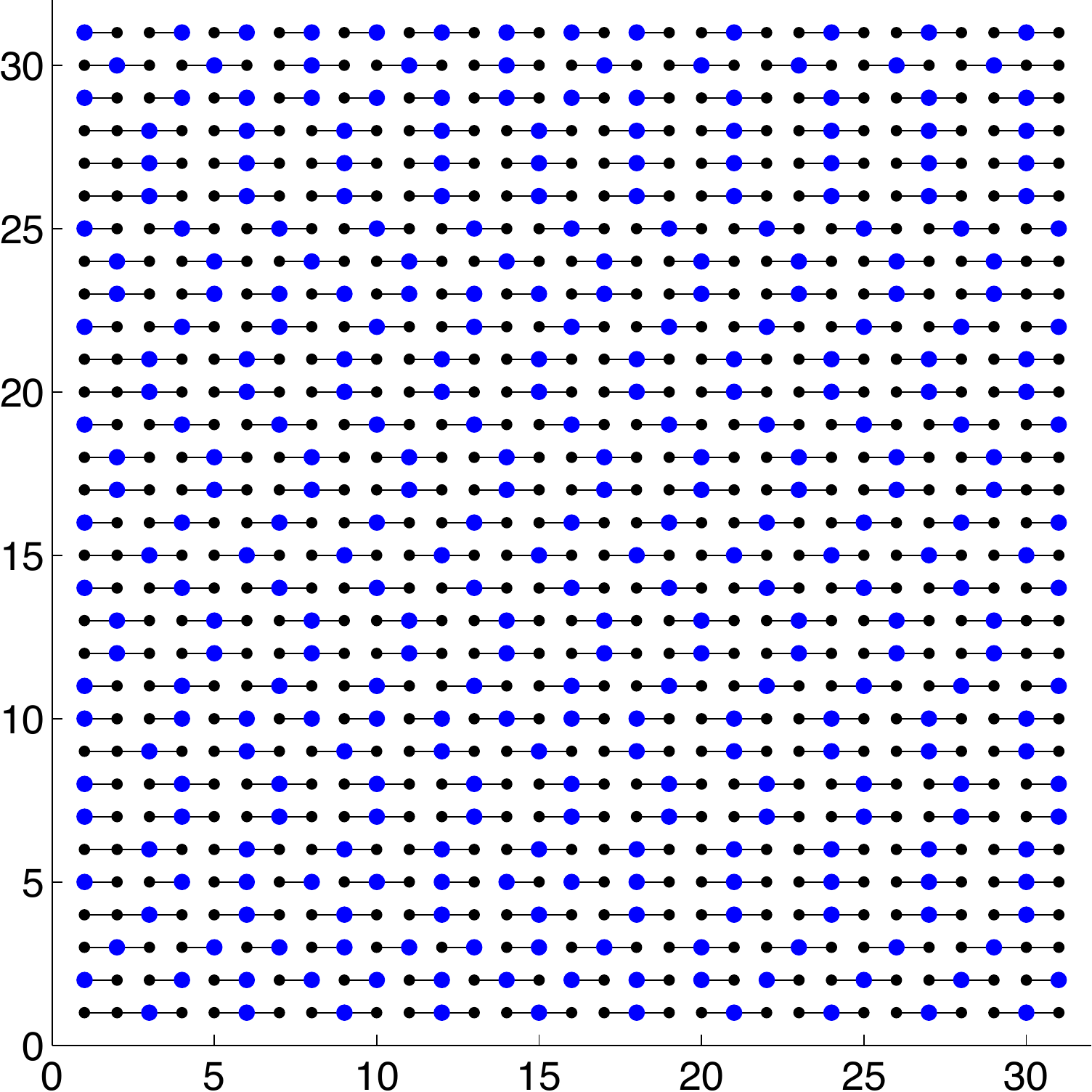}}
  \hfill
  \subfigure[$\alpha = \pi/4, \rho = .14, \rho_f = .65, \gamma_o =1.311, \gamma_g = .332$\label{fig:aggs2p4FD}]{\includegraphics[scale = 0.40]{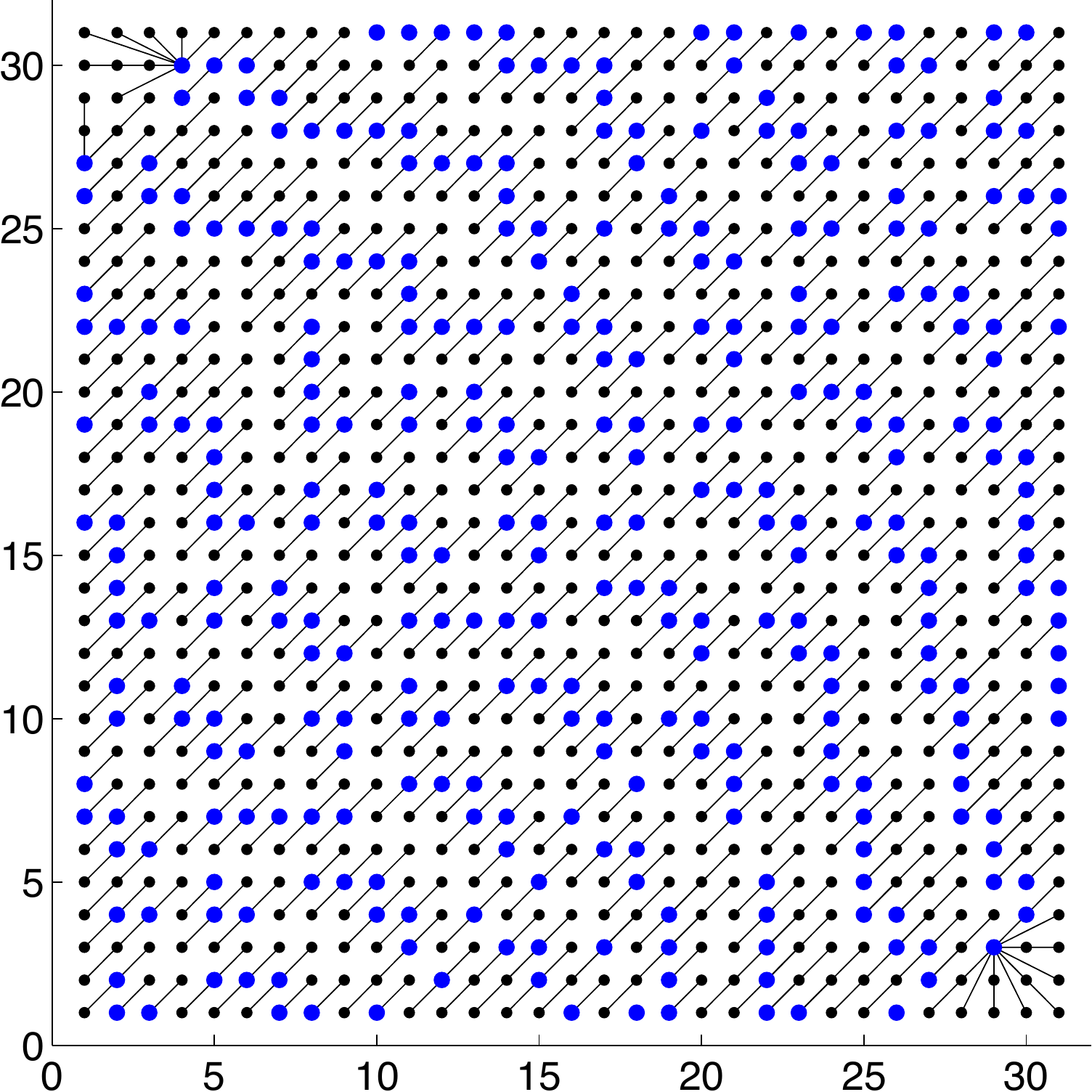}}
\hfill
  \subfigure[$\alpha = -\pi/4, \rho = .48, \rho_f = .60, \gamma_o = 1.497, \gamma_g = .445 $ \label{fig:aggs2mp4FD}]{\includegraphics[scale = 0.40]{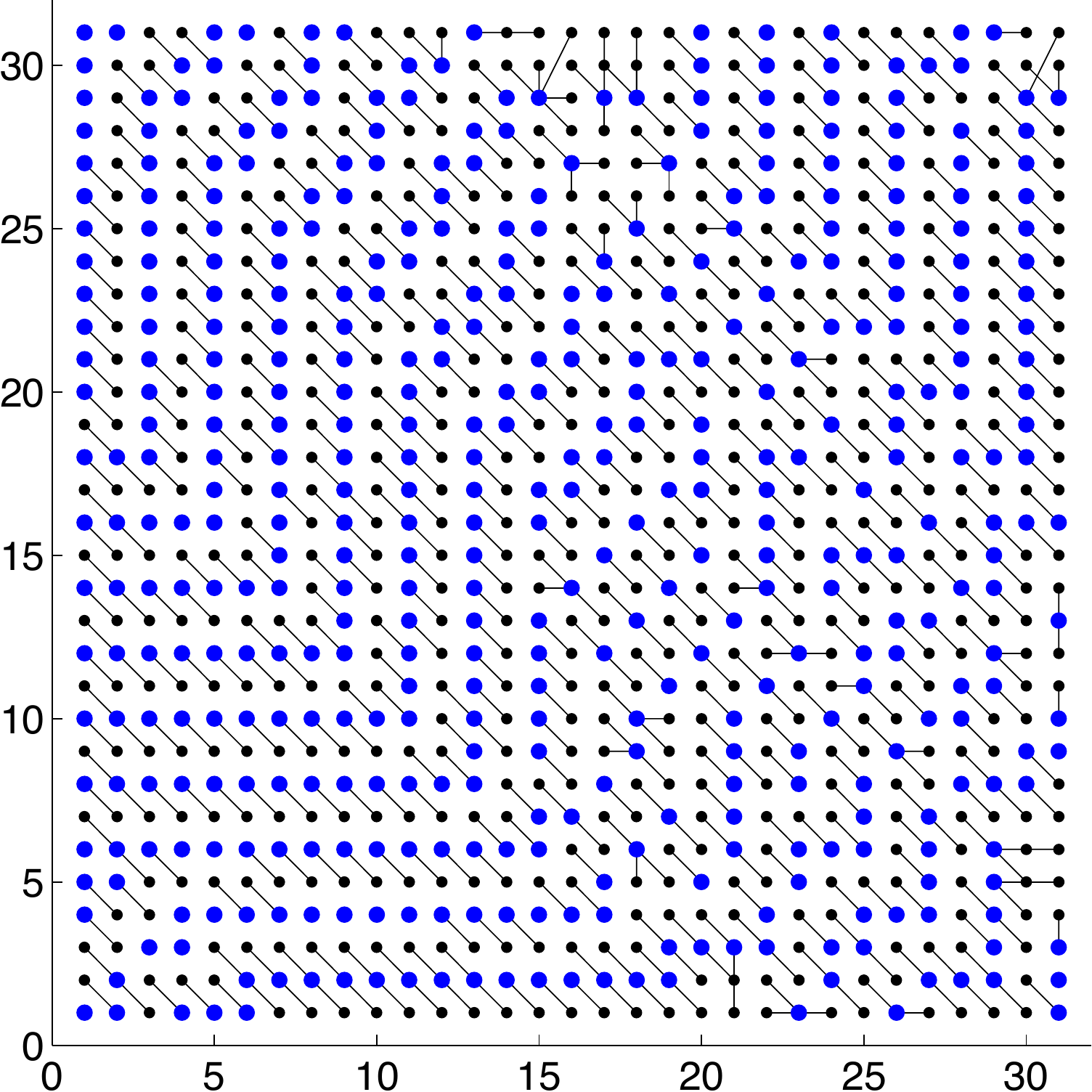}}
 \hfill
    \subfigure[$\alpha = \pi/8,   \rho = 0.21, \rho_f = .64, \gamma_o  =  1.481, \gamma_g = .431
 $\label{fig:aggs2p16FD}] {\includegraphics[scale = 0.40]{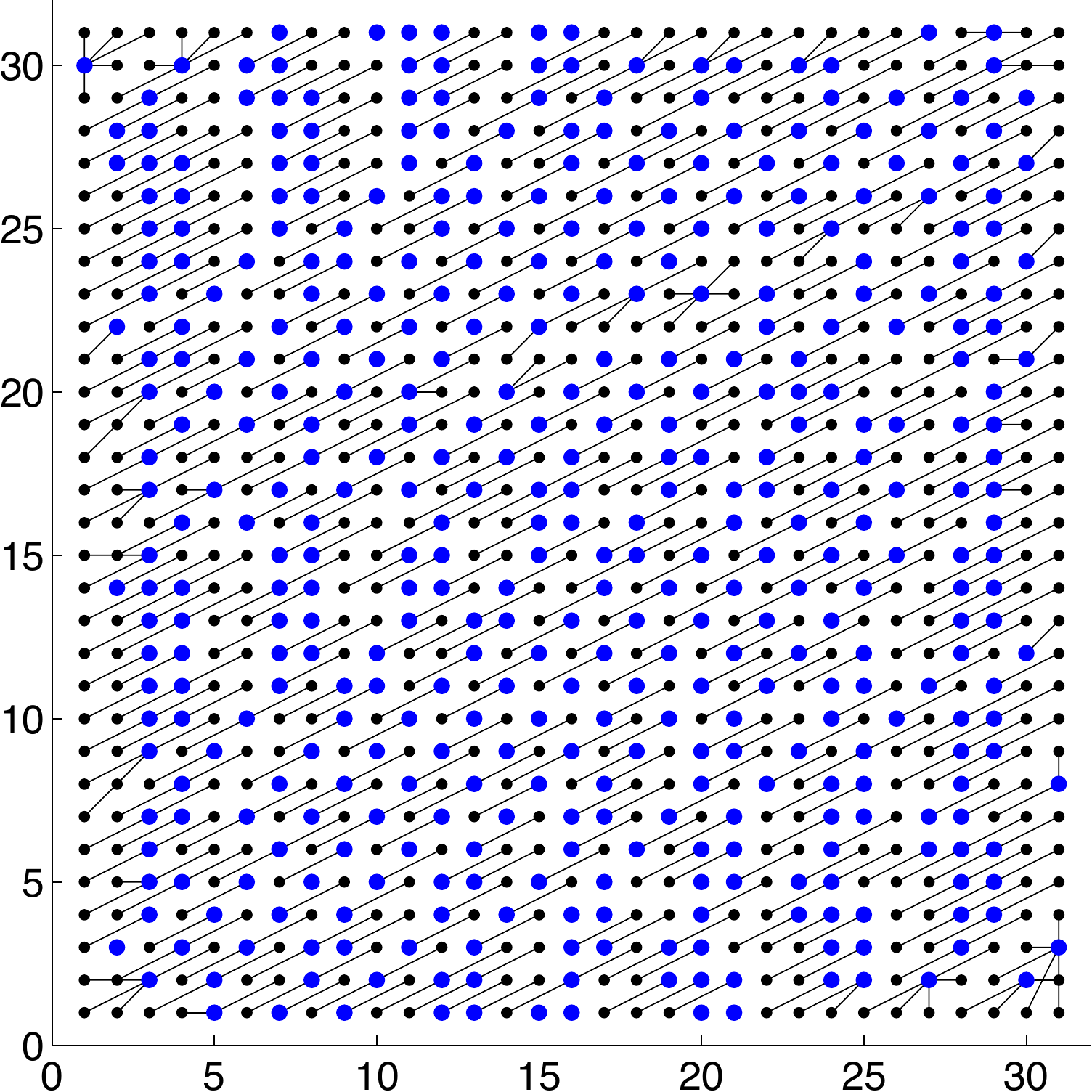}}
    \caption{Coarse grids and caliber $c=1$ interpolation patterns for the finite difference discretization with $h= 1/32$ for various choices of $\alpha$, using the graph of $A^2$, i.e., $d = 2$, and $d_{LS} = 4$, to define the strength matrix.  Here, the smaller circles are $F$-points and the larger circles are $C$-points \label{fig:aggs2FD}}
\end{figure}

\begin{figure}[htb!]
  \subfigure[$\alpha = 0, \rho = .31, \rho_f = .52, \gamma_o = 1.306, \gamma_g = .368 $\label{fig:s20agg}]{\includegraphics[scale = 0.40]{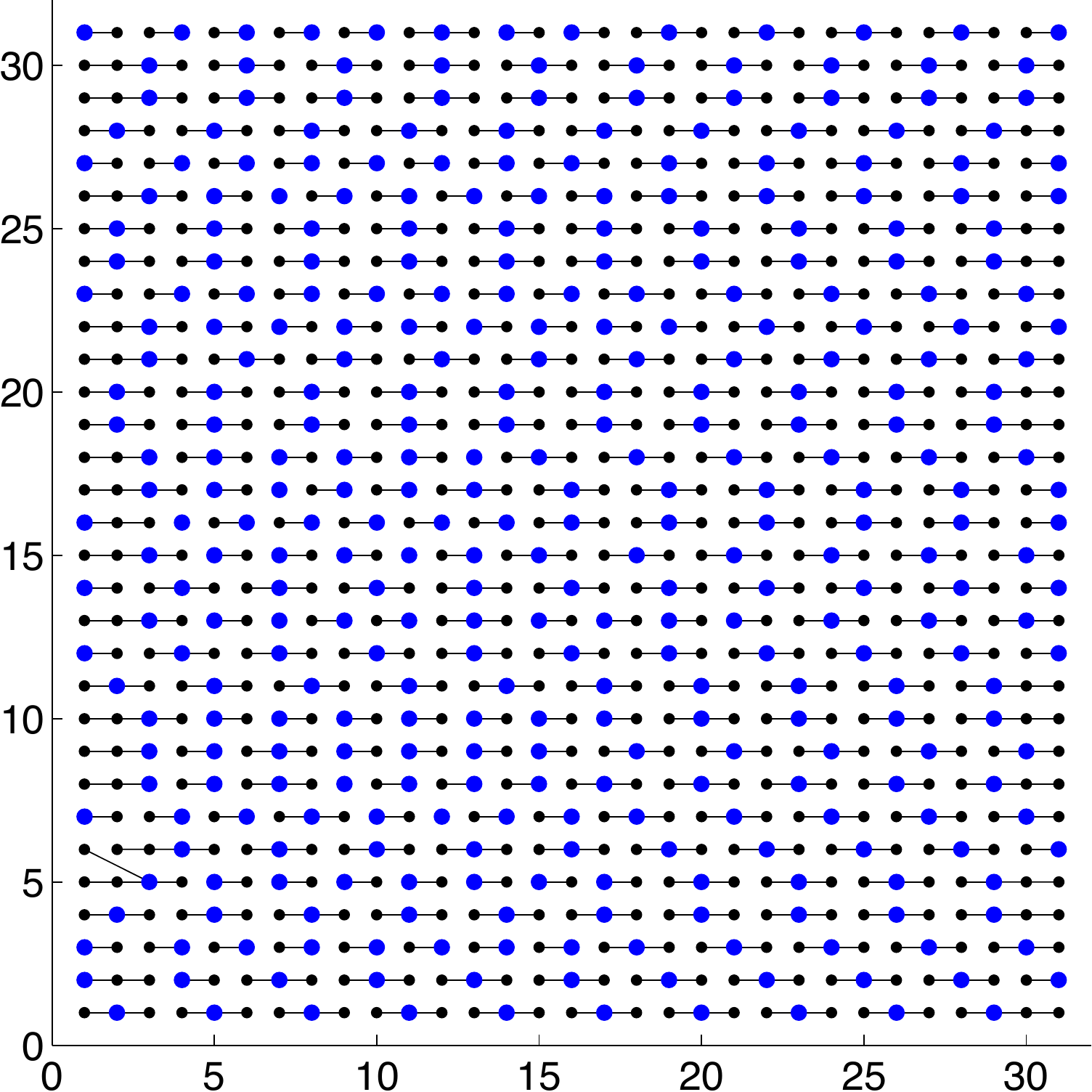}}  \hfill
  \subfigure[$\alpha = \pi/4, \rho = .25, \rho_f = .39, \gamma_o = 1.377, \gamma_g = .348 $\label{fig:s2p4agg}]{\includegraphics[scale = 0.40]{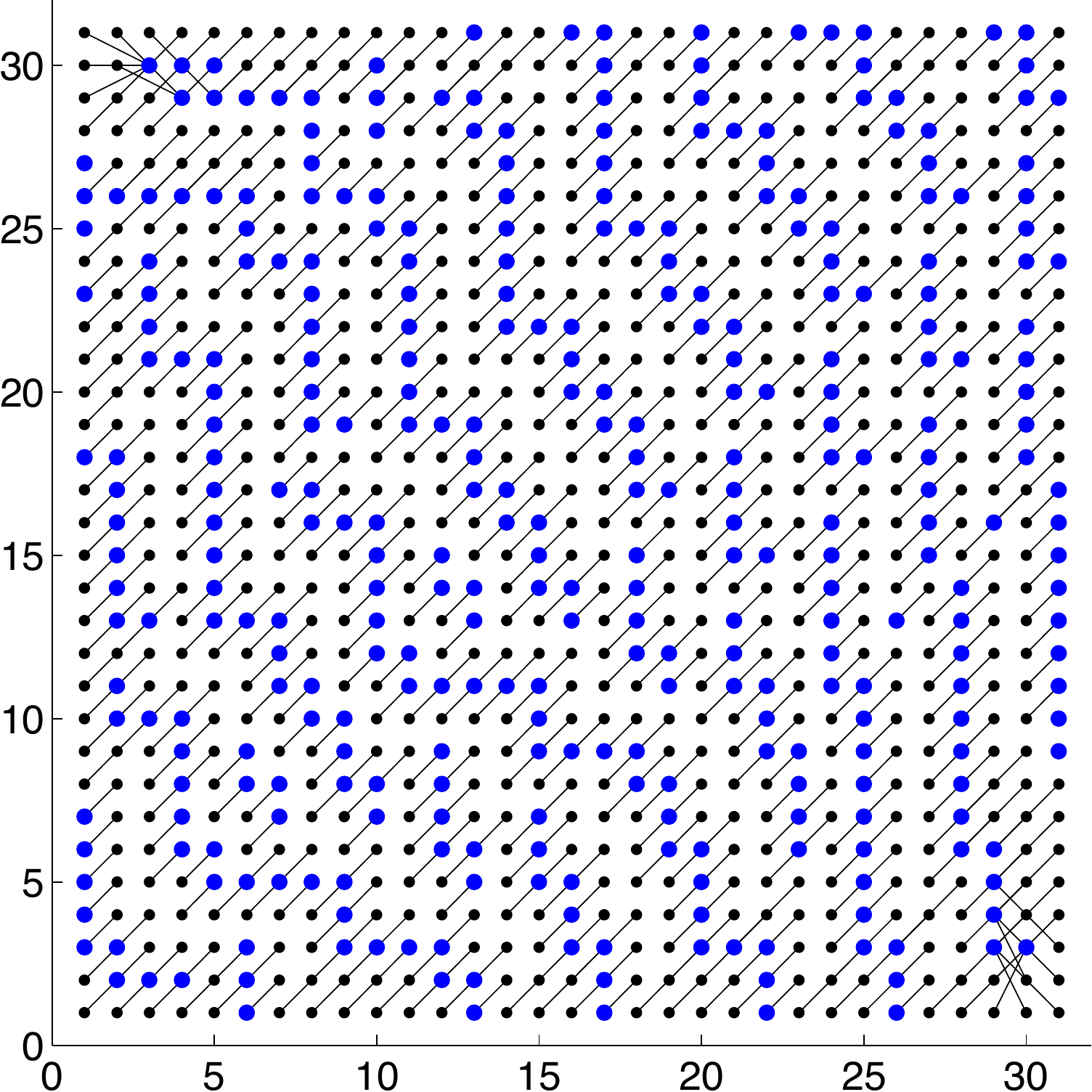}}
\hfill
 \subfigure[$\alpha = -\pi/4, \rho = .25, \rho_f = .65, \gamma_o = 1.376, \gamma_g = .348 $ \label{fig:s2mp4agg}]{\includegraphics[scale = 0.40]{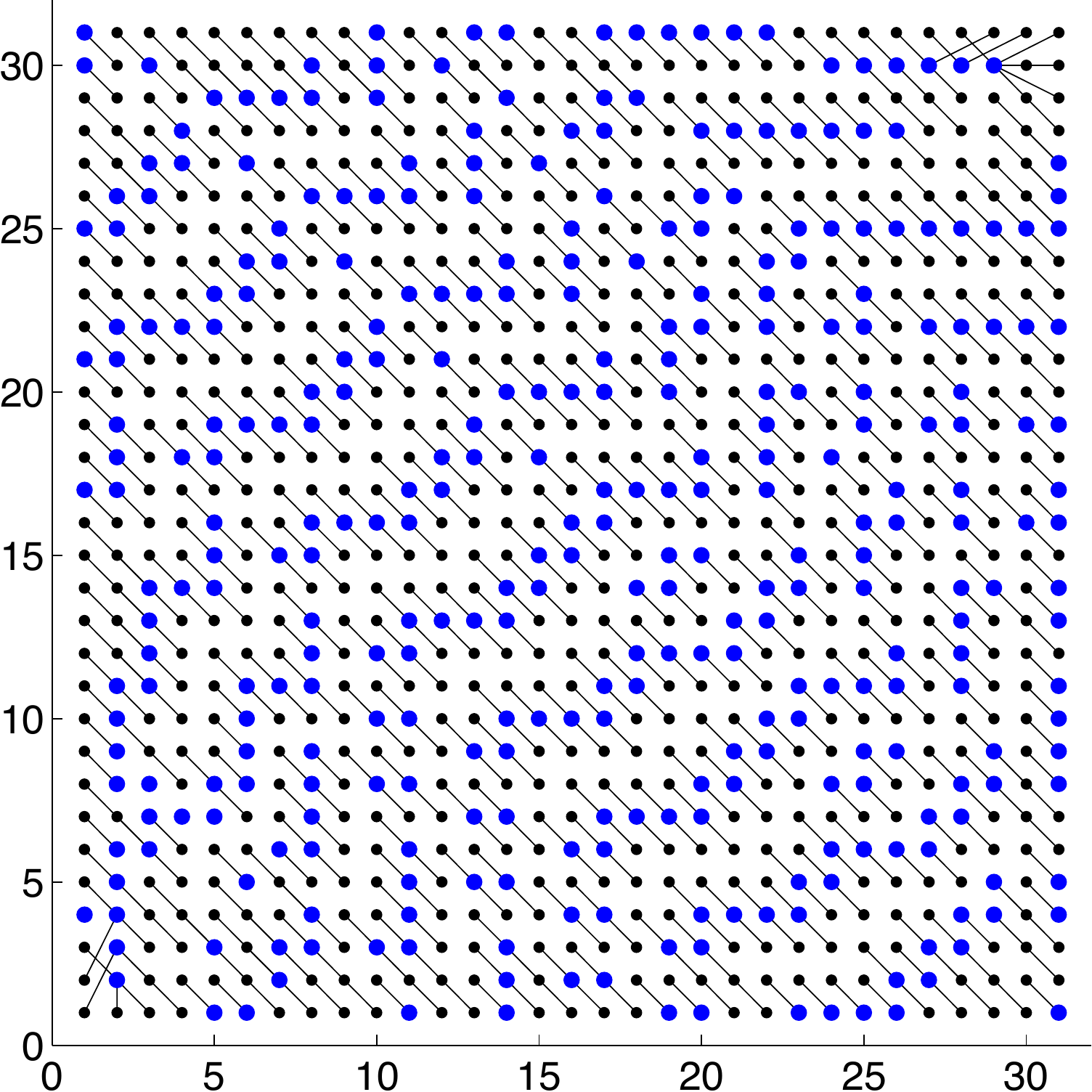}}
 \hfill
 \subfigure[$\alpha = \pi/8, \rho = .39, \rho_f = .64, \gamma_o = 1.5, \gamma_g = .418 $\label{fig:s2p16agg}]{\includegraphics[scale = 0.40]{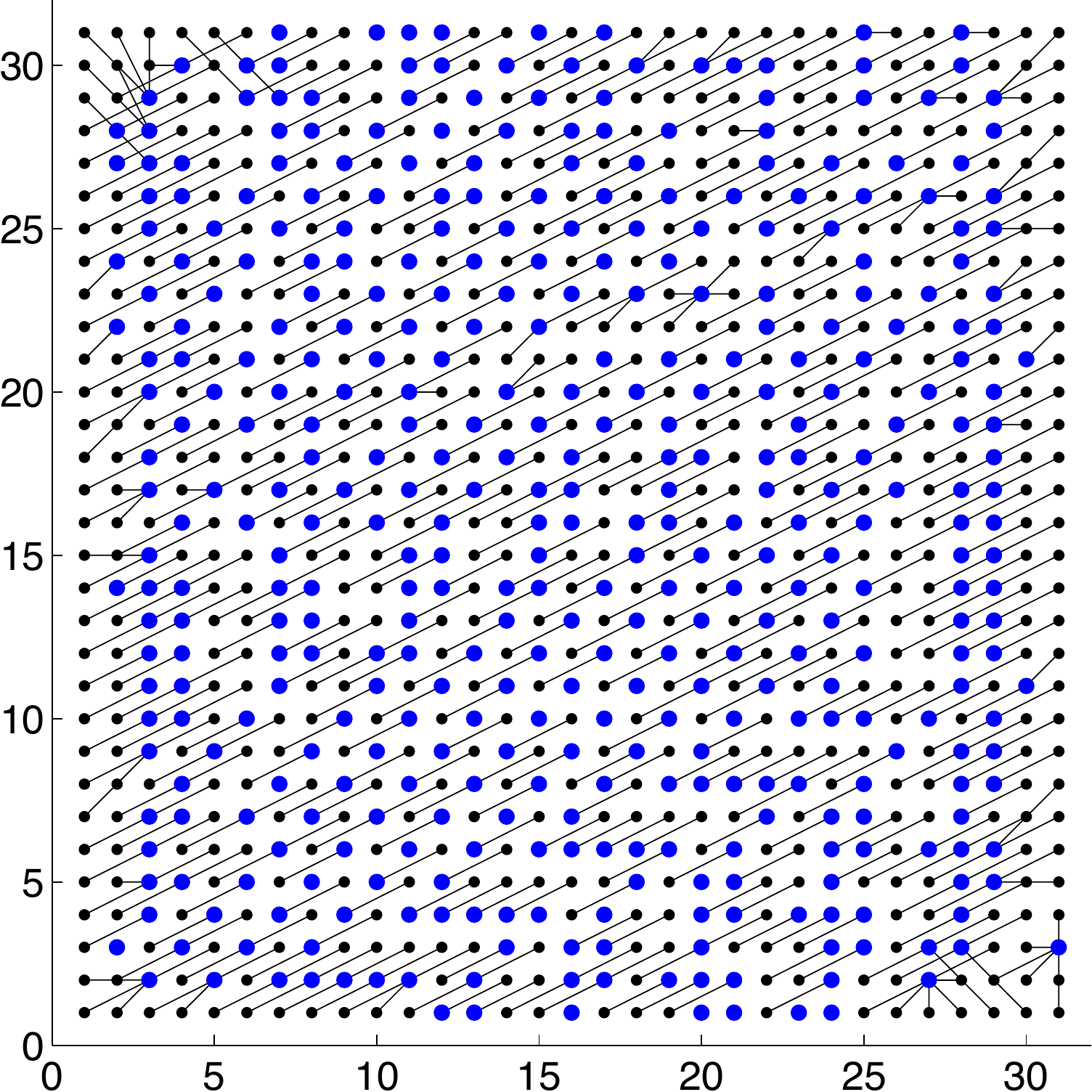}}  
    \caption{Coarse grids and caliber $c=1$ interpolation patterns for the bilinear finite element discretization with $h= 1/32$ for various choices of $\alpha$, using the graph of $A^2$, i.e., $d  = 2$ and $d_{LS} = 4$,  to define the strength matrix. \label{fig:s2agg}}
\end{figure}

\begin{figure}[htb!]
  \subfigure[$\alpha = 0, \rho = .01, \rho_f = .25, \gamma_o = 1.452, \gamma_g = .352 $\label{fig:s20}]{\includegraphics[scale = 0.40]{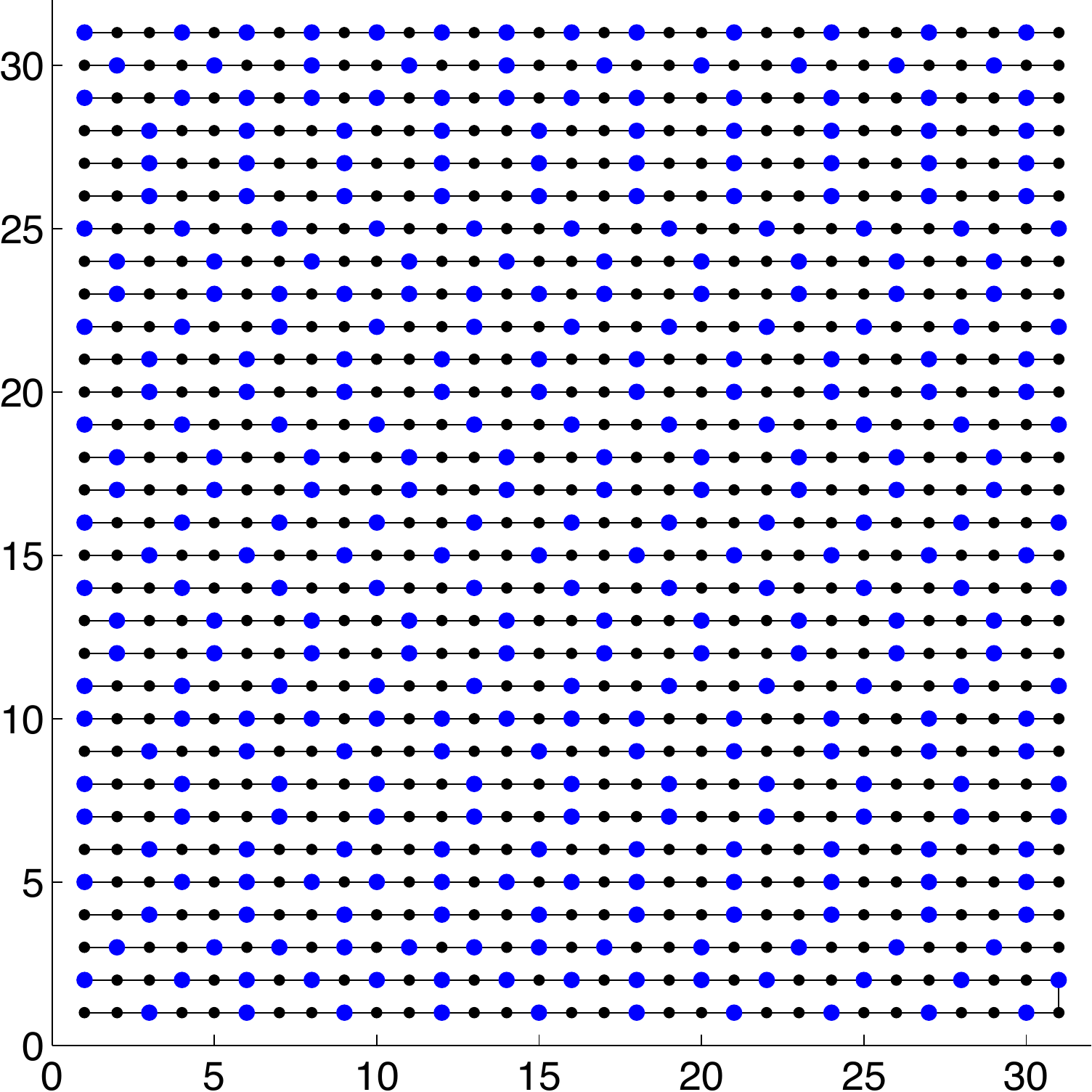}}  \hfill
  \subfigure[$\alpha = \pi/4, \rho = .03, \rho_f = .65, \gamma_o = 1.393, \gamma_g = 0.332 $\label{fig:s2p4}]{\includegraphics[scale = 0.40]{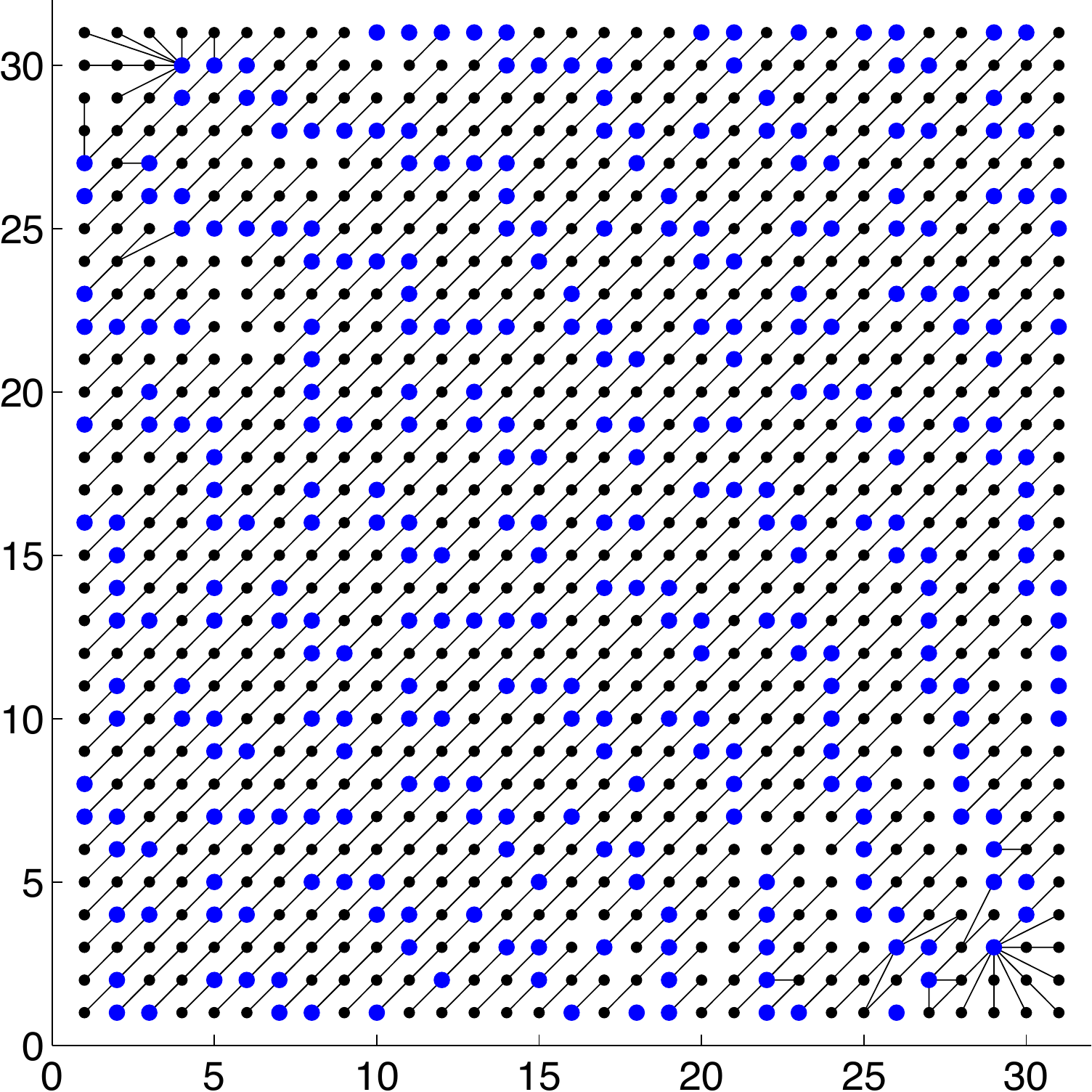}}
\hfill
  \subfigure[$\alpha = -\pi/4, \rho = .31, \rho_f = .60, \gamma_o = 1.772, \gamma_g = .445 $ \label{fig:s2mp4}]{\includegraphics[scale = 0.40]{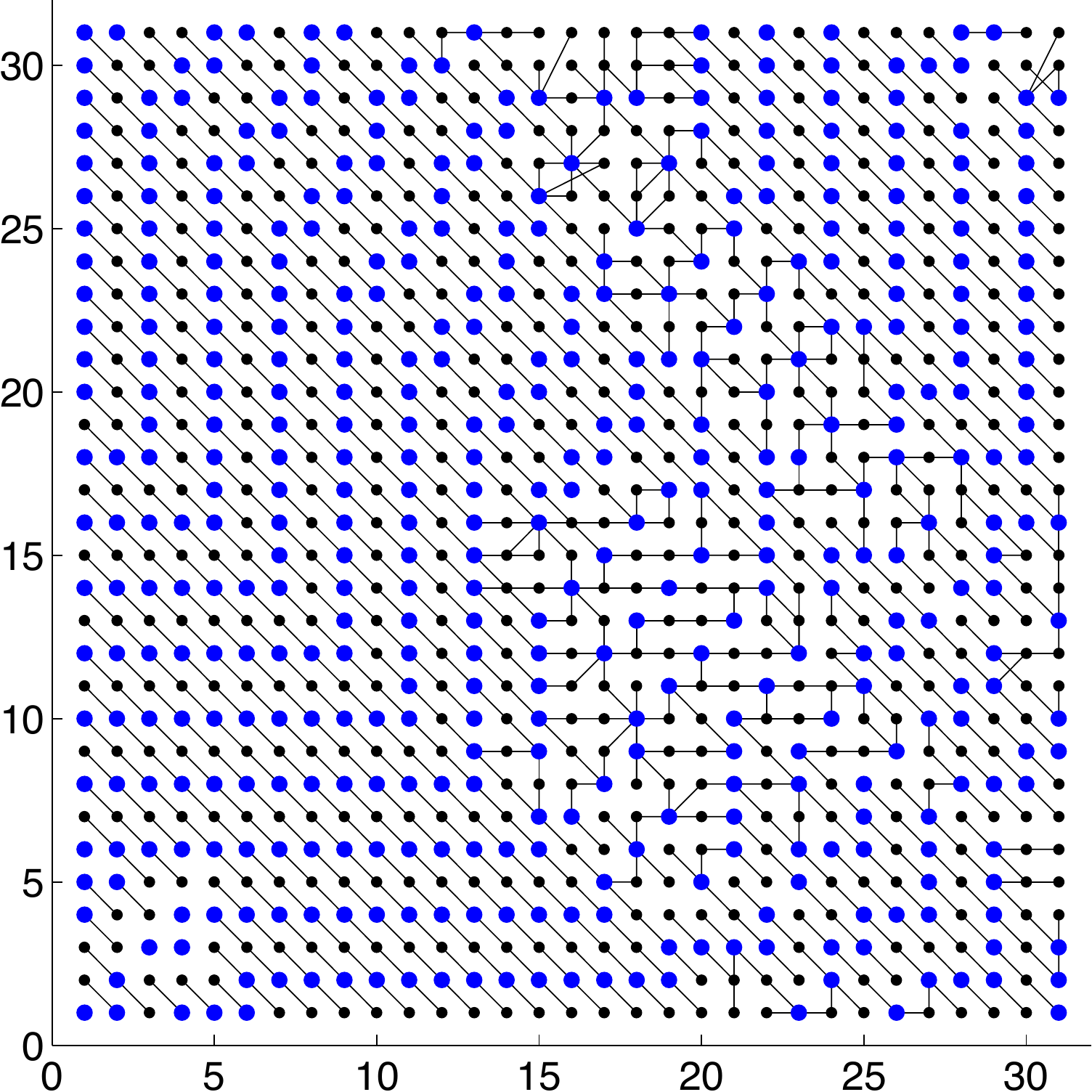}}
 \hfill
 \subfigure[$\alpha = \pi/8, \rho = .13, \rho_f = .64, \gamma_o = 1.788, \gamma_g = .431 $\label{fig:s2p16}]{\includegraphics[scale = 0.40]{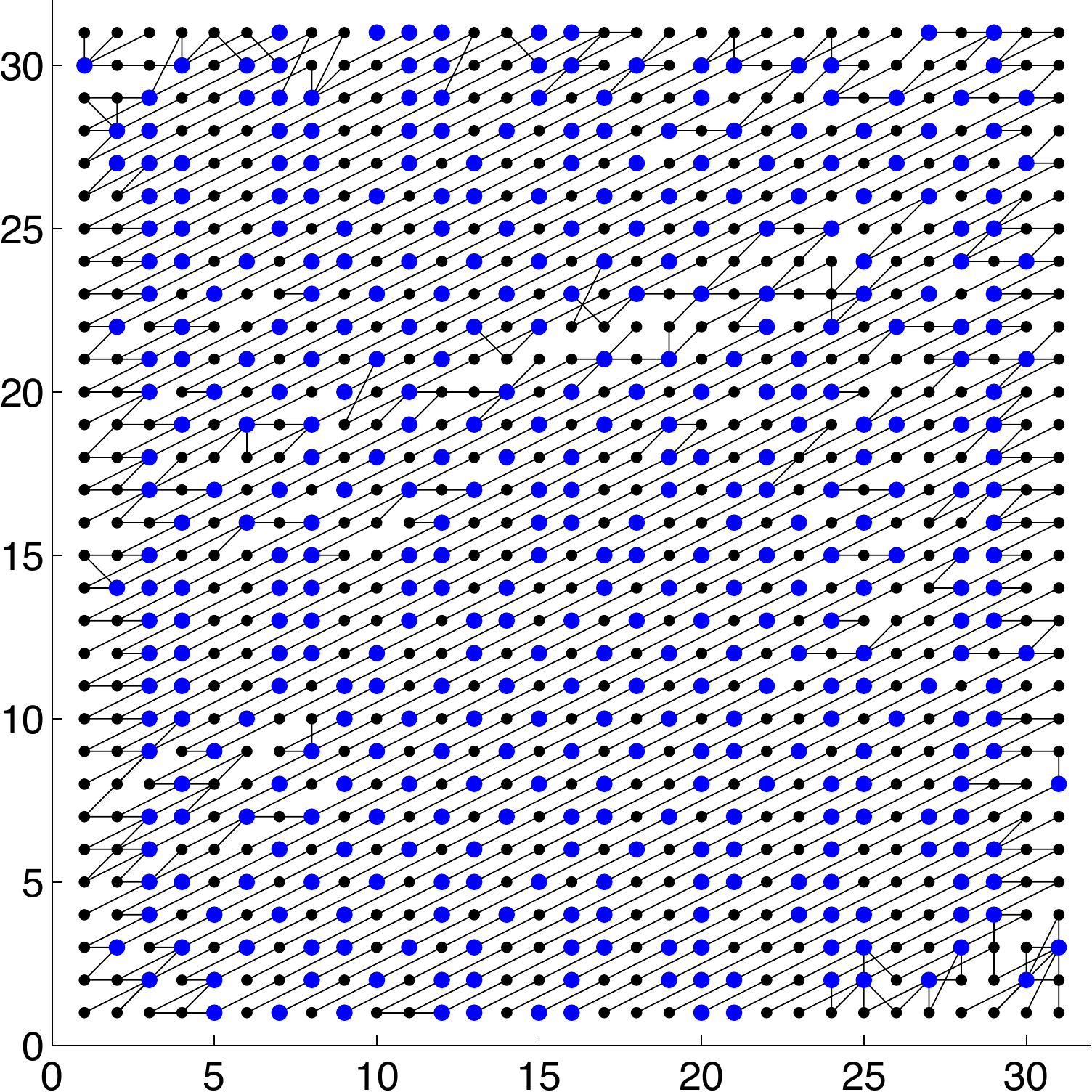}}  
    \caption{Coarse grids and caliber $c=2$  interpolation patterns for the finite difference discretization with $h= 1/32$ for various choices of $\alpha$, using the graph of $A^2$, i.e., $d  = 2$ and $d_{LS} = 4$,  to define the strength matrix. \label{fig:s2}}
\end{figure}

\begin{figure}[htb!]
  \subfigure[$\alpha = 0, \rho = .04, \rho_f = .52, \gamma_o = 1.396, \gamma_g = .368$ \label{fig:s20FE}]{\includegraphics[scale = 0.40]{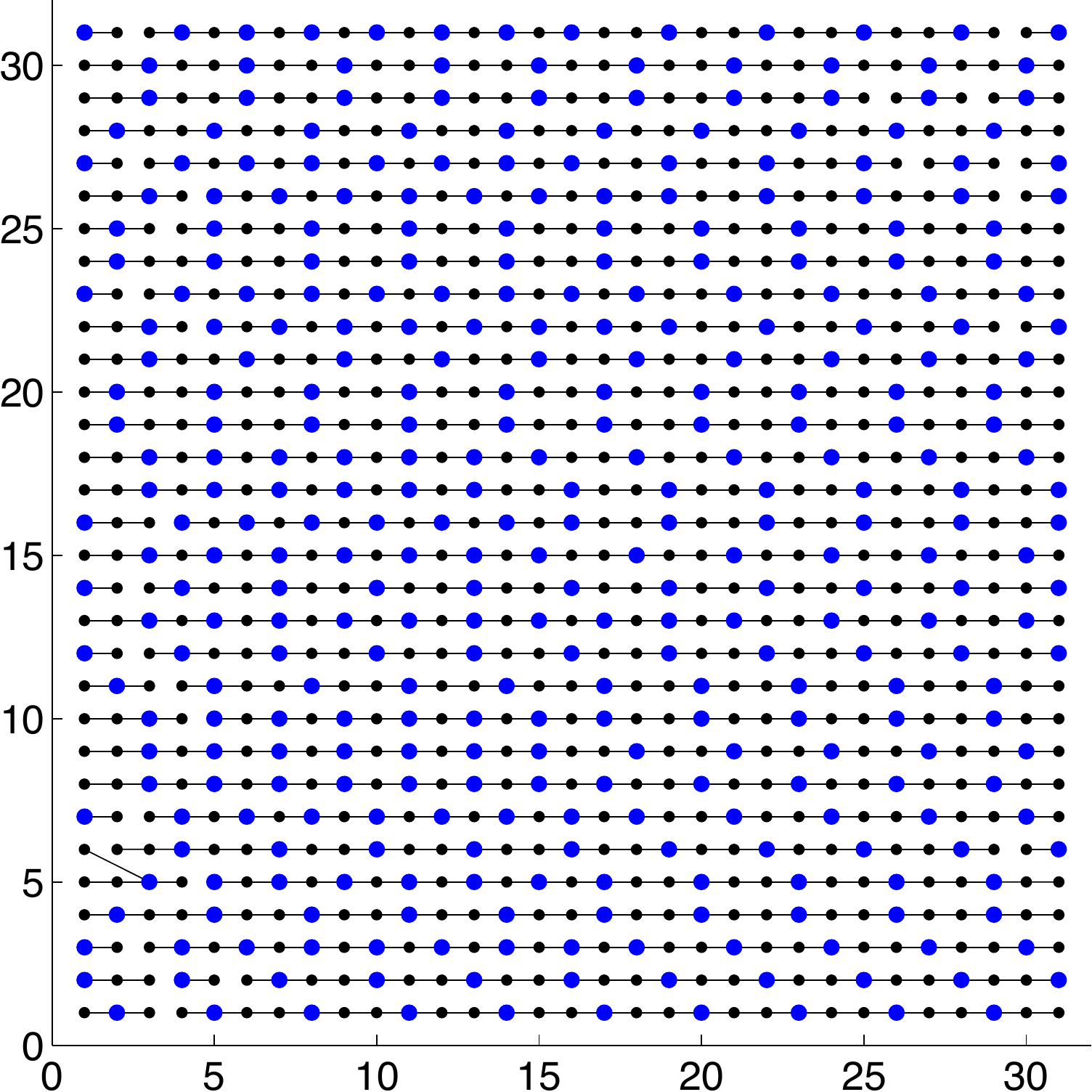}}
  \hfill
  \subfigure[$\alpha = \pi/4, \rho = .01, \rho_f = .39, \gamma_o =1.496, \gamma_g = .348$\label{fig:s2p4FE}]{\includegraphics[scale = 0.40]{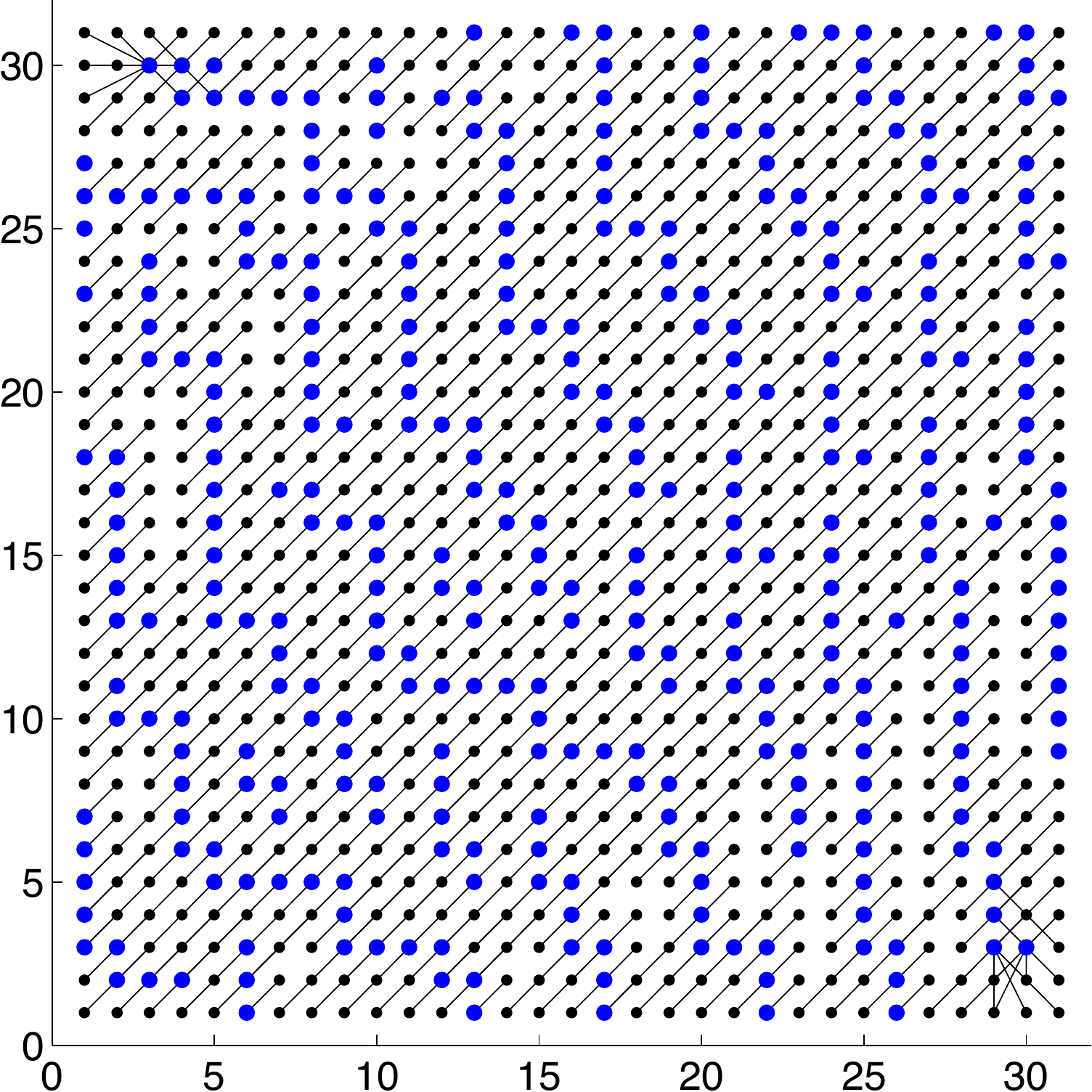}}
\hfill
  \subfigure[$\alpha = -\pi/4, \rho = .01, \rho_f = .64, \gamma_o = 1.498, \gamma_g = .348 $ \label{fig:s2mp4FE}]{\includegraphics[scale = 0.40]{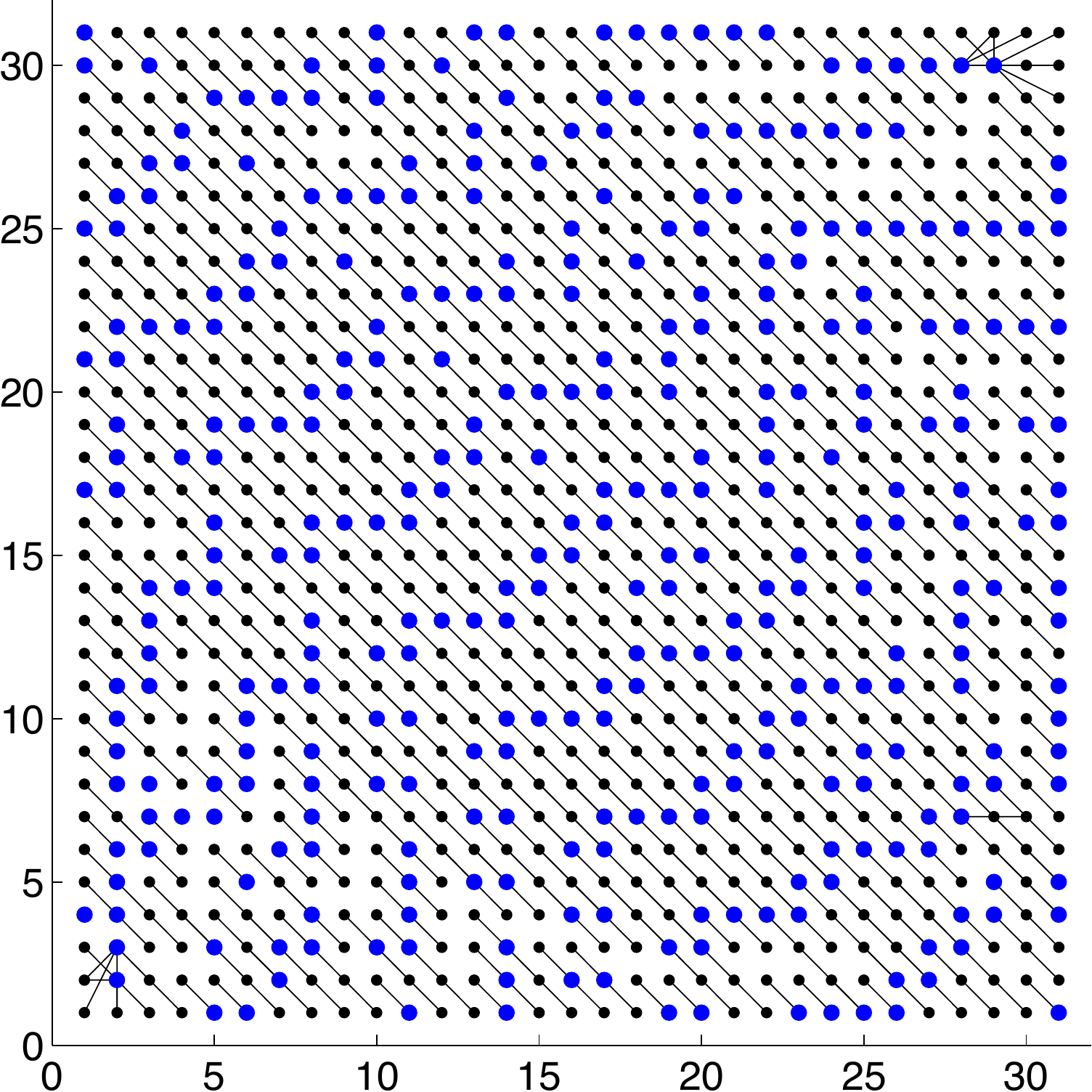}}
 \hfill
    \subfigure[$\alpha = \pi/8,   \rho = .22, \rho_f = .64, \gamma_o  = 1.774, \gamma_g = .418 $\label{fig:s2p16FE}]{\includegraphics[scale = 0.40]{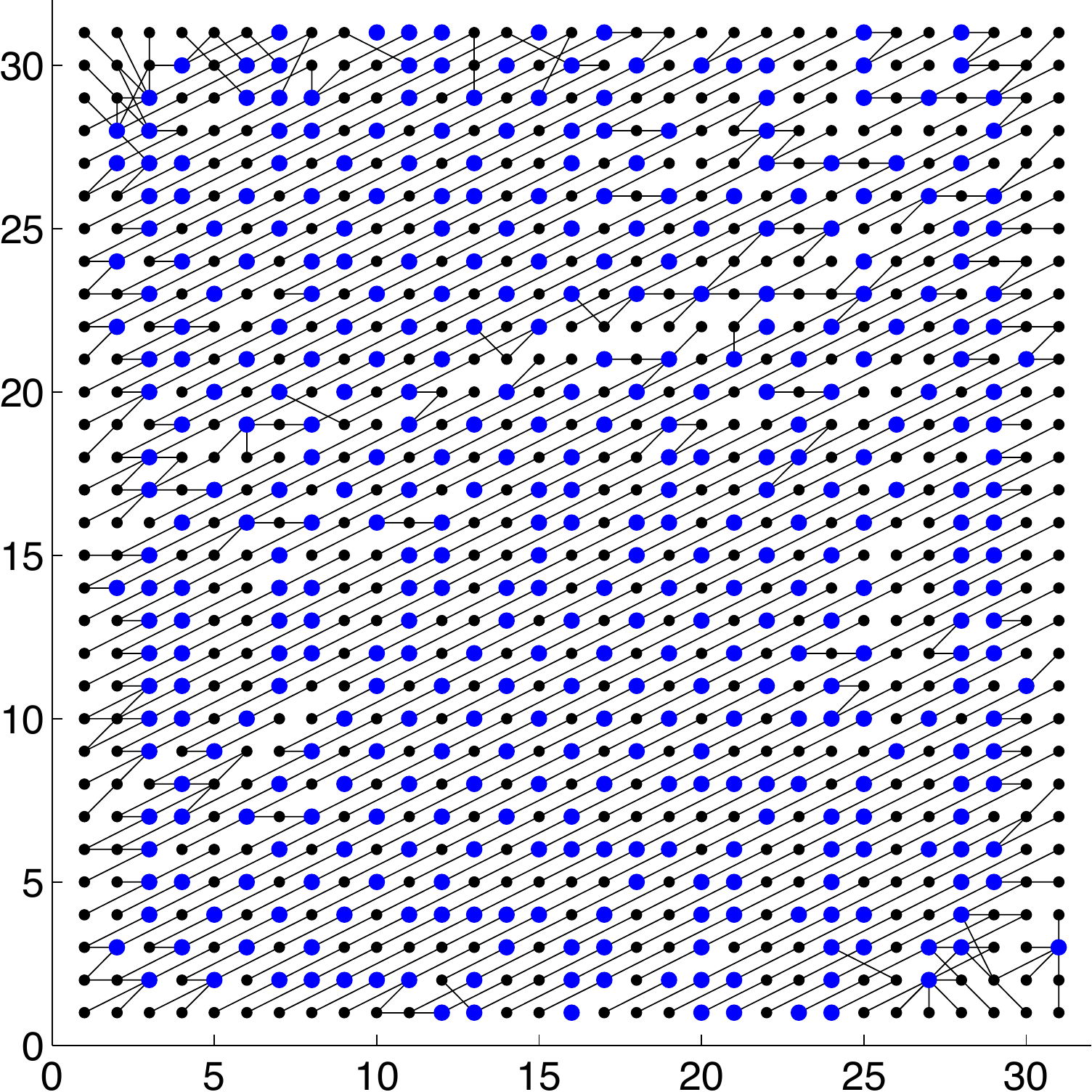}}
    \caption{Coarse grids and caliber $c=2$ interpolation patterns for the bilinear finite element discretization with $h= 1/32$ for various choices of $\alpha$, using the graph of $A$, i.e., $d = 2$ and $d_{LS} = 4$, to define the strength matrix.  \label{fig:s2FE}}
\end{figure}

In this section, we illustrate the choice of coarse grids and interpolation patterns the BAMG setup algorithm constructs when applied to the finite difference and bilinear finite element approximations of the anisotropic test problem for $\epsilon = 10^{-4}$ with various choices of the anisotropic angle $\alpha$ in \eqref{eq:a}.  The tests of BAMG setup we consider are for interpolation calibers $c=1,2$ and varying values of the search depth $d$. 
In the plots in Figures~\ref{fig:aggs1FD}-\ref{fig:s2FE}, the black lines depict the interpolation pattern for each $\finevar$-point (denoted by the smaller  circles) from its neighboring $\coarsevar$-points (denoted by the larger  circles).  
Generally, we observe that the coarsening and interpolation pattern follow the anisotropy when the choice of the discretization allow it.  For example, in Figures~\ref{fig:aggs1FD} and \ref{fig:aggs1FE} the coarsening perfectly follows 
the direction of the anisotropy for both discretizations when using a search depth $d=1$ and setting $\alpha = 0$ and $\alpha = \pi/4$. 
On the other hand, we see that this is not true for the case for $\alpha = -\pi/4$ and the seven-point finite difference approximation and for $\alpha =  \pi/8$ for both discretizations.  For these angles the direction of strong coupling is not able to capture the direction of the anisotropy for the given discretization of the problem when setting 
the search depth to $d=1$ since the matrix entries in these directions are zero.   Hence, for these problems it is not possible to form a strength matrix from immediate algebraic neighbors which 
produces a coarse grid that allows the interpolation pattern to exactly
follow the anisotropy.
These observations in fact lead naturally to 
the use of the graph $\mathcal{G}_d$ of $A^d$, $d>1$, to form the strength
matrix.  Overall, we see that the algebraic
distance strength measure leads to least squares interpolation that follows the general direction of anisotropy to the extent that the value of he search depth $d$ allows it.  

\textcolor{black}{Another interesting observation here is that for the $d=2$ tests reported in Figures~\ref{fig:aggs2FD}-~\ref{fig:s2FE}, we see that by using the graph of $A^2$ to form the algebraic-distance-based strength graph, the coloring algorithm is now able to coarsen in the direction of anisotropy for the finite difference discretization and $\alpha = -\pi/4$.  We note that the use of larger choices of the search depth in defining the strength matrix have the additional benefit --  it allows more aggressive coarsening while at the same time maintaining the characteristic directions induced by more general anisotropic directions.  The ability of the setup to produce aggressive coarsening is seen in the tests for grid-aligned anisotropy when comparing the grids obtained for the search depths $d=1$ to $d=2$.}

Another difficult choice of $\alpha$ occurs for the anisotropy direction $\pi/8$ for which
a longer-range interpolation and an extended search depth for coarse-grid candidates is required to properly capture the direction of the anisotropy.
For this example, the reference angles made by coarse-fine connections are much closer to $\pi/8$ when using $d=2$ (see Figures~\ref{fig:aggs2FD}-\ref{fig:s2FE}) than they are when using $d=1$ (see Figures~\ref{fig:aggs1FD}-\ref{fig:s1FE}) in the setup for both the finite difference and finite element discretizations.   

To conclude the section we consider the performance of the algorithm 
for two different choices of the interpolation caliber $c=1,2$.  Here, the algorithm chooses the same coarse grids independent 
of the choice of caliber, but the sparsity structure and also the values of the interpolation operators 
change.  
Generally, we have seen that using a larger value of $c$ improves the convergence rates of the two-level methods, $\rho$, in all cases, however, 
the complexity of solving the coarse-level system also increases.  

\subsection{The coarse-level operator}
Another interesting  deliverable of the proposed BAMG setup algorithm, in particular of its implementation of the compatible relaxation and the algebraic distances, is the pattern of the resulting coarse-grid stencil. Discretizations involving the finite difference discretization 
(\ref{fd_seven})  favor $\alpha = \pi/4$ and with the same argument result in the worst possible discretization for  $\alpha = -\pi/4$ for which the use of the upper-left and lower-right  grid-point neighbors in the discretization of $\partial_{xy}$ would be appropriate.     
We consider both cases $\alpha = \pi/4$ and $\alpha = -\pi/4$ for the seven-point finite difference discretization given in (\ref{fd_seven}). For both cases, we assume $\epsilon = 10^{-10}$,  $d = 2$, and $d_{LS} = 4$.

\begin{figure}
  \subfigure[Coarse-grid equation pattern for $\alpha = \pi/4$.\label{fig:good}]{\includegraphics[width=.4\textwidth]{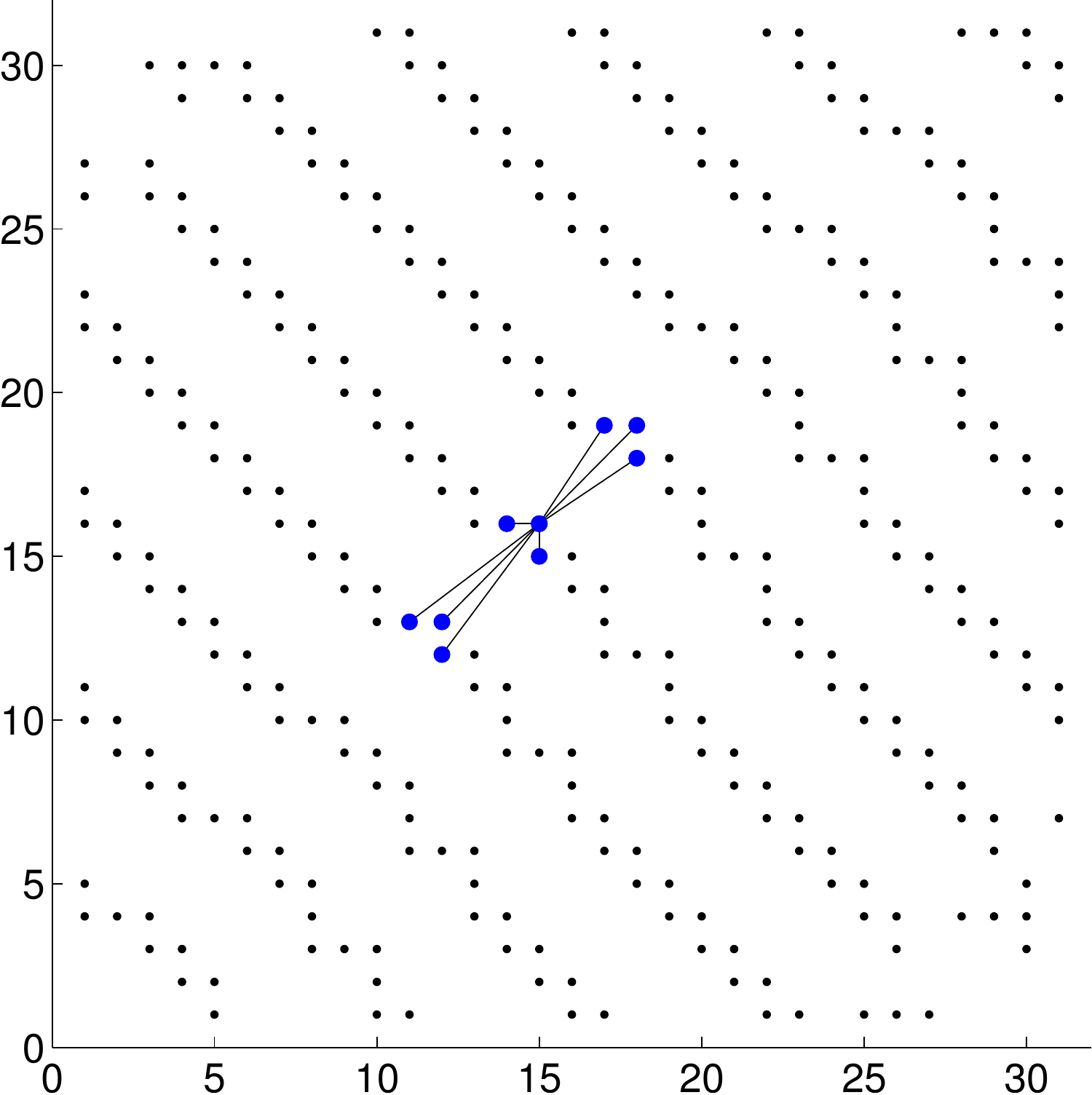}}\hfill
  \subfigure[Coarse-grid equation pattern for $\alpha = -\pi/4$ \label{fig:bad}]{\includegraphics[width=.4\textwidth]{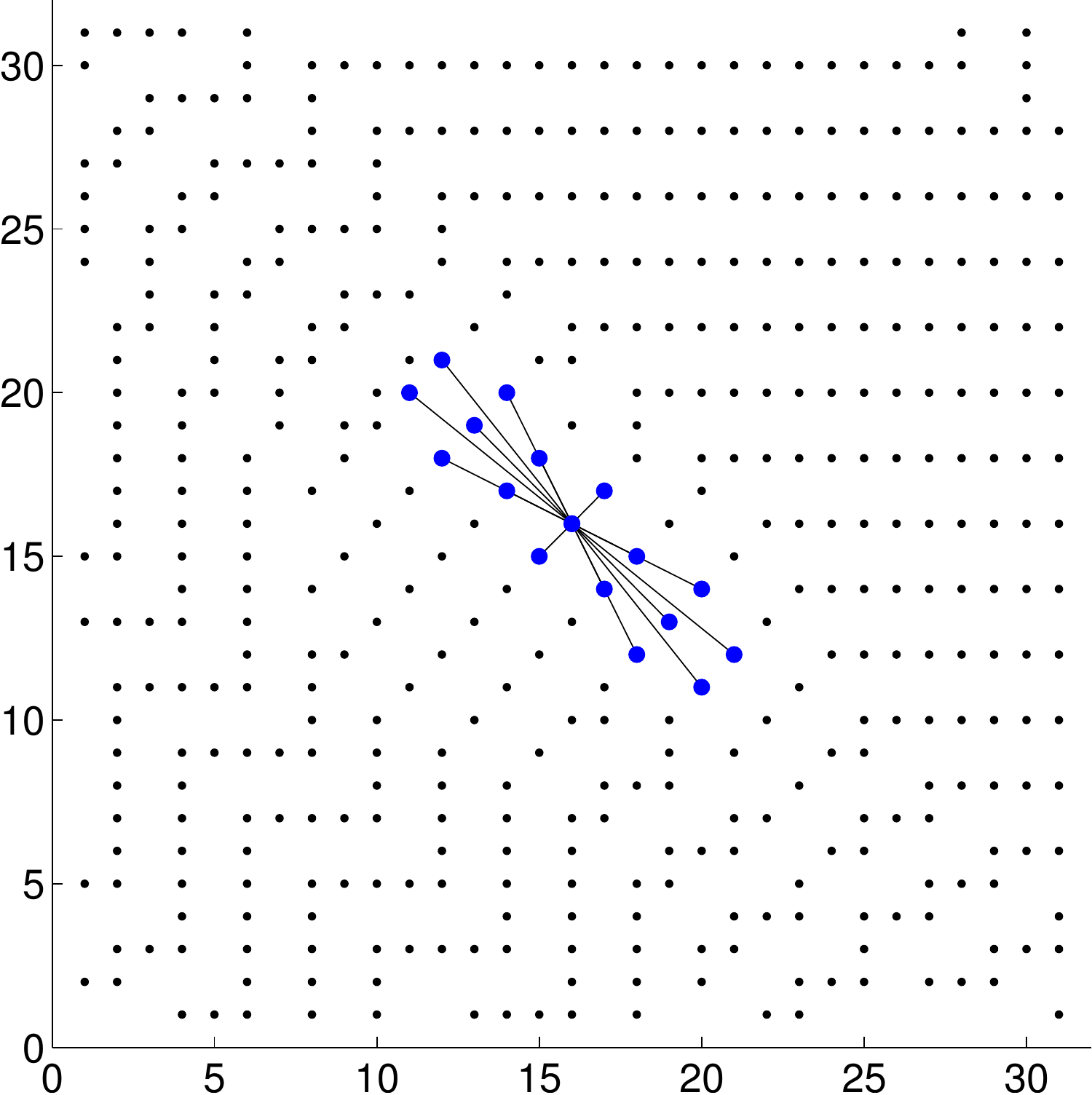}}
  \caption{Non-zero pattern of one of the stencils of the coarse-grid operator $A_c$, centered at $i \in C$ and connected with 
  $j \in C$ such that $(A_c)_{ij} \neq 0$, for for $h=1/32$ and $\alpha = \pi/4$. The smaller dots in the graph are all other coarse-grid variables.    
  \label{fig:pis}}
\end{figure}

 We first confirm that for $\alpha = \pi/4$, the coarse-grid operator $A_c = P^T A P$ preserves the intrinsic strength of connections inherited from the fine-grid operator $A$.  A typical example of the stencil of $A_c$  is given in Figure \ref{fig:good}.  Here, the details of configurations of each stencil depend on the coarsening pattern in the neighborhood of the considered coarse-grid equation.
 
 The results for the more challenging $\alpha = -\pi/4$ case are provided in Figure \ref{fig:bad}.  Here, we observe that although the discretization on the fine grid does not follow the anisotropy whatsoever, the non-zero pattern of the coarse-grid operator correctly aligns with the direction of anisotropy.  This result demonstrates the ability of the algorithm to overcome, if needed, the disadvantage of a poorly chosen fine-grid discretization and regain a more favorable discretization on the first coarse grid.  Further, the results for the $\alpha =  \pi/4$  indicate, that all consecutive coarse grids (though not employed in the two-level algorithm) are likely  to maintain a similar favorable discretization that too accurately reflects the anisotropy.

Coefficients of the coarse-grid stencils,  presented in
Figure \ref{fig:pis}, are given next. Here $S_{cg}^{+}$ corresponds to
$\alpha = \pi/4$ (entries denoted by $*$  are  negligible,  with absolute values below $10^{-11}$)

\[ \displaystyle{
S_{cg}^{+} = \tiny{
\begin{pmatrix} 
     &   & & &  * &  -0.17\\
             &        & & &  & * &\\
           &  & * & 0.33 & &    \\
     &  &  & * & &  \\
        * & -0.17 &  &  & & \\
             & * & &  &  & 
                     \end{pmatrix}, } }
    \]
    and $S_{cg}^{-} $ corresponds to  $\alpha = -\pi/4$:
        \[
                   S_{cg}^{-} =
\left(
\tiny{ \begin{array}{*{12}{r}} 
            & -0.11 & & & & & & & &  &  &  \\
            -0.12 &  &  &  \phantom{-}0.23 & & & &  &   &  \\
            & & \phantom{-}1.37 & & & & & & &  &   &  \\
            &  0.23 & &  & -2.94  & & & &  &   & \\
           &  & &   -2.91 &   & & 1.06 & &   &  \\
            &  & &    &   &    \phantom{-}6.36 &   & &  &   &  \\
            &  & &    &  1.06  &  &   & -2.90 &        &   \\              
            &  & &    &  &  &  -2.93 &   &  &    0.23   &   \\           
                &  & &    &  &  &   &   &   \phantom{-} 1.38  &     &  \\
                       &  & &    &  &  &   &    0.23 &  & &  -0.12    \\ 
           &  & &    &  &  &   &    &  &  -0.12  &   
 \end{array}} 
 \right).
\]
\noindent
In both stencils, all  entries are rounded to the nearest hundredth.  

The distribution of the stencils' coefficients further illustrates the ability of the algebraic distance based 
strength of connections measure to choose the 
correct coarse-grid points for problems with anisotropic coefficients.  In both cases, the non-zero sparsity pattern and dominant coefficients of the resulting coarse-grid operators follow the direction of anisotropy.   Clearly, larger search depths,
both for the CR and the LS procedures, yield an additional fill-in of the resulting coarse-level operator, with this effect becoming more profound in a multilevel setting. In such cases, it becomes especially important to choose accurate coarse grids based on a rich set of test functions. Strategies for  computing the latter are discussed in the following sections.

\subsection{Two-level convergence}
Here, we present experiments with tests of the proposed AMG setup algorithm applied to \eqref{eq:diff} for various choices of the anisotropy angle, $\alpha$, the anisotropy coefficient, 
$\epsilon$, and the mesh spacing $h = 1/N$.  \textcolor{black}{In Tables~\ref{tab:2gridFD}-~\ref{tab:2gridFE} we  present test results  for the bootstrap setup with calibers $c=1$ and $c=2$ interpolation for the finite difference and finite element discretizations.  We note that multilevel results are given in the next section.}
In the tables, the asymptotic convergence rates, $\rho$, of the two-grid solver produced by the BAMG setup algorithm  are reported, along with the corresponding coarsening factors $\gamma_g$ and  operator complexity ratios
 $\gamma_o$.  Here, for $c=1$, we observe a strong dependence of the computed convergence rates and grid and operator complexities on the problem parameters, $\epsilon$, 
 $\alpha$, and $h$.  For $c=2$, this dependence is less pronounced with only a slight dependence
 on $h$ that is restricted mostly to the non-grid aligned cases. The exception is $\alpha = 0$ and $\epsilon = .1$, where we see a slight increase in $\rho$ as the problem size is increased from $N=64$ to $N=128$.   Moreover, in all cases, the convergence rates and complexities are uniformly bounded with respect to $\epsilon$ and $\alpha$ for fixed $h$.   {\color{black}These results are promising when considering that all tests were performed with the same strength parameter $\theta_{ad} = .5$.  In fact, all  parameters in the setup algorithm were fixed, illustrating that the individual components of the BAMG setup are robust for the targeted anisotropic problems, even those leading to non $M$-matrix systems as in the $\alpha = -\pi/4$ and $\alpha = \pi/8$ cases.}
 \textcolor{black}{Further, we note that the bootstrap setup with caliber $c=2$ interpolation handles the isotropic case when $\alpha = 0$ and $\epsilon =1$ with similar efficiency, producing a two-grid method with convergence rate $\rho = .28$ and complexities $\gamma_g = .25$ and $\gamma_o = 1.6$ for $h=1/32,1/64,1/128$.}

\begin{table}[!ht]
  \begin{center}
 
 {\small \begin{tabular}{|c|c|c|c|c|c|c|c|}
      \hline
      $\alpha$\\
      \hline
      \hline
       $0$  \\
       \hline
      $\pi/4$\\
      \hline
       $-\pi/4$\\
       \hline
       $\pi/8$ \\ 
         \hline
         \end{tabular}
           \begin{tabular}{|c|c|c|c|c|c|c|c|}
    \hline
      $c=1$ and  $\epsilon =  .1$ \\ 
       \hline      \hline
         $   .23 (.33,1.3) / .36 (.36,1.4)  / .45 (.38,1.4)$ \\
       \hline 
     $     .14 (.37,1.4) / .37 (.33,1.4)  / .49 (.35,1.4) $  \\
       \hline
    $      .36 (.37,1.4) /  .60 (.33,1.4)  / .77 (.32,1.4) $ \\
        \hline
     $   .14 (.28,1.3) / .42 (.26,1.3)  / .67 (.25,1.3)$ \\
      \hline

              \end{tabular}
 \begin{tabular}{|c|c|c|c|c|c|c|c|}
    \hline
      $c=2$ and $\epsilon =  .1$\\ 
         \hline      \hline
     $    .04 (.33,1.6) / .13 (.36,1.5)  / .20 (.38,1.5)$ \\
       \hline 
      $    .01 (.37,1.5) / .04 (.33,1.5)  / .05 (.35,1.5) $  \\
       \hline
        $  .07 (.37,1.6) /  .27 (.33,1.5)  / .31 (.32,1.5) $ \\
        \hline
        $ .01 (.28,1.4) /  .12 (.26,1.4) / .15 (.25,1.4) $ \\
      \hline         \end{tabular}
\bigskip

 \begin{tabular}{|c|c|c|c|c|c|c|c|}
      \hline
       $\alpha$\\
      \hline
      \hline
       $0$  \\
       \hline
      $\pi/4$\\
      \hline
       $-\pi/4$\\
       \hline
       $\pi/8$ \\ 
        \hline
         \end{tabular}
    \begin{tabular}{|c|c|c|c|c|c|c|c|}
    \hline
  $c=1$ and  $\epsilon = .0001$ \\
        \hline      \hline
     $    .18 (.35,1.3) /  .41 (.38,1.4)  / .61 (.40,1.4)$  \\
       \hline 
$        .14 (.33,1.3) /   .33 (.36,1.4)  /  .49 (.37,1.4) $  \\
       \hline
  $         .48 (.45,1.5) / .69 (.42,1.5) / .84 (.40,1.5) $   \\
        \hline
   $     .21 (.43,1.5) /  .68 (.40,1.5)   / .89 (.38,1.5)  $ \\
      \hline
         \end{tabular} 
  \begin{tabular}{|c|c|c|c|c|c|c|c|}
    \hline
  $c=2$ and  $\epsilon = .0001$ \\
        \hline      \hline
     $    .01 (.35,1.5) /  .05 (.38,1.6)  / .05 (.40,1.5)$  \\
       \hline 
$        .03 (.33,1.4) /   .05 (.36,1.5)  /  .06 (.37,1.6) $  \\
       \hline
  $         .31 (.45,1.8) / .38 (.42,1.7) / .42 (.40,1.7) $   \\
        \hline
   $     .13 (.43,1.8) /  .35 (.40,1.7)   / .43 (.38,1.7)  $ \\
      \hline
         \end{tabular} 
        \bigskip

 \begin{tabular}{|c|c|c|c|c|c|c|c|}
      \hline
      $\alpha$\\
      \hline
      \hline
       $0$  \\
       \hline
      $\pi/4$\\
      \hline
       $-\pi/4$\\
       \hline
       $\pi/8$ \\ 
        \hline
         \end{tabular}
    \begin{tabular}{|c|c|c|c|c|c|c|c|}
    \hline
      $c=1$ and  $\epsilon = 0$\\
        \hline      \hline
        $    .18 (.35,1.3) /  .62 (.38,1.4)  / .66 (.40,1.4)$  \\
       \hline 
$        .10 (.33,1.3) /   .33 (.36,1.3)  /  .42 (.38,1.4) $  \\
       \hline
  $         .48 (.45,1.5) / .70 (.38,1.5) / .87 (.36,1.6) $   \\
        \hline
   $     .38 (.43,1.5) /  .76 (.35,1.6)   / .91 (.38,1.7)  $ \\
      \hline
         \end{tabular}
  \begin{tabular}{|c|c|c|c|c|c|c|c|}
    \hline
        $c=2$ and $\epsilon = 0$\\
        \hline      \hline
     $    .01 (.35,1.3) / .08 (.38,1.5)  / .09 (.40,1.5)$\\
       \hline 
   $    .04 (.33,1.4) / .05 (.36,1.5) / .06 (.38,1.6)  $  \\
       \hline
 $        .31 (.45,1.8) / .37 (.38,1.7) / .41 (.36,1.8)$ \\
        \hline
 $        .12 (.43,1.8) /  .35 (.35,1.7)  / .40 (.38,1.8)$  \\
      \hline
         \end{tabular}

}
        \bigskip

   \caption[Two-grid LS interpolation for FD Laplacian]{Approximate asymptotic convergence rates of the two-grid solver applied to the seven point FD Anisotropic Laplace problem with Dirichlet boundary conditions for various choices of $\alpha$, $\epsilon$ and $h$.  Here, the proposed setup algorithm is applied with search depth for the CR coarsening algorithm set as $d=2$, and with the search depth for caliber $c=1$ and $c=2$ least squares interpolation set as $d_{LS} = 4$.  The reported results correspond to convergence rates $\rho$, and, in parenthesis, coarsening factors, $\gamma_g$, and operator complexity factors, $\gamma_0$, computed for $h = 1/32$,  $h=1/64$, $h=1/128$. }\label{tab:2gridFD}  
  \end{center}
\end{table} 

\begin{table}[!ht]
  \begin{center}
 
 {\small \begin{tabular}{|c|c|c|c|c|c|c|c|}
      \hline
      $\alpha$\\
      \hline
      \hline
       $0$  \\
       \hline
      $\pi/4$\\
      \hline
       $-\pi/4$\\
       \hline
       $\pi/8$ \\ 
         \hline
         \end{tabular}
           \begin{tabular}{|c|c|c|c|c|c|c|c|}
    \hline
      $c=1$ and  $\epsilon =  .1$ \\ 
       \hline      \hline
         $   .15 (.35,1.3) / .33 (.34,1.3)  / .45 (.34,1.4)$ \\
       \hline 
     $     .18 (.36,1.4) / .43 (.32,1.4)  / .49 (.35,1.4) $  \\
       \hline
    $      .17 (.32,1.4) /  .42 (.30,1.3)  / .77 (.32,1.4) $ \\
        \hline
     $   .36 (.29,1.3) / .49 (.26,1.3)  / .67 (.25,1.3)$ \\
      \hline

              \end{tabular}
 \begin{tabular}{|c|c|c|c|c|c|c|c|}
    \hline
      $c=2$ and $\epsilon =  .1$   \\
         \hline      \hline
     $    .05 (.35,1.4) / .18 (.34,1.4)  / .21 (.34,1.5)$ \\
       \hline 
      $    .04 (.36,1.5) / .09 (.32,1.4)  / .11 (.35,1.5) $  \\
       \hline
        $  .02 (.32,1.5) /  .19 (.30,1.4)  / .24 (.32,1.6) $ \\
        \hline
        $ .22 (.29,1.8) /  .26 (.26,1.3) / .33 (.25,1.5) $ \\
      \hline         \end{tabular}
\bigskip

 \begin{tabular}{|c|c|c|c|c|c|c|c|}
      \hline
       $\alpha$\\
      \hline
      \hline
       $0$  \\
       \hline
      $\pi/4$\\
      \hline
       $-\pi/4$\\
       \hline
       $\pi/8$ \\ 
        \hline
         \end{tabular}
    \begin{tabular}{|c|c|c|c|c|c|c|c|}
    \hline
  $c=1$ and  $\epsilon = .0001$ \\
        \hline      \hline
     $    .31 (.37,1.3) /  .55 (.39,1.3)  / .61 (.40,1.4)$  \\
       \hline 
$        .25 (.35,1.4) /   .46 (.36,1.4)  /  .58 (.37,1.4) $  \\
       \hline
  $         .25 (.35,1.4) / .52 (.36,1.4) / .84 (.37,1.5) $   \\
        \hline
   $     .39 (.42,1.5) /  .68 (.39,1.5)   / .89 (.38,1.5)  $ \\
      \hline
         \end{tabular}          
  \begin{tabular}{|c|c|c|c|c|c|c|c|}
    \hline
  $c=2$ and  $\epsilon = .0001$ \\
        \hline      \hline
     $    .04 (.37,1.4) /  .05 (.39,1.6)  / .05 (.40,1.5)$  \\
       \hline 
$        .01 (.35,1.5) /   .19 (.36,1.5)  /  .22 (.37,1.5) $  \\
       \hline
  $         .01 (.35,1.5) / .20 (.36,1.6) / .25 (.37,1.6) $   \\
        \hline
   $     .22 (.42,1.8) /  .29 (.39,1.6)   / .36 (.38,1.6)  $ \\
      \hline
         \end{tabular} 
        \bigskip

 \begin{tabular}{|c|c|c|c|c|c|c|c|}
      \hline
      $\alpha$\\
      \hline
      \hline
       $0$  \\
       \hline
      $\pi/4$\\
      \hline
       $-\pi/4$\\
       \hline
       $\pi/8$ \\ 
        \hline
         \end{tabular}
    \begin{tabular}{|c|c|c|c|c|c|c|c|}
    \hline
      $c=1$ and  $\epsilon = 0$\\
        \hline      \hline
        $    .32 (.37,1.3) /  .60 (.39,1.3)  / .66 (.39,1.4)$  \\
       \hline 
$        .25 (.33,1.4) /   .46 (.36,1.4)  /  .72 (.38,1.4) $  \\
       \hline
  $         .25 (.35,1.5) / .51 (.36,1.4) / .74 (.37,1.4) $   \\
        \hline
   $     .40 (.42,1.5) /  .68 (.39,1.5)   / .79 (.38,1.5)  $ \\
      \hline
         \end{tabular}
  \begin{tabular}{|c|c|c|c|c|c|c|c|}
    \hline
        $c=2$ and $\epsilon = 0$\\
        \hline      \hline
     $    .05 (.37,1.4) / .10 (.39,1.4)  / .13 (.39,1.4)$\\
       \hline 
   $    .04 (.33,1.4) / .21 (.36,1.6) / .23 (.38,1.7)  $  \\
       \hline
 $        .01 (.35,1.5) / .20 (.36,1.6) / .26 (.37,1.6)$ \\
        \hline
 $        .22 (.42,1.8) /  .33 (.39,1.6)  / .45 (.38,1.6)$  \\
      \hline
         \end{tabular}

}
        \bigskip

   \caption[Two-grid LS interpolation for FD Laplacian]{Approximate asymptotic convergence rates of the two-grid solver applied to the nine point bilinear FE Anisotropic Laplace problem with Dirichlet boundary conditions for various choices of $\alpha$, $\epsilon$ and $h$.  Here, the proposed setup algorithm is applied with search depth for the CR coarsening algorithm set as $d=2$, and with the search depth for caliber $c=1$ and $c=2$ least squares interpolation set as $d_{LS} = 4$.  The reported results correspond to convergence rates $\rho$, and, in parenthesis, coarsening factors, $\gamma_g$, and operator complexity factors, $\gamma_0$, computed for $n = 32^2$/  $h = 1/32$,  $h=1/64$, $h=1/128$. }\label{tab:2gridFE}  
  \end{center}
\end{table} 
\subsection{Multilevel convergence}  Next, we consider the performance of a multilevel BAMG algorithm.  We report results of a nonlinear Algebraic Multilevel Iteration (AMLI) $W$-cycle preconditioner constructed by using recursively the proposed bootstrap setup algorithm applied to the same anisotropic test problems discretized using finite differences and bilinear finite elements.  
{\color{black}Though numerical results for the stand-alone $V$ cycle solver and preconditioner are not reported, we note that the convergence rates of both approaches deteriorate for increasing problem sizes and strength of anisotropy in the non-aligned cases.  This observation in fact motivated our
use of the AMLI $W$-cycle preconditioner.}
For details of the nonlinear AMLI solve cycle we refer to~\cite{Vassilevski_2005}.
Here, we limit the numerical tests to caliber $c=2$ interpolation since this choice produced the best results in 
the tests of the two-grid method given the previous section.   
Figure \ref{fig:boot:setupcycle} provides an schematic outline of the bootstrap V- and W-cycle setup algorithms.  
A main ingredient of the bootstrap setup is its use of a multilevel generalized eigensolver to compute the bootstrapped test vectors once an initial multigrid hierarchy has been constructed. The goal is to enrich the set of TV by using approximations of the lowest eigenmodes of the finest level matrix $A_0= A$, obtained by computing eigenvectors on the coarsest level $L$, and then transferring them, with some additional local processing, to the finest grid. The main ideas of the multilevel eigensolver are as follows.

Given the initial Galerkin operators $A_{0}, A_{1}, \ldots, A_{L}$ on each level
and the corresponding interpolation operators $P_{l+1}^{l}, l =
0,\ldots,L-1$, define the composite interpolation operators as
$P_{l} = P_{1}^{0}\cdot \ldots \cdot P_{l}^{l-1},\ l = 1, \ldots, L$.  Then, 
for any given vector $x_{l} \in \mathbb{C}^{n_l}$ we have
$
  \innerprod[A_{l}]{x_{l}}{x_{l}} = \innerprod[A]{P_{l}x_{l}}{P_{l}x_{l}}.
$ 
Furthermore, defining $T_{l} = P_{l}^{H}P_{l}$ we obtain
\begin{equation*}
  \frac{\innerprod[A_{l}]{x_{l}}{x_{l}}}{\innerprod[T_{l}]{x_{l}}{x_{l}}} = \frac{\innerprod[A]{P_{l}x_{l}}{P_{l}x_{l}}}{\innerprod{P_{l}x_{l}}{P_{l}x_{l}}}.
\end{equation*} 
This observation in turn implies that 
on any level $l$, given a vector $v^{(l)} \in \mathbb{C}^{n_l}$ and $\lambda^{(l)}\in
  \mathbb{C}$ such that
$
  A_l v^{(l)} = \lambda^{(l)}T_l v^{(l)},
$
we have the Raleigh quotient
\begin{equation}\label{eq:evalapprox}
  \mbox{rq}(P_l v^{(l)}) :=
  \frac{\innerprod[A]{P_{l}v^{(l)}}{P_{l}v^{(l)}}}{\innerprod{P_{l}v^{(l)}}{P_{l}v^{(l)}}}
  = \lambda^{(l)}.
\end{equation}
In this way, the eigenvectors and eigenvalues
of the operators in the multigrid hierarchy on all levels are related directly
to the eigenvectors and eigenvalues of the finest-grid operator $A$.  
Note that the eigenvalue approximations in~\eqref{eq:evalapprox} are continuously updated within the algorithm so that the overall approach resembles an inverse Rayleigh-Quotient
iteration found in eigenvalue computations (cf.~\cite{JWilkinson_1965}).  
For additional details of the algorithm and its implementation we refer to the paper \cite{BAMG2010}.  

{\color{black}The cost per iteration of a single $V(\mu_1,\mu_2)$ setup cycle
with $\mu_1$ pre- and $\mu_2$ post-smoothing steps can be roughly estimated in terms of the cost of a single fine grid relaxation  
step, i.e., one work unit, which we define as the number of non zero entries in the fine-level matrix denoted by $nnz_0(A_0)$.  
Letting, as before, the operator complexity ratio be the total number of nonzero entries of the matrices on all levels of the multilevel 
hierarchy:
$$\gamma_o =  \dfrac{\sum_l nnz(A_l)}{nnz(A)},$$
the cost of a single $V(\mu_1,\mu_2)$ is then roughly given by
\begin{equation}\label{eq:costW}
(\mu_1 + \mu_2) \times \gamma_o  \times k \quad \text{work units,}
\end{equation}
where $k$ denotes the number of test vectors computed in the BAMG setup.
For example, the cost of a single $V(4,4)$ cycle which computes eight relaxed vectors $v \in \mathcal{V}^r$ and 
eight eigenvector approximations $v \in \mathcal{V}^e$
requires roughly
$64 \gamma_o$ work units.  We note that here we are neglecting the cost of 
recomputing the LS interpolation operators and coarse-level operators on all levels.  
}

\begin{figure}
\begin{center}
     \tikzstyle{greenpoint}=[circle,inner sep=0pt,minimum size=2mm,draw=black!100,fill=black!30]
    \tikzstyle{whitepoint}=[rectangle,inner sep=0pt,minimum size=2mm,draw=black!100,fill=black!80]
    \tikzstyle{redpoint}=[diamond,inner sep=0pt,minimum size=2.55mm,draw=black!100,fill=black!80]
    \tikzstyle{blackpoint}=[circle,inner sep=0pt,minimum size=2mm,draw=black!100,fill=black!100]
    \tikzstyle{bluepoint}=[circle,inner sep=0pt,minimum size=2mm,draw=black!100,fill=black!0]
    \resizebox{\textwidth}{!}{\begin{tikzpicture}
      \draw [sharp corners] (0,0) node[blackpoint] {} --
      ++(300:1cm) node[blackpoint] {} --
      ++(300:1cm) node[blackpoint] {} --
      ++(300:1cm) node[blackpoint] {} --
      ++(300:1cm) node[redpoint] (L4) {} --
      ++(60:1cm) node[bluepoint] (L3) {} --
      ++(60:1cm) node[bluepoint] (L2) {} --
      ++(60:1cm) node[bluepoint] (L1) {} --
      ++(60:1cm) node[whitepoint] (L0) {} --
      ++(0:1cm) node[greenpoint] {} --
      ++(300:1cm) node[greenpoint] (end) {};
      \draw[dashed] (end) --  ++(300:1cm);
      
      \draw (L0) + (3cm,0) node[blackpoint,label=right:{\small Relax on $Av=0, v \in \mathcal{V}^r$, compute $P$}] {};
      \draw (L1) + (3.5cm,0) node[redpoint,label=right:{\small Compute $v$, s.t., $Av=\lambda Tv$, update $\mathcal{V}^e$}] {};
      \draw (L2) + (4.0cm,0) node[bluepoint,label=right:{\small Relax on $Av = \lambda T v, v \in \mathcal{V}^e$}] {};
      \draw (L3) + (4.5cm,0) node[greenpoint,label=right:{\small Relax on $Av=0, v \in \mathcal{V}^r$ and $Av = \lambda T v, v \in \mathcal{V}^e$, recompute $P$}] {};
      \draw (L4) + (5cm,0) node[whitepoint,label=right:{\small Test
        MG method, update $\mathcal{V}$}] {};

      \draw (L3) +(0,4cm) node (start) {};
      \draw[sharp corners] (start) ++(300:-3cm) node[blackpoint] {} --
      ++(300:1cm) node[blackpoint] {} --
      ++(300:1cm) node[blackpoint] {} --
      ++(300:1cm) node[redpoint] {} --
      ++(60:1cm) node[greenpoint] {} --
      ++(300:1cm) node[redpoint] {} --
      ++(60:1cm) node[bluepoint] {} --
      ++(60:1cm) node[greenpoint] {} --
      ++(300:1cm) node[greenpoint] {} --
      ++(300:1cm) node[redpoint] {} --
      ++(60:1cm) node[greenpoint] {} --
      ++(300:1cm) node[redpoint] {} --
      ++(60:1cm) node[bluepoint] {} --
      ++(60:1cm) node[greenpoint] {} --
      ++(300:1cm) node[greenpoint] {} --
      ++(300:1cm) node[redpoint] {} --
      ++(60:1cm) node[bluepoint] {} --
      ++(60:1cm) node[bluepoint] {} --
      ++(60:1cm) node[whitepoint] {} -- 
      ++(0:1cm) node[greenpoint] {} --     
      ++(300:1cm) node[greenpoint] {} --
      ++(300:1cm) node[greenpoint] {} --
      ++(300:1cm) node[redpoint] {} --
      ++(60:1cm) node[greenpoint] (end) {};
      \draw[dashed] (end) --  ++(300:1cm);
    \end{tikzpicture}}
  \caption{Galerkin Bootstrap AMG W cycle and V cycle setup schemes.\label{fig:boot:setupcycle}}
\end{center}
\end{figure}

We apply two W-cycle bootstrap cycles using four pre- and post-smoothing steps to compute the set of relaxed vectors and set of bootstrap vectors,  
with $|\mathcal{V}^r| = |\mathcal{V}^e|=8$, which together are then used in computing the algebraic distance measure for defining strength of connection and the least squares interpolation operator on each level.   
The compatible relaxation algorithm is applied with search depth $d=2$,  and the LS interpolation is formed
using $d_{LS} = d+2$ as before.  

In the solve phase, the nonlinear AMLI is used as a preconditioner
to the flexible conjugate gradient iteration.  
We note that there is a mild dependence of the iteration counts of the AMLI $W$-cycle preconditioner on the problem size for the non-grid aligned cases, though the coarse grids and interpolation stencils follow the direction of the anisotropy on all levels of the multigrid hierarchy.  This in turn suggests that a higher-caliber ($c>2$) interpolation is needed to obtain scalable multilevel results.    

  \begin{table}[!ht]
{\small    \begin{center}
      \begin{tabular}{|c|c|c|c|c|c|c|c|}
        \hline
        $\alpha$\\
        \hline
        \hline
        $0$\\
        \hline
        $\pi/4$\\
        \hline
        $-\pi/4$ \\ 
        \hline
        $\pi/8$  \\
        \hline \hline
        Levels\\
        \hline
      \end{tabular}
      \begin{tabular}{|c|c|c|c|c|c|c|c|}
        \hline
        \multicolumn{3}{|c|}{FD and $\epsilon =  .1$} \\ 
        \hline      \hline
        4 (1.4,1.8) & 4 (1.5,1.9)  & 4 (1.5,1.9) \\
        \hline 
        4 (1.4,1.5) & 4 (1.5,1.6)  & 4 (1.5,1.6)   \\
        \hline
         5 (1.4,1.9) & 6 (1.5,2.0)  & 7 (1.5,2.0)   \\
        \hline
        4 (1.3,1.5) &  4 (1.3,1.6)   & 5 (1.4,1.7)   \\
        \hline \hline
        3 & 4 & 5\\
        \hline
      \end{tabular}
                     \bigskip
 \begin{tabular}{|c|c|c|c|c|c|c|c|}
        \hline
        \multicolumn{3}{|c|}{FE and $\epsilon =  .1$} \\ 
        \hline      \hline
        5 (1.4,1.4) & 5 (1.4,1.4)  & 6 (1.5,1.4) \\
        \hline 
        5 (1.4,1.5) & 5 (1.5,1.5)  & 5 (1.5,1.6)   \\
        \hline
        5 (1.5,1.6) & 5 (1.5,1.7)  & 6 (1.6,1.7)   \\
        \hline
        6 (1.4,1.3) &  6 (1.4,1.4)   & 7 (1.4,1.4)   \\
        \hline \hline
        3 & 4 & 5\\
        \hline
      \end{tabular}
      \bigskip
      \begin{tabular}{|c|c|c|c|c|c|c|c|}
        \hline
        $\alpha$\\
        \hline
        \hline
        $0$\\
        \hline
        $\pi/4$\\
        \hline
        $-\pi/4$ \\ 
        \hline
        $\pi/8$  \\
        \hline \hline
        Levels\\
        \hline
      \end{tabular}
      \begin{tabular}{|c|c|c|c|c|c|c|c|}
        \hline
        \multicolumn{3}{|c|}{FD and $\epsilon =  .0001$} \\ 
        \hline      \hline
        3 (1.4,1.7) & 4 (1.4,1.9)  & 4 (1.4,2.0) \\
        \hline 
        4 (1.4,1.5) & 4 (1.5,1.6)  & 4 (1.5,1.7)   \\
        \hline
         6 (1.6,2,1) & 8 (1.6,2.2)  & 9 (1.6,2.2)   \\
        \hline
        6 (1.3,1.3) &  6 (1.5,1.9)   & 6 (1.5,2.0)   \\
        \hline \hline
        3 & 4 & 5\\
        \hline
      \end{tabular}
 \begin{tabular}{|c|c|c|c|c|c|c|c|}
        \hline
        \multicolumn{3}{|c|}{FE and $\epsilon =  .0001$} \\ 
        \hline      \hline
        5 (1.4,1.7) & 6 (1.4,1.7)  & 6 (1.4,1.8) \\
        \hline 
        7 (1.4,1.8) & 8 (1.4,.1.9)  & 9 (1.4,2.0)   \\
        \hline
                6 (1.5,1.8) & 7 (1.5,1.9)  & 8 (1.6,1.9)   \\
        \hline
        6 (1.5,1.9) &  6 (1.5,1.9)   & 6 (1.5,2.0)   \\
        \hline \hline
        3 & 4 & 5\\
        \hline
      \end{tabular}
          \end{center}}
      \caption{Nonlinear AMLI W-cycle preconditioned flexible CG applied to the seven-point finite difference anisotropic Laplace problem with Dirichlet boundary conditions for $\epsilon = .1,.0001$, $n = 32^2,64^2,128^2$, and various choices of $\alpha$.  The reported results correspond to the number of iterations needed to reduce the residual by $10^8$ and the grid and operator complexities for the method constructed using a W-cycle bootstrap setup algorithm, again with the same parameters  $d=2$ and $d_{LS} = 4$ that were used for the tests of the two-grid method reported in Table~\ref{tab:2gridFD}.}
  \end{table}    

\section{Concluding remarks}
\label{sec:conclusions}
The LS functional gives a flexible and robust tool for measuring AMG strength of connectivity via algebraic distances. It is computed for  pairs of points to define a strength graph used to choose coarse points and for  sets of points to  determine interpolatory sets.  
The proposed coarsening approach combines algebraic distances, compatible relaxation, and least squares interpolation.  It provides an effective mechanism  for the non-grid aligned anisotropic diffusion problems considered. 
The approach chooses suitable coarse-grid variables and prolongation operators for a wide range of anisotropies, without the need for parameter tuning.  
 Moreover, even when the initial fine-grid discretization is unfavorable, i.e., chosen in the direction opposite to the one defined by the  anisotropy (as in the $\alpha = -\pi/4$ case), the method constructs a suitable interpolation operator and, further, produces a coarse-grid operator which better captures the anistropy directions, correcting the deficiency of the fine grid operator.  Moreover, we have shown using caliber $c=2$ LS interpolation leads to a nearly optimal multilevel method for the targeted constant coefficient anisotropic diffusion problems.  {\color{black}As noted earlier, the main challenge faced in fully extending the approach to an optimal multilevel one for variable coefficient anisotropic problems is that of designing an algorithm capable of constructing long-range interpolation with caliber $c>2$, as needed to accurately capture general anisotropies, that at the same time maintains low grid and operator complexities. This is a research topic that we are currently investigating.  }

\bibliography{bamg_Aniso}

\providecommand{\noopsort}[1]{}
\begin{thebibliography}{10}

\bibitem{B83}
A.~Brandt.
\newblock Algebraic multigrid theory: The symmetric case.
\newblock {\em Appl. Math. Comput.}, 19(1-4):23--56, 1986.

\bibitem{ABrandt_2000a}
A.~Brandt.
\newblock General highly accurate algebraic coarsening.
\newblock {\em Elect. Trans. Numer. Anal.}, 10:1--20, 2000.

\bibitem{B00}
A.~Brandt.
\newblock Multiscale scientific computation: Review 2001.
\newblock In {\em Multiscale and Multiresolution Methods}, pages 1--96.
  Springer Verlag, 2001.

\bibitem{BAMG2010}
A.~Brandt, J.~Brannick, K.~Kahl, and I.~Livshits.
\newblock Bootstrap {AMG}.
\newblock {\em SIAM J. Sci. Comput.}, 33:612--632, 2011.

\bibitem{oAMG}
A.~Brandt, S.~McCormick, and J.~Ruge.
\newblock Algebraic multigrid ({AMG}) for sparse matrix equations.
\newblock In D.~J. Evans, editor, {\em Sparsity and Its Applications}.
  Cambridge University Press, Cambridge, 1984.

\bibitem{BMR83}
A.~Brandt, S.~McCormick, and J.~W. Ruge.
\newblock Algebraic multigrid ({AMG}) for automatic multigrid solution with
  application to geodetic computations.
\newblock Technical report, Colorado State University, Fort Collins, Colorado,
  1983.

\bibitem{JBrannick_2005a}
J.~Brannick.
\newblock {\em Adaptive algebraic Multigrid coarsening strategies}.
\newblock PhD thesis, University of Colorado at Boulder, 2005.

\bibitem{Bran_EBSOC}
J.~Brannick, M.~Brezina, S.~Maclachlan, T.~Manteuffel, S.~Mccormick, and
  J.~Ruge.
\newblock An energy-based {AMG} coarsening strategy.
\newblock {\em Numer. Linear Algebra Appl.}, 13:133--148, 2006.

\bibitem{brannick_amli_2011}
J.~Brannick, Y.~Chen, J.~Krauss, and L.~Zikatanov.
\newblock An algebraic multigrid method based on matching of graphs.
\newblock In {\em Proceedings of the 20th International Conference on Domain
  Decomposition Methods}, Lecture Notes in Computational Science and
  Engineering. Springer, Submitted May 16, 2011.

\bibitem{brannick_local_stab_2011}
J.~Brannick, Y.~Chen, and L.~Zikatanov.
\newblock An algebraic multilevel method for anisotropic elliptic equations
  based on subgraph matching.
\newblock {\em Numer. Linear Algebra Appl.}, To appear, accepted for
  publication November 17, 2011.

\bibitem{JBrannick_RFalgout}
J.~Brannick and R.~Falgout.
\newblock Compatible relaxation and coarsening in algebraic multigrid.
\newblock {\em SIAM J. Sci. Comput.}, 32(3):1393--1416, 2010.

\bibitem{Brannick_Trace_06}
J.~Brannick and L.Zikatanov.
\newblock Algebraic multigrid methods based on compatible relaxation and energy
  minimization.
\newblock In {\em Proceedings of the 16th International Conference on Domain
  Decomposition Methods}, volume~55 of {\em Lecture Notes in Computational
  Science and Engineering}, pages 15--26. Springer, 2007.

\bibitem{MBrezina_etal_2000a}
M.~Brezina, A.J. Cleary, R.D. Falgout, V.E. Henson, J.E. Jones, T.A.
  Manteuffel, S.F. McCormick, and J.W. Ruge.
\newblock Algebraic multigrid based on element interpolation ({AMG}e).
\newblock {\em SIAM J. Sci. Comp.}, 22:1570--1592, 2000.

\bibitem{MBrezina_etal_2005a}
M.~Brezina, R.~Falgout, S.~MacLachlan, T.~Manteuffel, S.~McCormick, and
  J.~Ruge.
\newblock Adaptive amg (a{AMG}).
\newblock {\em SIAM J. Sci. Comput.}, 26:1261--1286, 2005.

\bibitem{WLBriggs_VEHenson_SFMcCormick_2000a}
W.~L. Briggs, V.~E. Henson, and S.~F. McCormick.
\newblock {\em A Multigrid Tutorial}.
\newblock SIAM Books, Philadelphia, 2000.
\newblock Second edition.

\bibitem{TChartier_etal_2003a}
T.~Chartier, R.D. Falgout, V.E. Henson, J.E. Jones, T.A. Manteuffel, S.F.
  McCormick, J.W. Ruge, and P.S. Vassilevski.
\newblock Spectral {AMG}e ($\rho${AMG}e).
\newblock {\em SIAM J. Sci. Comp.}, pages 1--26, 2003.

\bibitem{AJCleary_RDFalgout_VEHenson_JEJones_1998a}
A.J. Cleary, R.D. Falgout, V.E. Henson, and J.E. Jones.
\newblock Coarse-grid selection for parallel algebraic multigrid.
\newblock In {\em Proc. of the Fifth International Symposium on Solving
  Irregularly Structured Problems in Parallel}, volume 1457 of {\em Lecture
  Notes in Computer Science}, pages 104--115. Springer--Verlag, New York, 1998.

\bibitem{HDeSterck_UMYang_JJHeys_2005a}
H.~{De Sterck}, U.M. Yang, and J.J. Heys.
\newblock Reducing complexity in parallel algebraic multigrid preconditioners.
\newblock {\em SIAM J. on Matrix Analysis and Applications}, 27(4):1019--1039,
  2006.
\newblock Also available as {LLNL} Technical Report {UCRL}-{JRNL}-206780.

\bibitem{PanayotRob_2003}
R.~Falgout and P.~Vassilevski.
\newblock On generalizing the amg framework.
\newblock {\em SIAM J. Numer. Anal.}, 42(4):1669--1693, 2004.

\bibitem{VEHenson_UMYang_2002a}
V.E. Henson and U.M. Yang.
\newblock Boomer{AMG}: a parallel algebraic multigrid solver and
  preconditioner.
\newblock {\em Applied Numerical Mathematics}, 41:155--177, 2002.

\bibitem{OLivne_2004a}
O.~Livne.
\newblock Coarsening by compatible relaxation.
\newblock {\em Num. Lin. Alg. Appl.}, 11(2):205--227, 2004.

\bibitem{iBAMG}
T.~Manteuffel, S.~McCormick, M.~Park, and J.~Ruge.
\newblock Operator-based interpolation for bootstrap algebraic multigrid.
\newblock {\em Numer. Linear Algebra Appl.}, 17(2-3):519--537, 2010.

\bibitem{Olson_10}
L.~Olson, J.~Schroder, and R.~Tuminaro.
\newblock A new perspective on strength measures in algebraic multigrid.
\newblock {\em Numerical Linear Algebra with Appl.}, 17:713--733, 2010.

\bibitem{DB_11}
D.~Ron, I.~Safro, and A.~Brandt.
\newblock Relaxation based coarsening and multiscale graph organization.
\newblock {\em SIAM J. on Multi. Mod. and Sim.}, 9:407--423, 2011.

\bibitem{JWRuge_KStuben_1987a}
J.~W. Ruge and K.~St{\"u}ben.
\newblock Algebraic multigrid ({AMG}).
\newblock In S.~F. McCormick, editor, {\em Multigrid Methods}, volume~3 of {\em
  Frontiers in Applied Mathematics}, pages 73--130. SIAM, Philadelphia, PA,
  1987.

\bibitem{Jacob}
J.~B. Schroder.
\newblock Smoothed aggregation solvers for anisotropic diffusion.
\newblock {\em Numer. Linear Algebra Appl.}, 19:296--312, 2012.

\bibitem{PVanek_JMandel_MBrezina_1995a}
P.~Van\v{e}k, J.~Mandel, and M.~Brezina.
\newblock ~{A}lgebraic multigrid by smoothed aggregation for second and fourth
  order elliptic problems.
\newblock {\em Computing}, 56:179--196, 1996.

\bibitem{Vassilevski_2005}
P.~Vassilevski.
\newblock {\em Multilevel Block Factorization Preconditioners: Matrix-based
  Analysis and Algorithms for Solving Finite Element Equations}.
\newblock Springer, 2009.

\bibitem{JWilkinson_1965}
J.~Wilkinson.
\newblock {\em The Algebraic eigenvalue problem}.
\newblock Clarendon Press, Oxford, 1965.

\end{thebibliography}
\bibliographystyle{plain}
\end{document}